\documentclass[a4paper,12pt]{article}
\usepackage{amssymb,amsmath,amsthm,dsfont,mathtools,stmaryrd}

\usepackage{pdfpages}

\usepackage{newclude}

\usepackage[utf8]{inputenc}
\usepackage[english]{babel}

\usepackage{lmodern}
\usepackage{csquotes} 
\usepackage{color} 

\usepackage{geometry}
\usepackage{soul}
\setstcolor{red}
\usepackage{graphicx}
\graphicspath{{Images/}}

\usepackage[export]{adjustbox}

\usepackage{float}

\usepackage{subcaption}

\usepackage[
backend=biber,
style=apa,
sorting=nyt,
natbib=true
]{biblatex}
\usepackage{proba}

\mathtoolsset{showonlyrefs}

\usepackage{enumitem}

\usepackage{algorithm}
\usepackage{algpseudocode}

\usepackage{hyperref}
\hypersetup{colorlinks,
breaklinks,
urlcolor=blue,
linkcolor=black,
citecolor=blue}

\usepackage{arydshln}

\newcommand{\PL}{\mathcal{PL}}

\renewcommand{\N}{\mathbb{N}}
\renewcommand{\Z}{\mathbb{Z}}
\renewcommand{\R}{\mathbb{R}}
\renewcommand{\X}{\mathcal{X}}

\renewcommand{\d}{\mathrm{d}}
\renewcommand{\c}{\mbox{Cov}}
\newcommand{\e}{\mathrm{e}}

\newcommand{\ind}{\mathds{1}}

\newcommand{\vertii}[1]{{\left\vert\kern-0.25ex\left\vert #1 
    \right\vert\kern-0.25ex\right\vert}}
\newcommand{\vertiii}[1]{{\left\vert\kern-0.25ex\left\vert\kern-0.25ex\left\vert #1 
    \right\vert\kern-0.25ex\right\vert\kern-0.25ex\right\vert}}

\newcommand{\argmax}{\mathrm{argmax}}

\usepackage{relsize}
\usepackage{exscale}

\makeatletter
\let\amsmath@bigm\bigm

\renewcommand{\bigm}[1]{%
  \ifcsname fenced@\string#1\endcsname
    \expandafter\@firstoftwo
  \else
    \expandafter\@secondoftwo
  \fi
  {\expandafter\amsmath@bigm\csname fenced@\string#1\endcsname}%
  {\amsmath@bigm#1}%
}

\newcommand{\DeclareFence}[2]{\@namedef{fenced@\string#1}{#2}}
\makeatother

\DeclareFence{\mid}{|}

\newtheorem{lem}{Lemma}[section]

\newtheorem{prop}[lem]{Proposition}

\newtheorem{theo}[lem]{Theorem}

\newtheorem{cor}[lem]{Corollary}

\theoremstyle{definition}

\newtheorem*{note}{Note}

\newenvironment{lemproof}[1][]{\emph{Proof#1.}}{\hfill $\square$ \\}

\usepackage[toc,page,titletoc]{appendix}

\providecommand{\keywords}[1]
{
  \small	
  \textbf{\textit{Keywords---}} #1
}

\usepackage{multicol}

\title{Theory and inference for multivariate autoregressive binary models with an application to absence-presence data in ecology}
{\small
\author{Guillaume Franchi {\it guillaume.franchi@ensai.fr} \\
           Univ Rennes, Ensai, CNRS, CREST -- UMR 9194, F-35000 Rennes, France\\
        Lionel Truquet {\it lionel.truquet@ensai.fr} \\
             Univ Rennes, Ensai, CNRS, CREST -- UMR 9194, F-35000 Rennes, France
             }}
\date{}

\begin{document}

\maketitle

\begin{abstract}
We introduce a general class of autoregressive models for studying the dynamic of multivariate binary time series with stationary exogenous covariates. Using a high-level set of assumptions, we show that existence of a stationary path for such models is almost automatic and does not require parameter restrictions when the noise term is not compactly supported. We then study in details statistical inference in a dynamic version of a multivariate probit type model, as a particular case of our general construction. To avoid a complex likelihood optimization, we combine pseudo-likelihood and pairwise likelihood methods for which asymptotic results are obtained for a single path analysis and also for panel data, using ergodic theorems for multi-indexed partial sums.  The latter scenario is particularly important for analyzing absence-presence of species in Ecology, a field where data are often collected from surveys at various locations. Our results also give a theoretical background for such models which are often used by the practitioners but without a probabilistic framework.
\end{abstract}

\medskip


\section{Introduction}

As a benchmark model for studying several correlated binary responses, multivariate logit or probit regression models have been extensively used in various domains. For instance, they have been applied to labor supply in economics \citep{chib1998analysis}, to the modeling of multiple claims in insurance \citep{young2009multivariate}, for evaluating customer satisfaction \citep{kenne2016pairwise}, for studying the response of biological organisms exposed to stimulus \citep{ashford1970multi} or for studying absence-presence of several species in an ecosystem \citep{wilkinson2019comparison}, \citep{poggiato2021interpretations}. \citet{sebastian2010testing} also applied a multivariate logistic model for analyzing such absence-presence time series in several ecosystems and \citet{guanche2014autoregressive} studied atmospheric circulation with similar models.
 
Though a few contributions considered dynamic multivariate probit models for time series analysis of binary vectors (see for instance \citet{candelon2013multivariate}, \citet{chaubert2008multivariate}, \citet{manner2016modeling}), theoretical contributions to statistical inference in such models are rather scarce.
A theory for the univariate dynamic probit model can be found in \citet{de2011dynamic} while \citet{tuzcuoglu2023composite} investigated the AR-probit model introduced by \citet{muller2005autoregressive} in the case of panel data.

Let us also mention another type of multiple binary autoregressive time series models based on random binary coefficients, an approach developed recently by \citet{jentsch2022generalized} who introduced a rigorous framework for defining stationary solutions and estimating model parameters. This class of model seems to be quite flexible but a possible inclusion of exogenous covariates in the corresponding model equations is unclear because the link function cannot map continuous data to multiple binary vectors. 
\medskip

In this paper, we consider the following model. For $i \in \left\lbrace 1, \ldots , k \right\rbrace$, let $\left(Y_{i,t}\right)_{t\in \Z}$ be a univariate  binary time series. Following our initial motivation for applying these models to Ecology, we will often refer to a specific application in this context, where $Y_{i,t}$ simply measures the absence-presence of a specific species $i$ at time $t$ in a given area. This analogy will be often useful to interpret model parameters, though our generic results can be used for analyzing arbitrary multiple binary time series.  Setting $Y_t=\left(Y_{1,t},\ldots,Y_{k,t}\right)\in \{0,1\}^k$, we assume that the dynamic of the multivariate time series $\left(Y_t\right)_{t \in \Z}$ is given by
	\begin{equation} 
	Y_{i,t}=\ind_{(0,+\infty)} \left( g_i(Y_{t-1},\ldots ,Y_{t-p},\zeta_t) \right), \quad 1 \leqslant i \leqslant k, \label{eq_model}
	\end{equation}
where $g =(g_1,\ldots ,g_k)':\{0,1\}^{kp}\times \R^n\rightarrow \R^k$ is a measurable mapping and $(\zeta_t)_{t\in\Z}$ is a stationary stochastic process taking values in $\R^n$.

Hence the multivariate binary vector $Y_t$ depends on its values at previous times and can be interpreted as the discretization of a multivariate continuous response $Y^{*}_t=g(Y_{t-1},\ldots , Y_{t-p},\zeta_t)$. Allowing dependence between the variables $\left( \zeta_t\right)_{t \in \Z}$ is interesting for including exogenous covariates, i.e. $\zeta_t=\left(X_{t-1},\varepsilon_t\right)$, where $X_{t-1}$ is a vector of covariates determined at previous time (and then observable before time $t$) and $\varepsilon_t$ is an unobserved noise component (that in some cases will be assumed to be independent from $X_{t-1}$). 
When $g$ is linear, i.e.
	\begin{equation} 
	g\left(Y_{t-1},\ldots,Y_{t-p},\zeta_t\right)=\sum_{j=1}^p A_j Y_{t-j}+B X_{t-1}+\varepsilon_t, \label{eq_model_lin}
	\end{equation}
where $A_1,\ldots,A_p$ are square matrices of size $k\times k$, $B$ is a matrix of size $k\times d$, $\varepsilon_t$ is independent of $(Y_{t-1},\ldots,Y_{t-p},X_{t-1})$ 
and follows a Gaussian distribution with mean $0$ and covariance matrix $R$, we obtain a dynamic version of the multivariate probit regression model.

\medskip

More generally, model \eqref{eq_model} is a natural time series version of the multivariate discrete choice models popular in econometrics (see for instance \citet{greene2009discrete}). For $k=1$, it extends the dynamic binary choice model of \citet{de2011dynamic}. When $F$ is the CDF of the logistic distribution, we get the multivariate logistic model used for instance by \citet{sebastian2010testing} for studying absence-presence data. Note that our version has some natural interpretations in term of dynamics for a multivariate binary response. For instance, when $p=1$ for simplicity, $A_1(i,j)$ is positive if and only if the presence of a species $j$ at time $t-1$ influences positively the presence of the species $i$ at time $t$.
Furthermore, the sign of $R_{i,j}$ indicates if positive or negative dependence occurs at time $t$ (positivity or negativity of the conditional covariances between $Y_{i,t}$ and $Y_{j,t}$ meaning that the values $1$ are generally produced simultaneously or not).

\medskip
This paper provides a threefold contribution.

First,  we give simple and very general conditions guarantying the existence of a strictly stationary solution for \eqref{eq_model}. This point, which is not straightforward because of the non-linearity of the link function, is always an important step in time series analysis since it is very useful to derive asymptotic properties of statistical inference procedures. As we will see, stationarity is almost automatic for model \eqref{eq_model_lin}, without any restriction on the coefficients. Note that from the presence of exogenous covariates, the model is not Markovian in general, leading to an additional difficulty for studying stationarity properties. 

Secondly, we investigate statistical inference for the autoregressive parameters $A_1,\ldots,A_p,B$ and the correlation matrix $R$ in model \eqref{eq_model_lin}, which is quite challenging. Likelihood estimators are generally intractable due to the presence of multiple integrals involving the Gaussian density in the expression of the likelihood function. There exists many approaches developed for inference in multivariate probit models. Some of them are based on a Bayesian framework \citep{chib1998analysis}, \citep{talhouk2012efficient}, or on importance sampling such as the GHK algorithm \citep{hajivassiliou1994classical}, or on composite likelihood methods \citep{zhao2005composite}, \citep{ting2022fast}. In the time series context, we will combine univariate and pairwise composite likelihood methods to estimate all the parameters of the model. The advantage of the proposed method is its fast implementation with respect to the methods based on the true likelihood function. We are not aware of a theoretical analysis of some inference procedures for multivariate logistic/probit models such as the double-step approach mentioned above.

Finally, the case of data collected at different sites is even more challenging due to a double asymptotic when both the number of sites $n$ and the horizon time $T$ go to infinity. The benefit of this work is to provide very simple assumptions for getting asymptotic results in this context. In particular, we use specific versions of the law of large numbers for partial sums with a double summation index. All the results are formulated using realistic assumptions on the covariates which can contain individual effects, time effects or both. Let us also mention that our results, which are based on a multi-parameter ergodic theorem not very well known in the time series literature, can be also interesting for studying more general panel type models.

The paper is organized as follows. In Section \ref{2}, we address the problem of stationarity for these dynamics.
Inference methods for the multivariate autoregressive probit model are introduced and discussed for a single path in Section \ref{sec_param_inf_single}, while Section \ref{4} is devoted to the problem of multiple paths with panel data. The results are illustrated in Section \ref{5} with an application to the dynamics of fisheries data. The proofs of all our results are given in a supplementary material available online.

\section{Existence of stationary solutions}\label{2}

\subsection{A general result for multiple binary autoregressive model}

In this section, we first state a general result ensuring the existence of a stationary solution for model \eqref{eq_model}. The assumptions used below will be much more general than what is needed for studying the multivariate linear model \eqref{eq_model_lin}. In what follows, we denote by $\mathcal{F}_t$ the sigma-field generated by the random variables $\zeta_s$, $s\leqslant t$:
	\begin{equation}
	\mathcal{F}_t = \sigma \left( \zeta_s \mid s \leqslant t \right).
	\end{equation}
We consider the following set of assumptions.

	\begin{enumerate}[label=\textbf{A\arabic*}]
	\item The process $\left(\zeta_t\right)_{t\in\Z}$, taking values in $\R^n$, is strictly stationary.
	\item For $i=1,\ldots,k$, let us introduce the two events	
		\begin{equation}
		C_i^{0}=\left\lbrace \underset{1\leq t\leq p}{\max} \ \underset{y_1,\ldots,y_p\in \left\lbrace 0,1\right\rbrace^k}{\max} g_i\left(y_1,\ldots,y_p,\zeta_t\right)\leqslant 0\right\rbrace,
		\end{equation}		
	and	
		\begin{equation}
		C_i^{1}=\left\lbrace \underset{1\leq t\leq p}{\min} \ \underset{y_1,\ldots,y_p\in \left\lbrace 0,1\right\rbrace^k}{\min} g_i\left(y_1,\ldots,y_p,\zeta_t\right)>0\right\rbrace.
		\end{equation}		
We assume that there exist $j_1,\ldots,j_k$ in $\left\lbrace 0,1 \right\rbrace$ such that $\P\left(\cap_{i=1}^kC_i^{j_i}\right)>0$.
	\item There exist $j_1,\ldots,j_k$ in $\{0,1\}$ such that $\P\left(\cap_{i=1}^kC_i^{j_i}\vert \mathcal{F}_0\right)>0$ almost surely.
	\end{enumerate}
	
Let us first comment these assumptions. Assumption \textbf{A1} is the minimal assumption needed for working with stationary time series. The technical assumption \textbf{A2} has an important interpretation. Assume that the stochastic recursions in \eqref{eq_model} are initialized with some values $Y_t=y_t$ for $-p+1\leqslant t\leqslant 0$. If the event $\cap_{i=1}^k C_i^{j_i}$ occurs, whatever the initial values $y_{-p+1},\ldots,y_0$, we will obtain the value $Y_{i,t}=j_i$ for $1\leqslant t\leqslant p$ and $1\leqslant i\leqslant k$. As a consequence, the values of $Y_t$ for $t\geqslant p+1$ will no more depend on the initial values. This kind of condition is analog to the standard irreducibility and aperiodicity conditions imposed to finite-state Markov chains though the model studied here does not have the Markov property.

From these two assumptions, it is possible to prove a stochastic convergence property for the iterations of the recursions \eqref{eq_model}. The proof of the following result is based on a general coalescence argument for iterated random maps acting on finite sets, an idea already used by \citet{debaly2021iterations} for the binary choice model and extended by \citet{truquet2023strong} to more general model, in the case of ergodic inputs. This argument has some flavors with a similar one used by \citet{propp1996exact} to simulate perfectly a realization of the invariant probability of a finite-state Markov chain. However, our aim here is also to discuss the existence of stationary but not necessarily ergodic paths (as it will be required for studying panel data with individual effects). We have to adapt the techniques used in these references, replacing Assumption {\bf A2} by a more stringent condition given in Assumption \textbf{A3}.

	\begin{theo} \label{theo_stationarity} \hfill
		\begin{enumerate}[label=\arabic*)]
		\item Suppose that assumptions \emph{\textbf{A1-A2}} hold true and the process $\left(\zeta_t\right)_{t\in\Z}$ is ergodic. Then there exists a unique strictly stationary solution $(Y_t)_{t\in \Z}$ for the recursions \eqref{eq_model}.
		Moreover, the unique solution has a Bernoulli shift representation, i.e. there exists a measurable mapping $H:\left(\R^n\right)^{\N}\rightarrow \left\lbrace 0,1\right\rbrace$ such that	
			\begin{equation}\label{eq_Bernoulli_shift}
			Y_t=H\left(\zeta_t,\zeta_{t-1},\ldots\right),\quad t\in\Z.
			\end{equation}		
		The process $\left((Y_t,\zeta_t)\right)_{t\in\Z}$ is then ergodic.
		\item Suppose now that assumptions \emph{\textbf{A1-A3}} hold true with $(\zeta_t)_{t\in\Z}$ not necessarily ergodic. There still exists a unique strictly stationary solution $(Y_t)_{t\in\Z}$ for the recursions \eqref{eq_model} and the Bernoulli shift representation is still valid. 
		\end{enumerate}
	\end{theo}
	
	\paragraph{Notes.}
	
		\begin{enumerate}[label=\arabic*.]
		\item When $\zeta_t=(X_{t-1},\varepsilon_t)$ and $X_t$ is a vector of covariates observed at time $t$ and $\varepsilon_t$ is an unobserved noise component, we point out that the non-ergodic case is important in some applications. Indeed, in panel data, some of the components of the vector $X_t$ may be not time-varying. For instance, $X_t=(Z,\Gamma_t)$ with $Z$ being a random vector independent of a stationary process $\left((\Gamma_t,\varepsilon_t)\right)_{t\in\Z}$. In this case, the process $(\zeta_t)_{t\in\Z}$ is still stationary but obviously not ergodic. Let us also mention that stationarity of the process $(\zeta_t)_t$ is equivalent to stationarity for the process $(X_t,\varepsilon_t)_t$ and a similar remark is valid for ergodicity.
		\item In the ergodic case, Theorem \ref{theo_stationarity}, point \emph{1)} is useful for applying Birkhoff's ergodic theorem to some partial sums of the form $\frac{1}{T}\sum_{t=1}^T f\left(Y_t,X_t,Y_{t-1},X_{t-1},\ldots\right)$ where $f$ is a measurable function such that $\E\left[f\left(Y_0,X_0,Y_{-1},X_{-1},\ldots\right)\right]<+\infty$. 
This ergodic theorem is generally sufficient to derive consistency and asymptotic normality of conditional likelihood or pseudo-likelihood estimators in smooth parametric time series models.  
\end{enumerate}	
	
\subsection{Specific results for models with covariates and an additive noise component}

In this part, we assume that $\zeta_t=\left(X_{t-1},\varepsilon_t\right)$ takes values in $\R^d\times \R^k$ (then $n=d+k$) and we consider the more specific dynamic
	\begin{equation}\label{intereq}
	g\left(Y_{t-1},\ldots,Y_{t-p},\zeta_t\right)=h\left(Y_{t-1},\ldots,Y_{t-p},X_{t-1}\right)+\varepsilon_t,
	\end{equation}	
where $h:\left( \left\lbrace 0,1\right\rbrace^k\right)^p\times \R^d\rightarrow \R^k$ is a measurable mapping. Of course, one could also assume a dependence with respect to $X_{t-2},\ldots, X_{t-p}$ but in this case it is possible to replace $X_{t-1}$ in \eqref{intereq} by the vector $Z_{t-1}:=(X_{t-1},\ldots,X_{t-p})$ so that model \eqref{intereq} contains implicitly such an extension.

\medskip

The following condition is sufficient for getting \textbf{A3} and then \textbf{A2}.

	\begin{enumerate}
\item[\textbf{A4}] Almost surely, the conditional distribution of $\varepsilon_t$ given $\sigma
\left( \varepsilon_s , X_s \mid s \leqslant t-1  \right)$ has a full support $\R^k$. 
	\end{enumerate}
	
\medskip

Assumption \textbf{A4} is sufficient but not necessary for checking \textbf{A3}, other more sophisticated conditions on the support of the previous conditional distribution are possible, but \textbf{A4} is simpler to formulate. From Assumption \textbf{A4}, we get a simple and meaningful condition for recurrence of the states in the dynamic. Indeed, whatever the previous values of the process and the covariates are, the process is allowed to move to an arbitrary element of $\left\lbrace 0, 1 \right\rbrace^k$ (i.e. we have positive conditional probabilities).

Note that Assumption \textbf{A4} does not require a specific form for the mapping $h$. Moreover, it is valid when $\varepsilon_t$ is independent from $\sigma
\left( \varepsilon_s , X_s \mid s \leqslant t-1  \right)$ (i.e. the covariates are sequentially exogenous) and its probability distribution has a full support (as for the multivariate Gaussian for instance). But assumption \textbf{A4} also holds true with non i.i.d. noises. For instance, it is satisfied for infinite moving averages $\varepsilon_t=\eta_t+\sum_{j=1}^{\infty}A_j \eta_{t-j}$ where $(A_j)_{j\geqslant 0}$ is a summable sequence of square matrices and $\left((X_t,\eta_t)\right)_{t\in\Z}$ is a stationary process such that $\eta_t$ has a probability distribution with full support and is independent from $\sigma\left((\eta_s,X_s) \mid s\leqslant t-1\right)$. Let us recall that from the definition of the support of a probability measure, i.e. the intersection of all the closed sets with probability $1$, the full support assumption entails that $\P\left(\eta_t\in O\right)>0$ for any non-empty open subset $O$.

	\begin{cor} \label{cor_stationarity}
	Suppose that Assumptions \emph{\textbf{A1-A4}} hold true.
	
	Then the conclusions of Theorem \ref{theo_stationarity} are valid. 
	\end{cor}

 \paragraph{Note.} When there is no covariates (i.e. $h$ does not depend on $X_{t-1}$) and the $\varepsilon_t$'s are i.i.d., the process is simply an irreducible $p-$order Markov chain and we already know that there exists a unique strictly stationary solution to \eqref{eq_model} and \eqref{eq_model_lin}. Although the process is not necessarily Markovian here due to the presence of the covariates, one can observe that the stationarity conditions are similar and do not impose any specific condition on $h$.
 Things become more complicated if the noise $\varepsilon_t$ has a bounded support even when no covariates are included in the dynamics. For instance, suppose that $p=k=1$ with $\varepsilon_t$ supported on $[0,1]$. When $h(y)=2y-1$, the chain is not irreducible and points $0$ and $1$ are absorbing states. On the other hand, when $h(y)=-y$, the chain is irreducible but periodic as it alternates between state $0$ and state $1$.
 One can then see on this simple example that the behavior of the chain strongly depends on the parameters of the model. 

\section{Theory and inference for a multivariate probit autoregressive model} \label{sec_param_inf_single}

\subsection{Framework}

In this section, we consider a stochastic process $(Y_t)_{t\in \Z}$ satisfying equations \eqref{eq_model} and \eqref{eq_model_lin} with multivariate Gaussian errors i.e.
	\begin{equation} \label{eq_model_single}
	\forall t \in \Z, \ \forall i \in \left\lbrace 1,\ldots , k\right\rbrace, \ Y_{i,t} = \ind_{(0,+\infty)} \left( \lambda_{i,t} + \varepsilon_{i,t} \right),
	\end{equation}
where
	\begin{equation} \label{eq_lambda_single}
	\lambda_t = C+\sum_{l=1}^p A_l Y_{t-l} + B X_{t-1}.
	\end{equation}
The sequences $(X_t)_{t\in \Z}$ and $(\varepsilon_t)_{t\in \Z}$ of random vectors are valued in $\R^d$ and $\R^k$ respectively, $A_1,\ldots, A_p$ are matrices of size $k\times k$, $B$ is a matrix of size $k\times d$ and $C$ is a vector of $\R^k$. The process $(X_t)_{t\in \Z}$ is actually a process of covariates, while $(\varepsilon_t)_{t\in \Z}$ corresponds to an unobserved noise component.

Finally, we now denote by $\mathcal{F}_t$ the sigma-field generated by $X_s$ and $\varepsilon_s$ for $s \leqslant t$, i.e.
	\begin{equation}
	\mathcal{F}_t = \sigma \left( X_s,\varepsilon_s \mid s \leqslant t \right).
	\end{equation}	
The following assumptions will be used.

	\begin{enumerate}[label=\textbf{B\arabic*}]
	\item The process $(\zeta_t)_{t\in \Z}$ is strictly stationary (resp. strictly stationary and ergodic).
	\item For any $t \in \Z$, the random vector $\varepsilon_t$ is independent from $\mathcal{F}_{t-1}$.
	\item The random vectors $\varepsilon_t$ are normally distributed
		\begin{equation}
		\varepsilon_t \sim \mathcal{N}(0,R)
 		\end{equation}
 	where $R$ is a correlation matrix.
  \end{enumerate}

The following result is a straightforward consequence of Corollary \ref{cor_stationarity}.

\begin{cor}
Suppose that Assumptions {\bf B1-B2-B3} hold true. Then there exists a unique stationary (resp. stationary and ergodic) solution $(Y_t))_{t\in\Z}$ to the recursions (\ref{eq_model_single}) and (\ref{eq_lambda_single}).
\end{cor}

\paragraph{Notes.}
\begin{enumerate}
\item Assumptions \textbf{B1} and \textbf{B2} are verified in the particular case where there exists a measurable application $H$ such that $X_t = H(\eta_t,\eta_{t-1},\ldots)$, with the variables $(\eta_t,\varepsilon_t)_{t\in \Z}$ i.i.d., but $\eta_t$ not necessarily independent from $\varepsilon_t$. For example, $(X_t)_{t\in \Z}$ can be a $\mathrm{MA}(\infty)$ process. This assumption is sufficiently general and does not require a specific dynamic for the covariates such as a Markov property for instance.
		
\item 
When $A_1=\cdots=A_p=0$, the previous model reduces to the standard multivariate probit regression model studied for instance in \citet{chib1998analysis}. We also recover a specific case of the dynamic model studied in \citet{candelon2013multivariate} but without theory. Note that our version is particularly convenient for model interpretation. For instance, in the case of absence/presence data, a positive (resp. negative) sign for $A_j(i,\ell)$ indicates if the presence of an element $\ell$ at time $t-j$ promotes or not the presence of the element $i$ at time $t$. An interpretation of the elements of the correlation matrix $R$ can be also given in term of 
conditional dependencies. For instance, conditionally on $\mathcal{F}_{t-1}$, the conditional covariance between $Y_{i,t}$ and $Y_{\ell,t}$ have the same sign as $R(i,\ell)\in [-1,1]$. Indeed if $\Phi$ denotes the cdf of the standard Gaussian distribution and if $\Phi_{R(i,\ell)}$ denotes the cdf of a Gaussian vector with standard Gaussian marginals and correlation coefficient $R(i,\ell)$, we have
$$\c\left(Y_{i,t},Y_{\ell,t}\vert \mathcal{F}_{t-1}\right)=\Phi_{R(i,\ell)}\left(\lambda_{i,t},\lambda_{\ell,t}\right)-\Phi\left(\lambda_{i,t}\right)\Phi\left(\lambda_{\ell,t}\right).$$
From Slepian's lemma, this quantity increases with $R(i,\ell)$ and equals to $0$ when $R(i,\ell)=0$.
\item
Another very nice property of this model is its invariance if we are only interested in predicting or explaining a specific subset of the $k$ binary random variables. Indeed, if equation (\ref{eq_lambda_single}) is valid for the $k$ coordinates, then for any subset $I$ of $\left\{1,\ldots,k\right\}$, the time series $\left((Y_{i,t})_{i\in I}\right)_{t\in\Z}$ follows the same kind of dynamic by considering the $Y_{i,t-j}'$s for $i\notin I$ as exogenous regressors as we do for $X_{t-1}$. Including exogenous regressors satisfying {\bf B2}, without assuming complete independence with the entire noise sequence is then very important for such models.
\item 
Considering a correlation matrix $R$ instead of a covariance matrix for the noise term $\varepsilon_t$ is mainly motivated by identification issues.
One can always multiply $\lambda_t$ and $\varepsilon_t$ by a diagonal matrix with positive entries and this will not change the sign of the coordinates. Up to a change of parameters, one can always assume the covariance matrix of $\varepsilon_t$ to be a correlation matrix.
\end{enumerate}

\subsection{Statistical inference from a single path}\label{sec_asym_res_single}
In this part, we assume that the random vectors $Y_1,\ldots,Y_T$ and $X_1,\ldots,X_{T-1}$ are observed for some positive integer $T$. Our aim is to estimate regression parameters $C,A_1,\ldots,A_p,B$ but also 
the correlation matrix $R$ of the noise component $\varepsilon_t$. 
Likelihood inference in multivariate probit model is particularly tricky due to the presence of multiple integrals of the Gaussian density over some rectangles with boundaries depending on the unknown parameters. Bayesian inference combined with MCMC methods have been developed in \citet{chib1998analysis} and there also exist importance sampling strategies such as in \citet{liesenfeld2010efficient} or \citet{Varin}. In the dynamic setting, the conditional likelihood is based on the conditional distribution given below.
We first introduce some notations. We set
	\begin{equation}
	I_1=(0,+\infty) \quad \text{and} \quad I_0=(-\infty,0].
	\end{equation}
The vector of all unknown parameters is given by
	\begin{equation}
	\theta=\left(\gamma,r'\right)
	\end{equation}
where $\gamma$ is a column vector containing the column vectors of the matrices $C, A_1,\ldots,A_p,B$ and
	\begin{equation}
	r = \left( R(1,2),\ldots,R(1,k),R(2,3),\ldots,R(2,k),\ldots,R(k-1,k) \right).
	\end{equation}
The vector $\lambda_t$ is actually a function of $\gamma$, and without ambiguity, we will rather denote $\lambda_t$ by $\lambda_t(\gamma)$ and $\lambda_{i,t}$ by $\lambda_{i,t}(\gamma)$ for $1 \leqslant i \leqslant k$.

In what follows, we focus on two estimation procedures that can be used for the inference of the true parameter $\theta_0=(\gamma_0,r_0')$ satisfying \eqref{eq_model_single} and \eqref{eq_lambda_single}.
For $t \in \Z$, the distribution of $Y_t$ conditionally to $\mathcal{F}_{t-1}$ is given by
		\begin{equation}
		\mathcal{L} \left( Y_t \mid \mathcal{F}_{t-1}\right) = \sum_{s \in \left\lbrace 0,1 \right\rbrace^k } \left( \int_{\R^k} \prod_{i=1}^k \ind_{I_{s(i)}} \left( \lambda_{i,t}(\gamma_0) +x_i \right) \varphi_{R_0} (x) \d x \right) \delta_s,
		\end{equation}
	where $\varphi_R$ denotes the Gaussian density on $\R^k$ with mean 0 and covariance matrix $R$. When $k=1$, we simply denote $\varphi_R$ by $\varphi$.
	
	Therefore we define the pseudo conditional log-likelihood function by		
		\begin{equation} \label{eq_pseudo_like_1}
		\mathcal{PL} (\theta) = \sum_{t=p+1}^T \log \left( \int_{\R^k} \prod_{i=1}^k \ind_{I_{Y_{i,t}}} \left( \lambda_{i,t}(\gamma) +x_i \right) \varphi_R (x) \d x \right),
		\end{equation}		
	and we set $\hat{\theta}_{\PL} = \underset{\theta \in \Theta}{\argmax} \ \mathcal{PL}(\theta)$.
	
	Note that \eqref{eq_pseudo_like_1} is based on the density of $(Y_t)_{p+1 \leqslant t \leqslant T}$ conditionally on $Y_1,\ldots , Y_p$ and $X_1,\ldots , X_T$. However, such a conditional density coincides with \eqref{eq_pseudo_like_1} only if the $X_t$'s are independent from the $\varepsilon_t$'s, which is not necessarily the case in our framework, hence the terminology pseudo conditional log-likelihood. One main drawback of \eqref{eq_pseudo_like_1} is the computation of many iterated integrals with some random bounds depending on the autoregressive parameters. Hence, numerical evaluations of the conditional likelihood function is complicated.

We will thus study a two-step method by first estimating the autoregressive parameters $A_1,\ldots,A_p,B$ and then the correlation matrix $R$. We first define for a positive integer $T \geqslant p+1$		
			\begin{align} \label{eq_pseudo_like_2}
			\overline{\mathcal{PL}}(\gamma) & = \sum_{t=p+1}^T \left\lbrace \sum_{i=1}^k  Y_{i,t} \log \left[1-\Phi(-\lambda_{i,t}(\gamma))\right] + (1-Y_{i,t})\log \left[\Phi(-\lambda_{i,t}(\gamma))\right] \right\rbrace \nonumber \\
			 & = \sum_{t=p+1}^T \left\lbrace \sum_{i=1}^k Y_{i,t} \log \left[ \Phi (\lambda_{i,t}(\gamma)) \right] + (1-Y_{i,t}) \log \left[ \Phi (-\lambda_{i,t}(\gamma)) \right] \right\rbrace,
			\end{align}		
		where $\Phi$ denotes the CDF of the standard Gaussian distribution in $\R$.
		
		We also set $\hat{\gamma} = \underset{\gamma \in \Gamma}{\argmax} \ \overline{\mathcal{PL}} (\gamma)$.
		
		Note that \eqref{eq_pseudo_like_2} is simply a pseudo conditional log-likelihood function. It is obtained from \eqref{eq_pseudo_like_1} under an additional independence assumption on the coordinates of $\varepsilon_t$ for any $t \in \Z$.
		
		Note also that the maximization of \eqref{eq_pseudo_like_2} can be obtained equation by equation. Indeed, we have for $1\leqslant i\leqslant k$ and $t\in\Z$
			\begin{equation}
			\lambda_{i,t}(\gamma)=\sum_{l=1}^p\sum_{j=1}^k A_l(i,j)Y_{j,t-l}+\sum_{j=1}^d B(i,j)X_{j,t-1}
			\end{equation}
and only the lines number $i$ of the matrices are involved in the previous expression. Hence maximizing \eqref{eq_pseudo_like_2} is equivalent to maximize $k$ functions with different parameters.

	Once we have an estimator $\hat{\gamma}$ obtained by optimizing \eqref{eq_pseudo_like_2}, one can plug it in the expression of the conditional likelihood function \eqref{eq_pseudo_like_1} given in the previous point, and maximize it with respect to $R$. However, this still requires numerical evaluations of multiple integrals.
	
	We prefer to adopt a pairwise conditional likelihood approach similar to the one presented in \citet{de2005pairwise} (see \citet{varin2011overview} for a review about composite likelihood approaches). For $i < j$, we define
		\begin{equation}
		R_{i,j}=\begin{pmatrix}1&r_{i,j}\\r_{i,j}&1\end{pmatrix} \quad  \text{where} \ r_{i,j} = R(i,j),
		\end{equation}
	and
		\begin{align} \label{optim_pseudo_like_2}
		\hat{r}_{i,j} & = \underset{r_{i,j}\in \mathcal{R}_{i,j}}{\argmax} \ \sum_{t=p+1}^T \log \left( \int_{I_{Y_{i,t}}-\lambda_{i,t}(\hat{\gamma})} \int_{I_{Y_{j,t}}-\lambda_{j,t}(\hat{\gamma})}  \varphi_{R_{i,j}} (x_{i},x_{j}) \d x_{i} \d x_{j} \right) \nonumber \\
		 & = \underset{r_{i,j}\in \mathcal{R}_{i,j}}{\argmax} \ \sum_{t=p+1}^T \log \left\lbrace \int_{I_{Y_{i,t}}-\lambda_{i,t}(\hat{\gamma})} \Phi \left( \left( 2Y_{j,t}-1 \right) \dfrac{\lambda_{j,t}(\hat{\gamma})+r_{i,j}x_{i}}{\sqrt{1-r_{i,j}^2}} \right) \varphi (x_{i}) \d x_{i}  \right\rbrace.
		\end{align}
	
	Note that for any $t \in \Z$, and $s_{i}, s_{j} \in \left\lbrace 0 , 1 \right\rbrace$,
		\begin{equation}
		\int_{I_{s_{i}}-\lambda_{i,t}(\gamma_0)} \Phi \left( \left( 2s_{j}-1 \right) \dfrac{\lambda_{j,t}(\gamma_0)+r_{i,j}x_{i}}{\sqrt{1-r_{i,j}^2}} \right) \varphi (x_{i}) \d x_{i} = \P \left( Y_{i,t}=s_{i},Y_{j,t}=s_{j} \mid \mathcal{F}_{t-1} \right),
		\end{equation}
	which explains the terminology pairwise conditional likelihood. However, the optimization \eqref{optim_pseudo_like_2} does not correspond to the optimization of a likelihood function. Here, we simply sum the log-densities of a given pair $(Y_{i,t},Y_{j,t})$ along time $t$, and maximize it.
	
	Finally, we denote $\hat{r} = \left( \hat{r}_{1,2},\ldots ,\hat{r}_{1,k},\hat{r}_{2,3},\ldots,\hat{r}_{2,k},\ldots \hat{r}_{k-1,k} \right)$, and for simplicity of further notations, we set for $t\in\Z$ and $\theta=(\gamma,r')\in \Theta$
		\begin{equation}\label{eq_def_ell}
\ell_{t}(\theta)=\sum_{1\leq i<j\leq k} \log \left( \int_{I_{s_i}-\lambda_{i,t}(\gamma)}\Phi\left((2Y_{j,t}-1)\frac{\lambda_{j,t}(\gamma)+r_{i,j}x_i}{\sqrt{1-r_{i,j}^2}}\right)\varphi(x_i)\d x_i\right).
		\end{equation}	
	
For deriving asymptotic properties of these M-estimators, the following additional assumptions will be needed.

  \begin{enumerate}
 	\item[\textbf{B4}]   All matrices $A_1, \ldots , A_p$ belong to a compact set $\mathcal{A} \subset \R^{k\times k}$, $B$ belongs to a compact set $\mathcal{B} \subset \R^{k\times d}$ and $C$ belongs to compact set of $\R^k$. In addition, $R$ belongs to a compact set $\mathcal{R}$ included in the set of positive definite symmetric matrices in $\R^{k\times k}$.
 	\item[\textbf{B5}] For all $t \in \Z$ and $v \in \R^d \setminus \left\lbrace 0 \right\rbrace$, we have for any $\sigma (\varepsilon_t) \vee \mathcal{F}_{t-1}$-measurable application $G$
 		\begin{equation}
 		\P \left( v' \cdot X_t =G \right) <1
 		\end{equation}		
	\end{enumerate}

\paragraph{Notes.}
\begin{enumerate}	
 \item To construct a compact subset of the set of correlation matrices, one can use Cholesky decomposition  $R=LL^T$ with $L$ lower triangular and  $L(i,i)=\sqrt{1-\sum_{\ell=1}^{i-1}L(i,\ell)^2}$. One can then assume that $\left\vert L(i,\ell)\right\vert\leq \overline{L}(i,\ell)$ where the $\overline{L}(i,\ell)'$s satisfy $\sum_{\ell=1}^{i-1}\overline{L}(i,\ell)^2<1$ for $i=1,\ldots,k$ and $1\leq \ell <i$. The number of free parameters is $(k^2-k)/2$. Cholesky decomposition will be used in our numerical experiments.
 \item 
	   Assumption \textbf{B5} guarantees that any non-trivial linear combination of the stochastic covariates at time $t$ cannot be almost surely explained by the past values and the sigma-field $\sigma (\varepsilon_t)$. Note that the error term is not necessarily independent of $X_t$, i.e. contemporaneous dependence is allowed between both quantities. By stationarity, it is sufficient to verify assumption \textbf{B5} only for $t=0$.
		 Under the more restrictive assumption of complete independence between the two processes $(X_t)_{t\in \Z}$ and $(\varepsilon_t)_{t \in \Z}$, the identifiability assumption \textbf{B5} can be replaced by the following one which is the standard assumption used for regression models in the i.i.d. setting. Details are given in the proof of Theorem \ref{theo_cons_single_global}.
		\item[\textbf{B5'}] For any $t \in \Z$, the distribution of $X_t$ is not supported by an affine hyperplane.
	\end{enumerate}

We now derive asymptotic results for the estimators defined above.  
In what follows, for an arbitrary integer $n \geqslant 2$, we denote by $\nabla f$ the gradient of a real-valued mapping $f$ defined on a subset of $\R^n$. If $n=1$, we denote by $\dot{f}$ the derivative of the mapping $f$. If $f$ is defined on $\Theta$, we will denote by $\nabla_1 f(\gamma,r)$ and $\nabla_2 f(\gamma,r)$ the two gradients of the functions $\gamma\mapsto f(\gamma,r)$ and $r\mapsto f(\gamma,r)$ evaluated at point $\theta=(\gamma,r)$.

Similarly, we denote by $\nabla^2 f$ the hessian matrix of $f$, and if $n=1$, we simply denote by $\ddot{f}$ the second derivative of the mapping $f$. If $f$ is defined on $\Theta$, setting $d_{\Gamma}$ and $d_{\mathcal{R}}$ the dimensions of $\Gamma$ and $\mathcal{R}$ respectively, we will denote
	\begin{equation}
	\nabla_1^2 f = \left( \dfrac{\partial^2 f}{\partial \gamma_i \partial \gamma_j} \right)_{1 \leqslant i,j \leqslant d_{\Gamma}}, \quad \nabla_{1,2}^2 f = \left( \dfrac{\partial^2 f}{\partial \gamma_j \partial r_i} \right)_{\underset{1 \leqslant j \leqslant d_{\Gamma}}{1 \leqslant i \leqslant d_{\mathcal{R}}}} \quad \text{and} \quad \nabla_2^2 f = \left( \dfrac{\partial^2f}{\partial r_i \partial r_j} \right)_{1 \leqslant i,j \leqslant d_{\mathcal{R}}}.
	\end{equation}
Finally, we denote by $\| \cdot \|_2$ the Euclidean norm.

\medskip

The following theorem establishes the consistency and asymptotic normality of the estimators $\hat{\theta}_{\PL}, \hat{\gamma}$ and $\hat{r}$ defined previously.

	\begin{theo}\label{theo_cons_single_global}
	Suppose that Assumptions \emph{\textbf{B1}} to \emph{\textbf{B5}} or \emph{\textbf{B1}} to \emph{\textbf{B4}} and \emph{\textbf{B5'}} hold true and $(Y_t)_{t \in \Z}$ is a strictly stationary process satisfying \eqref{eq_model_single} and \eqref{eq_lambda_single}.
		\begin{enumerate}[label=\arabic*)]
			\item \begin{enumerate}[label=\roman*)]
			\item Assume that  $\E \left( \| X_0 \|_2^2 \right) < +\infty$, then the estimators $\hat{\theta}_{\PL}, \hat{\gamma}$ and $\hat{r}$ are strongly consistent, i.e.
				\begin{equation}
				\hat{\theta}_{\PL} \underset{T \to +\infty}{\longrightarrow} \theta_0, \quad \hat{\gamma} \underset{T \to +\infty}{\longrightarrow} \gamma_0 \quad \text{and} \quad \hat{r} \underset{T \to +\infty}{\longrightarrow} r_0.
				\end{equation}
			\item Assume furthermore that $\gamma_0$ is located in the interior of $\Gamma$. Setting for $y_i \in \{0,1\}$ and $s_i\in \R$, $h_{y_i}(s_i)=y_i\log\Phi(s_i)+(1-y_i)\log\Phi(-s_i)$, we have the stochastic expansion
			\begin{equation}
			\sqrt{T}\left(\hat{\gamma}-\gamma_0\right)=\dfrac{J_{\gamma_0}^{-1}}{\sqrt{T}}\sum_{t=p+1}^T \sum_{i=1}^k \dot{h}_{Y_{i,t}}\left(\lambda_{i,t}(\gamma_0)\right)\nabla\lambda_{i,t}(\gamma_0)+o_{\P}(1),
			\end{equation}
		where
			\begin{equation}
			J_{\gamma_0}=-\sum_{i=1}^k\E\left[\ddot{h}_{Y_{i,0}}\left(\lambda_{i,0}(\gamma_0)\right)\nabla\lambda_{i,0}(\gamma_0)\nabla\lambda_{i,0}(\gamma_0)'\right].
			\end{equation}
		As a consequence, we have the asymptotic normality
			\begin{equation}
			\sqrt{T}\left(\hat{\gamma}-\gamma_0\right)\overset{\mathcal{L}}{\longrightarrow} \mathcal{N}_Q\left(0,J_{\gamma_0}^{-1}I_{\gamma_0}J_{\gamma_0}^{-1}\right),
			\end{equation}
		where
			\begin{equation}
			I_{\gamma_0}=\sum_{i=1}^k\sum_{j=1}^k\E\left[\dot{h}_{Y_{i,0}}\left(\lambda_{i,0}(\gamma_0)\right)\dot{h}_{Y_{j,0}}\left(\lambda_{j,0}(\gamma_0)\right)\nabla\lambda_{i,0}(\gamma_0)\nabla\lambda_{j,0}(\gamma_0)'\right].
			\end{equation}
			\end{enumerate}
		\item  Assume there exists $\kappa >0$ such that for all $t \in \Z$
				\begin{equation}
				\E \left( \exp \left( \kappa \| X_t \|_2^2 \right) \right) < +\infty,
				\end{equation}
			and assume furthermore that $\theta_0$ is located in in the interior of $\Theta$.
				\begin{enumerate}[label=\roman*)]
				\item Setting for any $t \in \Z$
					\begin{equation}
					m_t (\theta) = \log \left( \int_{\R^k} \prod_{i=1}^k \ind_{Y_{i,t}}\left( \lambda_{i,t}(\gamma) + x_i \right) \varphi_R(x) \d x \right),
					\end{equation}
				we have the stochastic expansion
					\begin{equation}
					\sqrt{T} \left( \hat{\theta}_{\PL} - \theta_0 \right) = \dfrac{J_{\theta_0}^{-1}}{\sqrt{T}} \nabla \PL (\theta_0) + o_{\P}(1)
					\end{equation}
				where $J_{\theta_0} = - \E \left( \nabla^2 m_0(\theta_0) \right)$.
				
				As a consequence, we have the asymptotic normality
					\begin{equation}
					\sqrt{T}\left( \hat{\theta}_{\PL} - \theta_0 \right) \overset{\mathcal{L}}{\longrightarrow} \mathcal{N}\left(0 , J_{\theta_0}^{-1} I_{\theta_0} J_{\theta_0}^{-1} \right)
					\end{equation}
				where $I_{\theta_0} = \E \left( \nabla m_0(\theta_0) \cdot \nabla m_0(\theta_0)' \right)$.
				\item Setting 
					\begin{equation}
					\E \left( \nabla^2 \ell_{0}(\theta_0) \right) = \E \left(
					\begin{array}{c;{2pt/2pt}c}
		 			& \\
					\dfrac{\partial^2 \ell_{0}}{\partial \gamma^2} (\theta_0) & \dfrac{\partial^2 \ell_{0}}{\partial \gamma \partial r}(\theta_0) \\
		 			& \\
					\hdashline[2pt/2pt]
					& \\
					\dfrac{\partial^2 \ell_{0}}{\partial r \partial \gamma}(\theta_0) & \dfrac{\partial^2 \ell_{0}}{\partial r^2}(\theta_0) \\
					& \\
					\end{array}
					\right) = \left(
					\begin{array}{c;{2pt/2pt}c}
					& \\
					L_{1,1} & L_{2,1} \\
					& \\
					\hdashline[2pt/2pt]
					& \\
					L_{1,2} & L_{2,2} \\
					& \\
					\end{array}
					\right) 
					\end{equation}
				we have the stochastic expansion
					\begin{equation}
					\sqrt{T} \left( \hat{r} - r_0 \right) = -\dfrac{1}{\sqrt{T}} \sum_{t=p+1}^T L_{2,2}^{-1} \nabla_2 \ell_t(\theta_0) - \sqrt{T} L_{2,2}^{-1} L_{1,2} \left( \hat{\gamma} - \gamma_0 \right) + o_{\P} (1).
					\end{equation}
				In particular, $\sqrt{T} \left( \hat{\gamma} - \gamma_0 , \hat{r} - r_0 \right)$ also has a Gaussian limiting distribution with mean 0.
				\end{enumerate}
		\end{enumerate}
	\end{theo}

 \paragraph{Notes.}
 \begin{enumerate}
\item 
Our asymptotic results require some moment conditions on the exogenous covariates. For consistency, only a moment of order $2$ is required. One can note that the covariates do not need to be bounded almost surely, which is important to be compatible with most of the time series models used in Econometrics. On the other hand, asymptotic normality of conditional likelihood estimators or of pairwise likelihood estimators of the correlation coefficients require an exponential moment for $\Vert X_0\Vert_2^2$. This condition is required to control the variance of the score which is a ratio of multiple integrals. Though restrictive, this condition is still satisfied for any stationary Gaussian time series. In both cases, bounded covariates sequences are allowed. An example of bounded time series in Ecology concerns the absence-presence of other species than the species of interest, as discussed in Section \ref{applipli}.
\item
Although all the asymptotic variances may have a complicated expression and are not all defined explicitly in the previous theorem, all of them have a sandwich form. Note that confidence interval for model parameters can be also obtained by parametric bootstrap.
The main difficulty for deriving asymptotic normality of the previous estimates is the presence of some parameters in the integration bounds of the multivariate Gaussian density, which leads to lengthy and tedious asymptotic expansions.
 \end{enumerate}

\section{Parametric inference for panel type data}\label{4}

Our aim is now to generalize the results given in Section \ref{sec_param_inf_single} to the case of panel data. Ecology is indeed a field where abundance data is collected from surveys at various locations and most of the time one or two times per year. Taking into account the information coming from different time points $t=1,\ldots,T$ and locations $j=1,\ldots,n$ is then primordial to get reasonable sample sizes.  

Several time series $(X_{j,t},Y_{j,t})_{1\leq t\leq T}$ are observed at different locations $j=1,\ldots,n$, with $X_{j,t}$ being a vector of exogenous covariates and $Y_{j,t}$ is the absence/presence vector observed at time $t$ and location $j$. We will assume that the same model (\ref{eq_model_single}), (\ref{eq_lambda_single}) is satisfied for all the observed time series. $M-$estimators are generally $\sqrt{nT}-$consistent in this context and the problem of statistical inaccuracy in our multivariate model due to short horizon times in many Ecological time series can be overcome if several dynamics are observed for a significant number of <<individuals>>, which corresponds to ecosystems observed at different spatial locations in our context. Note that our models and results are by no means restricted to Ecological time series.
In our real data example, the orders of magnitude of $n$ and $T$ are similar. Then for the asymptotics, we will work under a large $T$ and large $n$ scenario. Our joint limits will be obtained under the most general setup, that is when both $T$ and $n$ go to infinity but without other restrictions such as the standard diagonal path limit used in the panel data literature for which $T=T_n$ is such that $n/T\rightarrow c\neq 0$ as the index $n\rightarrow \infty$. 
We defer the reader to \citet{hsiao2022analysis}, Appendix $3$B, for a description of various approaches for large $n$ and large $T$
asymptotics. 
We will however not follow the standard assumptions used for panels data with unobserved individual or time effects, though we will allow some observed random covariates to be only location-dependent or time-dependent. See \citet{hahn2011bias} or \citet{fernandez2016individual} for asymptotic results for estimating regression parameters in panel data of this type in presence of unobserved heterogeneity.
Additionally, in the latter references, weak law of large numbers for double-indexed partial sums are derived by combining mixing assumptions on the covariates and independence across individuals, conditionally on unobserved heterogeneity. Though general results for deriving mixing properties in discrete-valued time series models with covariates have been proposed recently in \citet{truquet2023strong}, we will consider some strong (instead of weak) law of large numbers available without any stong mixing assumption.

\subsection{Law of large numbers for partial sums of double-indexed sequences}

We first present a rigorous version of the strong law of large numbers which is adapted to data collected at different time points $t$ and different locations $j$.
Finding an appropriate result in the statistical literature is not straightforward since our aim here is to consider a complex modeling with covariates $X_{j,t}$ that can contain non-time varying components or some components not depending on the location $j$. 
Technically, the problem is to find an ergodic theorem adapted to multiparameter processes whose probability distributions are invariant under some specific transformations such as the time-shift. 

There exists various formulations of the ergodic theorem for multiparameter processes. See for instance \citet{nguyen1979ergodic} for a general ergodic theorem for spatial processes or \citet{krengel2011ergodic} for various ergodic theorems available in a general context. 
In this paper, we present here two versions of the strong law of large numbers which are based on an ergodic theorem given by \citet{giap2016multidimensional}, Lemma $3.1$. A related result concerning more specifically stationary and ergodic random fields can be also be found in \citet{banys2010remarks}. An application of these results in the specific context of panel data is an important contribution of our work.

In what follows, we say that an array $(u_{n,m})_{n\in \N, m\in \N}$ converges toward $\ell$ when both $n$ and $m$ tend to $+\infty$ if
			\begin{equation}
			\forall \varepsilon >0, \exists k_0 \in \N, \ \forall n,m \geqslant k_0, |u_{n,m} - \ell | < \varepsilon.
			\end{equation}
\medskip
We consider some strong laws of large numbers for an array $(Z_{j,t})_{j,t\geq 1}$ of the form 
  \begin{equation}\label{bernshift}
Z_{j,t}=g\left(\zeta_{j,t},\zeta_{j,t-1},\ldots\right),\quad \zeta_{j,t}=(X_{j,t-1},\varepsilon_{j,t})
  \end{equation}
for some measurable mapping $g:(\R^d\times\R^k)^{\N}\rightarrow \R^{d'}$, $d'$ being a positive integer. 
Note that if $Z_{j,t}=(Y_{j,t},Y_{j,t-1},\ldots,Y_{j,t-p})$ 
where $(Y_{j,t})_{t\in\Z}$, is solution of \eqref{eq_model_single} and \eqref{eq_lambda_single}, such a representation is valid.
From now on, for any $t\in\Z$ and $j\in\N^{*}$, we set 
\begin{equation}\label{field}
\mathcal{F}_{j,t}=\sigma\left((X_{j,s},\varepsilon_{j,s}): s\leq t\right).
\end{equation}
Our first result is based on the following assumption.

\begin{enumerate}[label=\textbf{C\arabic*}]
\item For any $j\geq 1$, the stochastic process $\left(\zeta_{j,t}\right)_{t\in\Z}$ is strictly stationary and its distribution does not depend on $j$. Moreover, setting $\mathcal{F}_{j,t}=\sigma\left((X_{j,s},\varepsilon_{j,s}):s\leq t\right)$ and $\mathcal{G}_j=\vee_{t\in\Z}\mathcal{F}_{j,t}$,  the family of $\sigma$-fields $\left( \mathcal{G}_j\right)_{j \in \N^*}$ is mutually independent and for any $(j,t)\in\N^{*}\times \Z$, $\varepsilon_{j,t}$ is independent of $\mathcal{F}_{j,t-1}$.
\end{enumerate}

	\begin{prop} \label{prop_law_large_numbers_panel}
	 Suppose that \eqref{bernshift} and {\bf C1} hold true. 
  Then for all measurable mapping $f:\R^{d'}\rightarrow \R$ such that $\E \left( |f(Z_{1,1})| \times \log^+ (|f(Z_{1,1}|) \right) < +\infty$,
  we have
	\begin{equation}
		\lim_{n,T \to +\infty} \dfrac{1}{nT} \sum_{j=1}^n \sum_{t=1}^T f(Z_{j,t}) = \E (f(Z_{1,1})) \quad \mathrm{a.s.}
		\end{equation}
	\end{prop}
	
\medskip

\paragraph{Notes.}
		\begin{enumerate}
\item 
Considering $\varepsilon_{j,t}=0$ and setting $Z_{j,t}=X_{j,t}$, Proposition \ref{prop_law_large_numbers_panel} shows that the strong law of large numbers for $(X_{j,t})_{j,t}$ is always valid when we have stationarity in the time direction and i.i.d. paths, up to an integrability condition which is slightly stronger than the usual one for stationary sequences. As discussed in \citet{klesov2017moment} in the independent setup, this moment condition is optimal and cannot be avoided.

		\item 
   Suppose that $X_{j,t}=\left(V_j,W_{j,t}\right)$ takes values in $\R^{d_1}\times \R^{d_2}$ with $d=d_1+d_2$, a useful scenario for considering both time-varying and non time-varying covariates. In the context of the absence/presence of waterbirds, \citet{sebastian2010testing} used ponds or vegetation characteristics for fixed covariates and absence/presence of other species across the time for time-varying covariates. Considering stationary but not necessarily ergodic paths $(Y_{j,t})_{t\in\Z}$ as in 
  Theorem \ref{theo_stationarity} is necessary if non-time varying covariates are included in the dynamic.  
To check {\bf C1}, one possibility is to set $\overline{Z}_{j,t}=\left(W_{j,t},\varepsilon_{j,t}\right)$ and assume that the two families of random variables $\left\{V_j:j\geq 1\right\}$ and 
   $\left\{\overline{Z}_{j,t}: j\geq 1, t\in\Z\right\}$ are independent, the $V_j'$s are i.i.d and $\left\{(\overline{Z}_{j,t})_{t\in\Z}: j\geq 1\right\}$ is a family composed of independent and identically distributed stationary stochastic processes.
   In a more subtle setup, the time-varying covariates $W_{j,t}$ may also depend on $V_j$. For instance, if $W_{j,t}=H\left(\eta_{j,t},\eta_{j,t-1},\ldots\right)$ with $\eta_{j,t}=G(W_j,\kappa_{j,t})$ and $\left\{(\kappa_{j,t},\varepsilon_{j,t}): j\geq 1, t\in\Z\right\}$ is a family of i.i.d. random vectors independent from the sequence $(V_j)_{j\geq 1}$. In particular, the time-varying covariates $W_{j,t}$ can be represented as an infinite moving average sequence with  innovation sequence $\eta_{j,t}$ that can still interact with $V_j$. One can then conclude that except the independence condition imposed for different $j$, the assumptions made on the covariates are sufficiently general to be compatible with practical applications.
\item
Note that the integrability of $f\log^{+}(\vert f\vert)$ is guaranteed as soon as $\vert f\vert^{1+\delta}$ is integrable for some $\delta>0$. For simplicity, this latter condition will be always assumed in the statements of our next results.  
\item
The strong law of large numbers might also be valid for some stationary random fields $\left(Z_{j,t}\right)_{j\in \N, t\in \Z}$, without assuming independence between the sigma-fields $\mathcal{G}_j$. However, such results would require to specify the geometry of the domain indexing the spatial locations $j$. 
In panel data, independence assumption across the $j$ index is relevant when the different individuals do not interact. For spatial locations, this independence condition can be only considered as an approximation of the reality when the different locations are far enough from each other.
\end{enumerate}

We also study another more specific scenario of the form \eqref{bernshift} but without assuming independence across $j$. In this setup, the covariates have a general dynamical system representation with some common time-varying covariates.

\begin{enumerate}[label=\textbf{C'\arabic*}]
\item The covariates write as $X_{j,t}=\left(\Gamma_t,V_j,W_{j,t}\right)\in \R^{d_1}\times \R^{d_2}\times \R^{d_3},\quad d=d_1+d_2+d_3$.
Moreover, there exist measurable mappings $H:\left(\R^{d_4}\right)^{\N}\rightarrow \R^{d_1}$, $L:\left(\R^{d_5}\right)^{\N}\rightarrow \R^{d_3}$ and $G:\R^{d_2}\times\R^{d_4}\times \R^{d_6}\rightarrow \R^{d_5}$ such that for $j\geq 1$ and $t\in\Z$,
\begin{equation}\label{new2}
 \Gamma_t= H\left(\eta_t,\eta_{t-1},\ldots\right),\quad W_{j,t}=L\left(\widetilde{\eta}_{j,t},\widetilde{\eta}_{j,t-1},\ldots\right),\quad \widetilde{\eta}_{j,t}=G\left(V_j,\eta_t,\kappa_{j,t}\right),
\end{equation}
where $\left((\kappa_{j,t},\varepsilon_{j,t})\right)_{j\in \N, t\in\Z}$, $(V_j)_{j\in\N}$ and $\left(\eta_t\right)_{t\in\Z}$ are three independent sequences of  i.i.d. random variables taking values respectively in $\R^{d_6}\times \R^k$, $\R^{d_2}$ and $\R^{d_4}$.
\end{enumerate}

\begin{prop}\label{completion}
Suppose that \eqref{bernshift} and {\bf C'1} hold true. 
 Then for any measurable mapping $f:\R^{d'}\rightarrow \R$ such that
$$\E\left(\vert f(Z_{1,1})\vert\cdot \log^{+}\left(\vert f(Z_{1,1})\vert\right)\right)<\infty,$$
we have
\begin{equation}
		\lim_{n,T \to +\infty} \dfrac{1}{nT} \sum_{j=1}^n \sum_{t=1}^T f(Z_{j,t}) = \E (f(Z_{1,1})) \quad \mathrm{a.s.}
		\end{equation}
\end{prop}

\paragraph{Note.} It is not difficult to check that {\bf C'1} entails stationarity of the sequence $(\zeta_{j,t})_{t\in\Z}$, $j\geq 1$. This second scenario is interesting to take into account both local characteristics $(V_j,W_{j,t})$ of the system as well as global time-dependent and strictly exogenous characteristics $\Gamma_t$ (for instance, $\Gamma_t$ can represent the influence of macroeconomics variables in Econometrics or time-varying characteristics of the atmosphere in Ecology). Note that local covariates $W_{j,t}$ can be correlated with $\Gamma_t$ (through the impact of $\eta_t$) and can be sequentially exogenous and not necessarily strictly exogenous (since $\varepsilon_{j,t}$ and $\kappa_{j,t}$ are not necessarily independent). As assumed in \eqref{new2}, time-varying covariates can be represented as general random dynamical systems with correlated noise terms. Specification \eqref{new2} includes general moving averages and many autoregressive processes. For instance $\Gamma_t=M \Gamma_{t-1}+\eta_t$, $W_{j,t}=N W_{j,t-1}+\widetilde{\eta}_{j,t}$, where $M$ and $N$ are square matrices of appropriate dimension and spectral radius less than $1$.

\subsection{Asymptotic results for panel data}

Before generalizing our estimation results from section \ref{sec_param_inf_single}, let us specify the framework of our statistical model for panel data.

In this entire section, we consider for $j \in \N^*$ two sequences $\left( X_{j,t} \right)_{t\in \Z}$ and $\left( \varepsilon_{j,t} \right)_{t\in \Z}$ of random vectors taking values in $\R^{d}$ and $\R^{k}$ respectively. We then define an array of random binary vectors $\left( Y_{j,t} \right)_{j \in \N^*, t \in \Z}$ such that for any $j \in \N^*$, $t\in \Z$, and $ i \in \left\lbrace 1,\ldots ,k\right\rbrace$
	\begin{equation}\label{eq_model_panel}
	Y_{i,j,t} = \ind_{(0,+\infty)} \left( \lambda_{i,j,t}  + \varepsilon_{i,j,t} \right),
	\end{equation}
with
	\begin{equation} \label{eq_lambda_panel}
	\lambda_{j,t} = C+\sum_{l=1}^p A_l Y_{j,t-l} + BX_{j,t-1}.
	\end{equation}
	
\medskip
	
	\begin{note}
	All matrices $A_1,\ldots,A_p$ belong to $\R^{k\times k}$, $B \in \R^{k\times d}$, and $C\in\R^k$. These matrices are fixed in equation \eqref{eq_lambda_panel} and do not depend on $j$.
	
	In the literature of panel data, it is classical to have heterogeneity across individual $j$ by adding an intercept $c_j$ to equation \eqref{eq_lambda_panel}. The estimation of parameters $A_1, \ldots , A_p$ and $B$ is then the central problem. See for instance \citet{hsiao2022analysis}, Chapter 7, for a classical reference on this topic. However, in the context of discrete data, there is no generic way for consistently estimating the fixed regression parameters. For instance, considering $(Y_{j,t})_{1\leqslant j \leqslant N , 1 \leqslant t \leqslant T}$ for small $T$ and large $N$ in univariate binary models, the MLE of regression parameters can be inconsistent and an alternative consistent estimator is known in the logistic case, but not in the probit case.
     In this paper, we do not investigate such unobserved heterogeneity problem, but we allow observable random covariates which can vary only across individuals or only across time and which can play the rule of an observed heterogeneity. 
	\end{note}
	
\medskip

We consider the following assumptions which are the analogs of the single path analysis.

	\begin{enumerate}[label=\textbf{C\arabic*}]
 \setcounter{enumi}{1}
		\item For any $j\geq 1$ and $t\in\Z$, the random vector $\varepsilon_{j,t}$ are normally distributed
		\begin{equation}
		\varepsilon_{j,t} \sim \mathcal{N}(0,R)
		\end{equation}
	where $R$ is a correlation matrix, i.e.
		\begin{equation}
		\forall i_1,i_2 \in \left\lbrace 1,\ldots k \right\rbrace, \ -1 \leqslant R(i_1,i_2) \leqslant 1
		\end{equation}
	and $R(i,i)=1$ for any $i \in \left\lbrace 1 , \ldots , k \right\rbrace$.
		\item All matrices $A_1, \ldots , A_p$ belong to a compact set $\mathcal{A} \subset \R^{k \times k}$, $B$ belongs to a compact set $\mathcal{B} \subset \R^{k \times d}$, $C$ belongs to a compact subset of $\R^k$ and $R$ belongs to a compact set $\mathcal{R}$ included in the set of positive definite symmetric matrices in $\R^{k \times k}$.
	\item For any couple $(j,t) \in \N^* \times \Z$ and all $v \in \R^d$, we have for any $\sigma (\varepsilon_{j,t}) \vee \mathcal{F}_{j,t-1}$-measurable application $G$
		\begin{equation}
		\P \left(v' \cdot X_{j,t}=G \right) <1.
		\end{equation}
	\end{enumerate}

Keeping the same notations of Section \ref{sec_asym_res_single}, we finally focus on the estimation of the true vector of parameters $\theta_0 = (\gamma_0,r_0')$ satisfying \eqref{eq_model_panel} and \eqref{eq_lambda_panel}. We rewrite the two estimation procedures mentioned in section \ref{sec_asym_res_single} in the specific framework of panel data.

Note that by assumption \textbf{C3}, $\theta$ belongs to a compact set $\Theta$, $\gamma$ belongs to a compact set $\Gamma$, and without ambiguity in the notations, $r$ belongs to a compact set $\displaystyle{\mathcal{R}=\prod_{1 \leqslant i_1 < i_2 \leqslant k} \mathcal{R}_{i_1,i_2}}$, where each set $\mathcal{R}_{i_1,i_2}$ is a compact set included in $(-1,1)$.

	\begin{enumerate}[label=\arabic*)]
	\item We can first consider the so-called pseudo conditional log-likelihood. 
  For any $t \in \Z$ and $n \in \N^*$, the distribution of $\left( Y_{1,t},\ldots Y_{n,t} \right)$ conditionally to $\mathcal{F}_{t-1}$ is given by
		\begin{equation}
		\mathcal{L}\left( \left(Y_{1,t},\ldots , Y_{n,t} \right) \mid \vee_{j=0}^n\mathcal{F}_{j,t-1} \right) = \sum_{s \in \left(\left\lbrace 0,1 \right\rbrace^k \right)^n} \left\lbrace \prod_{j=1}^n \left( \int_{\R^k} \prod_{i=1}^k \ind_{I_{s_{i,j}}}\left( \lambda_{i,j,t}(\gamma) +x_i \right) \varphi_{R}(x) \d x \right) \right\rbrace \delta_s.
		\end{equation}
	Similarly to the single path case, we can define the pseudo conditional log-likelihood function by
		\begin{equation}\label{pseudo_like_panel_1}
		\mathcal{PL}(\theta) = \sum_{t=p+1}^T \sum_{j=1}^n \log \left( \int_{\R^k} \prod_{i=1}^k \ind_{I_{Y_{i,j,t}}} \left( \lambda_{i,j,t}(\gamma)+ x_i \right) \varphi_{R}(x) \d x \right),
		\end{equation}
	and we set $\hat{\theta}_{\PL} = \underset{\theta \in \Theta}{\argmax} \ \mathcal{PL}(\theta)$.
	
	Once again, one can note that \eqref{pseudo_like_panel_1} is based on the density of $\left( (Y_{1,t},\ldots , Y_{n,t}) \right)_{p+1\leqslant t \leqslant T}$ conditionally on $(Y_{1,1},\ldots , Y_{n,1}),\ldots , (Y_{1,p},\ldots , Y_{n,p})$ and $\left(X_{1,1}, \ldots X_{n,1}\right), \ldots,\left(X_{1,T},\ldots,X_{n,T}\right)$, but such a density coincides with \eqref{pseudo_like_panel_1} only if the $X_{j,t}$'s are independent from the $\varepsilon_{j,t}$'s, which is not necessarily the case in our framework.
	\item It is also possible to consider a two-step method as in Section 3.2, by first estimating the vector of parameters $\gamma$, then the correlation matrix $R$. We thus define
		\begin{align}\label{pseudo_like_panel_2}
		\overline{\mathcal{PL}}(\gamma) & = \sum_{t=p+1}^T \sum_{j=1}^n \sum_{i=1}^k \left( Y_{i,j,t} \log \left[1- \Phi (-\lambda_{i,j,t}(\gamma)) \right]+ (1-Y_{i,j,t}) \log \left[ \Phi (-\lambda_{i,j,t}(\gamma)) \right]\right) \nonumber \\
		& = \sum_{t=p+1}^T \sum_{j=1}^n \sum_{i=1}^k \left( Y_{i,j,t} \log \left[\Phi (\lambda_{i,j,t}(\gamma)) \right]+ (1-Y_{i,j,t}) \log \left[ \Phi (-\lambda_{i,j,t}(\gamma)) \right]\right)
		\end{align}
	and we set $\hat{\gamma} = \underset{\gamma \in \Gamma}{\argmax} \ \overline{\mathcal{PL}}(\gamma)$.
	
	In order to obtain an estimate for $R$, it is still possible to plug the estimate $\hat{\gamma}$ of $\gamma$ in the expression of function \eqref{pseudo_like_panel_1}, and maximize it with respect to $R$, but it is still complicated numerically speaking.
	
	A pairwise conditional likelihood approach similar to the one in section 3.2 is also possible, and setting for $i_1 < i_2, \ R_{i_1,i_2} = \begin{pmatrix} 1 & r_{i_1,_2} \\ r_{i_1,i_2} & 1 \end{pmatrix}$ where $r_{i_1,i_2} = R(i_1,i_2)$, we can define
		\begin{small}
		\begin{align}\label{optim_pseudo_like_panel_2}
		\hat{r}_{i_1,i_2} & = \underset{r_{i_1,i_2} \in \mathcal{R}_{i_1,i_2}}{\argmax} \ \sum_{t=p+1}^T \sum_{j=1}^n \log \left( \int_{I_{Y_{i_1,j,t}}-\lambda_{i_1,j,t}(\hat{\gamma})} \int_{I_{Y_{i_2,j,t}}-\lambda_{i_2,j,t}(\hat{\gamma})} \varphi_{R_{i_1,i_2}}(x_{i_1},x_{i_2}) \d x_{i_1} \d x_{i_2} \right) \nonumber \\
		 & = \underset{r_{i_1,i_2}\in \mathcal{R}_{i_1,i_2}}{\argmax} \ \sum_{t=p+1}^T \sum_{j=1}^n \log \left( \int_{I_{Y_{i_1,j,t}}-\lambda_{i_1,j,t}(\hat{\gamma})} \Phi \left( (2 Y_{i_2,j,t}-1) \dfrac{(\lambda_{i_2,j,t}(\hat{\gamma}) + r_{i_1,i_2} x_{i_1})}{\sqrt{1-r_{i_1,i_2}^2}} \right) \varphi (x_{i_1}) \d x_{i_1} \right).
		\end{align}
		\end{small}
	
	Note that once again, under assumption \textbf{C4}, for any $t \in \Z$, any $n \in \N^*$ and \\ $s_{i_1,1},\ldots ,s_{i_1,n}, s_{i_2,1}, \ldots ,s_{i_2,n} \in \left\lbrace 0,1 \right\rbrace$,
		\begin{multline}
		\prod_{j=1}^n \int_{I_{s_{i_1,j}}- \lambda_{i_1,j,t}(\gamma_0)} \Phi \left( (2s_{i_2,j}-1) \dfrac{\lambda_{i_2,j,t}(\gamma)+r_{i_1,i_2}e_{i_1}}{\sqrt{1-r_{i_1,i_2}^2}} \right) \varphi (e_{i_1}) \d e_{i_1} \\
		= \P \left( (Y_{i_1,1,t},Y_{i_2,1,t})= (s_{i_1,1},s_{i_2,1}),\ldots , (Y_{i_1,n,t},Y_{i_2,n,t})= (s_{i_1,n},s_{i_2,n}) \mid \mathcal{F}_{t-1}  \right),
		\end{multline}
	which explains the terminology pairwise conditional likelihood, even it still doesn't correspond to a proper likelihood function.
	
	Similarly as in the single path case, we denote $\hat{r}=(\hat{r}_{1,2},\ldots , \hat{r}_{1,k}, \hat{r}_{2,3},\ldots ,\hat{r}_{2,k},\ldots \hat{r}_{k-1,k})$, and we set for any $j \in \N^*$, $t \in \Z$ and $\theta = (\gamma,r') \in \Theta$
		\begin{equation}
		\ell_{j,t}(\theta) = \sum_{1 \leqslant i_1 < i_2 \leqslant k} \log \left( \int_{I_{Y_{i_1,j,t}}-\lambda_{i_1,j,t}(\gamma)} \Phi \left( (2 Y_{i_2,j,t} -1) \dfrac{\lambda_{i_2,j,t}(\gamma) + r_{i_1,i_2} x_{i_1}}{\sqrt{1-r_{i_1,i_2}^2}} \right) \varphi(x_{i_1}) \d x_{i_1} \right).
		\end{equation}
	\end{enumerate}
	
	\begin{theo}\label{theo_cons_panel_global}
	Suppose that either Assumption \emph{\textbf{C1}} or Assumption \emph{\textbf{C'1}} is valid and Assumptions \emph{\textbf{C2-C3-C4}} hold true. For any $j \in \{1,\ldots,n\}$, there then exists a unique strictly stationary process $(Y_{j,t})_{t\in\Z}$  satisfying \eqref{eq_model_panel} and \eqref{eq_lambda_panel}. Moreover, the following assertions are valid.
	
		\begin{enumerate}[label=\arabic*)]
		\item \begin{enumerate}[label=\roman*)]
			  \item Assume that for all $j \in \N^*$ and $t \in \Z$, $\E \left( \|X_{j,t}\|_2^{2+\delta} \right) < +\infty$ for some $\delta >0$, then the estimators $\hat{\theta}_{\PL}, \hat{\gamma}$ and $\hat{r}$ are strongly consistent, i.e.
			  	\begin{equation}
			  	\hat{\theta}_{\PL} \underset{n,T \to \infty}{\longrightarrow} \theta_0, \quad \hat{\gamma} \underset{n,T \to \infty}{\longrightarrow} \gamma_0 \quad \text{and} \quad \hat{r} \underset{n,T \to \infty}{\longrightarrow} r_0.
			  	\end{equation}
			  \item Assume furthermore that $\gamma_0$ is located in the interior of $\Gamma$. Setting for $y_i \in \{ 0,1 \}$ and $s_i \in \R, \ h_{y_i}(s_i) = y_i \log \left( \Phi (s_i) \right) + (1-y_i) \log \left( \Phi(-s_i) \right)$, we have the stochastic expansion
			  	\begin{equation}
			  	\sqrt{nT} \left( \hat{\gamma} - \gamma_0\right) = \dfrac{J_{\gamma_0}^{-1}}{\sqrt{nT}} \sum_{j=1}^n \sum_{t=p+1}^T \dot{h}_{Y_{i,j,t}}(\lambda_{i,j,t}(\gamma_0)) \nabla \lambda_{i,j,t}(\gamma_0) + o_{\P}(1),
			  	\end{equation}
			  when $n,T \to \infty$, and where 
			  	\begin{equation}
			  	J_{\gamma_0} = - \sum_{i=1}^k \E \left[\ddot{h}_{Y_{i,1,1}}(\lambda_{i,1,1}(\gamma_0)) \nabla \lambda_{i,1,1}(\gamma_0) \nabla \lambda_{i,1,1}(\gamma_0)' \right].
			  	\end{equation}
			  As a consequence, we have the asymptotic normality
			  	\begin{equation}
			  	\sqrt{nT} \left( \hat{\gamma} - \gamma_0 \right) \overset{\mathcal{L}}{\longrightarrow} \mathcal{N} \left(0 , J_{\gamma_0}^{-1} I_{\gamma_0} J_{\gamma_0}^{-1} \right)
			  	\end{equation}
			  when $n,T \to \infty$, where
			  	\begin{equation}
			  	I_{\gamma_0} = \sum_{i=1}^k \sum_{j=1}^k \E \left[ \dot{h}_{Y_{i,1,1}}(\lambda_{i,1,1}(\gamma_0)) \dot{h}_{Y_{j,1,1}}(\lambda_{j,1,1}(\gamma_0)) \nabla \lambda_{i,1,1}(\gamma_0) \nabla \lambda_{j,1,1}(\gamma_0)' \right].
			  	\end{equation}
			  \end{enumerate}
		\item  Assume there exists $\kappa >0$ such that for all $j \in \N^*$ and $t \in \Z$
				\begin{equation}\label{eq_kappa_cons_1}
				\E \left( \exp \left( \kappa \| X_{j,t} \|_2^2 \right) \right) < +\infty,
				\end{equation}
			and assume furthermore that $\theta_0$ is located in in the interior of $\Theta$.
				\begin{enumerate}[label=\roman*)]
				\item Setting for any $j \in \N^*$ and $t \in \Z$
					\begin{equation}
					m_{j,t} (\theta) = \log \left( \int_{\R^k} \prod_{i=1}^k \ind_{Y_{i,j,t}}\left( \lambda_{i,j,t}(\gamma) + x_i \right) \varphi_R(x) \d x \right),
					\end{equation}
				we have the stochastic expansion
					\begin{equation}
					\sqrt{nT} \left( \hat{\theta}_{\PL} - \theta_0 \right) = \dfrac{J_{\theta_0}^{-1}}{\sqrt{nT}} \nabla \PL (\theta_0) + o_{\P}(1)
					\end{equation}
				when $n,T \to +\infty$, and where $J_{\theta_0} = - \E \left( \nabla^2 m_{1,1}(\theta_0) \right)$.
				
				As a consequence, we have the asymptotic normality
					\begin{equation}
					\sqrt{nT}\left( \hat{\theta}_{\PL} - \theta_0 \right) \overset{\mathcal{L}}{\longrightarrow} \mathcal{N}\left(0 , J_{\theta_0}^{-1} I_{\theta_0} J_{\theta_0}^{-1} \right)
					\end{equation}
				when $n,T \to +\infty$, where $I_{\theta_0} = \E \left( \nabla m_{1,1}(\theta_0) \cdot \nabla m_{1,1}(\theta_0)' \right)$.
				\item Setting
				\begin{equation}\label{eq_stoc_exp_3}
					\E \left( \nabla^2 \ell_{0}(\theta_0) \right) = \E \left(
					\begin{array}{c;{2pt/2pt}c}
		 			& \\
					\dfrac{\partial^2 \ell_{0}}{\partial \gamma^2} (\theta_0) & \dfrac{\partial^2 \ell_{0}}{\partial \gamma \partial r}(\theta_0) \\
		 			& \\
					\hdashline[2pt/2pt]
					& \\
					\dfrac{\partial^2 \ell_{0}}{\partial r \partial \gamma}(\theta_0) & \dfrac{\partial^2 \ell_{0}}{\partial r^2}(\theta_0) \\
					& \\
					\end{array}
					\right) = \left(
					\begin{array}{c;{2pt/2pt}c}
					& \\
					L_{1,1} & L_{2,1} \\
					& \\
					\hdashline[2pt/2pt]
					& \\
					L_{1,2} & L_{2,2} \\
					& \\
					\end{array}
					\right) 
					\end{equation}
				we have the stochastic expansion
					\begin{equation}
					\sqrt{nT} \left(\hat{r} - r_0 \right) = -\dfrac{1}{\sqrt{nT}} \sum_{t=p+1}^T \sum_{j=1}^n L_{2,2}^{-1} \nabla_2 \ell_{j,t}(\theta_0) - \sqrt{nT} L_{2,2}^{-1} L_{1,2} \left( \hat{\gamma} - \gamma_0 \right) + o_{\P} (1)
					\end{equation}
				when $n,T \to +\infty$.
				
				In particular, $\sqrt{nT} \left( \hat{\gamma} - \gamma_0 , \hat{r} - r_0 \right)$ also has a Gaussian limiting distribution with mean 0.
				\end{enumerate}
		\end{enumerate}
	\end{theo}

 \paragraph{Notes.}

 \begin{enumerate}
    \item 
    As expected, the convergence rate of the different estimators is $\sqrt{nT}$ and $(n,T)$ can grow to infinity with some arbitrary rates.
    \item 
    The law of large numbers derived for two-parameters partial sums plays a central role in the proof of Theorem \ref{theo_cons_panel_global}. With respect to the single path analysis, we simply need to reinforce the assumption of square- integrability of the covariates by the existence of a moment greater than $2$. In the supplementary file, we provide some general asymptotic results for parametric M-estimators, in the setup of longitudinal observations and under Assumption {\bf C1} or Assumption {\bf C'1}.   These generic results can be also applied to general autoregressive time series models with exogenous covariates.
 \end{enumerate}

\section{Numerical applications}\label{5}

In this section, we first compare the two inference strategies mentioned in the previous section, in the case of simulated data. We will then apply the two-steps method to a real dataset of fish catches in the North Sea, in order to get some insights about the interactions between some particular species.

\subsection{Simulated data and methods comparison}

Our aim here is to compare the performances of our two estimation strategies. To do so, we perform 1000 simulations of an absence-presence process of two species, depending on one exogenous variable.

For each simulation, the time dimension is fixed at $T=100$, and the number of time series paths is $n=50$. For simulation number $s \in \left\lbrace1, \ldots , 1000 \right\rbrace$, we thus denote $\left(Y_{j,t}^{(s)}\right)_{t\in \Z}$ the time series valued in $\left\lbrace 0, 1 \right\rbrace^{2}$ for path $j \in \left\lbrace 1, \ldots , 50 \right\rbrace$.

Furthermore, we denote $ \left( X_{j,t}^{(s)}\right)_{t \in \Z}$ the process of covariates valued in $\R$ that impacts the process $\left(Y_{j,t}^{(s)}\right)_{t\in \Z}$ in our study. For any $j \in \left\lbrace 1,\ldots,50\right\rbrace$, $X_{j,t}^{(s)}$ is a stationary ARMA$(3,1)$ process.

The covariates process is centered and standardized before simulating the process $\left( Y_{j,t}^{(s)} \right)_{1 \leqslant t \leqslant T}$ according to the following equation
	\begin{equation}
	\forall t \in \Z, \forall i \in \left\lbrace 1,2 \right\rbrace, \ Y_{i,j,t}^{(s)} = \ind_{(0,+\infty)} \left( \lambda_{i,j,t}^{(s)} + \varepsilon_{i,j,t}^{(s)} \right)
	\end{equation}
where $\lambda_{j,t}^{(s)} = A \cdot Y_{j,t-1}^{(s)} + B \cdot X_{j,t-1}^{(s)} + C$ with
	\begin{equation}
	A = \left( \begin{array}{cc}
	0.3 & -0.5 \\
	0.2 & 0.7 
	\end{array} \right), \quad B = \left( \begin{array}{c}
	-0.5 \\
	0.6
	\end{array} \right),  \quad 
	C = \left( \begin{array}{c}
	0.2 \\ 0.4 
	\end{array} \right)
	\end{equation}
and the variables $\varepsilon_{j,t}^{(s)}$ are independent and identically distributed, $\varepsilon_{j,t}^{(s)} \sim \mathcal{N}(0,R)$, where the correlation matrix $R$ is given by
	\begin{equation}
	R = \left( \begin{array}{cc}
	1 & -0.2 \\
	-0.2 & 1
	\end{array}
	\right).
	\end{equation}
One can find in Tables \ref{table_1step} and \ref{table_2steps} our estimation results with both 1-step and 2-steps methods discussed in Section 4.

\medskip

For the initialization of coefficient $A(k,l)$ in the optimization procedure, we actually compute over all 50 paths how many times the presence of species $l$ at time $t$ is followed by the presence of species $k$ at time $t+1$, and how many times it is not. These two latter sums are pondered by 1 or -1 whether we consider the positive or negative impact of species $l$ on species $l$. We then divide the obtained result by the number of time steps over all 50 sites to get an initial guess for coefficient $A(k,l)$.

The initialization of coefficient $B(k)$ is simply given by the mean over the 50 paths observed of the correlation estimated between $Y_{k,j,t}^{(s)}$ and $X_{j,t-1}^{(s)}$. For each estimation of the correlation, we proceed here as if for $1 \leqslant k \leqslant 50$, $ \left( Y_{k,j,t}^{(s)}, X_{j,t}^{(s)}\right)_{1 \leqslant t \leqslant 100}$ were i.i.d. random variables.

Coefficients of type $C(k)$ are simply initialized by the proportion $p_k$ of presence over time $t$ and space $j$ for species $k$, which is then transformed by $2 \times p_k -1$.

Finally, for a coefficient $R(k,l)$, we simply estimate the correlation between $Y_{k,j,t}^{(s)}$ and $Y_{l,j,t}^{(s)}$. Here, we also proceed as if the random variables $\left( Y_{k,j,t}^{(s)} , Y_{l,j,t}^{(s)} \right)_{\underset{1 \leqslant j \leqslant 50}{1 \leqslant t \leqslant 100}}$ were i.i.d.

\medskip

In order to compare the two methods, we fixed in both cases the same stopping criterion to end the numeric optimization, and used the same 100 CPUs. The estimates of $\hat{\gamma}$ in the 2-steps method were obtained using the Nelder-Mead simplex algorithm (see \citet{nelder1965simplex}), while the estimates of $\hat{R}$ in the 2-steps method, as well as the estimates of $\hat{\theta}$ in the 1-step method, were computed with the COBYLA algorithm (see \citet{powell1994direct}) due to the constraints on the correlation matrix. 

\medskip
	
 Table \ref{table_1step} and Table \ref{table_2steps} exhibit the bias and variance of the different estimates. Note that for this moderate sample size, the 2-steps method performed significantly better both in terms of bias and MSE, while the variances of the estimates are very similar for the two methods. We think that this larger bias can be due to numerical issues. Indeed, the 1-step method is very difficult to implement with a numerical approximation of multiple integrals. Moreover, the presence of many local maxima makes the optimization of \eqref{pseudo_like_panel_1} quite difficult. The computation time for the 2-steps method was about 21 minutes, while it was about 2 hours and 30 minutes for the 1-step method. Such simulation shows that the 2-steps method is an interesting choice since it reduces the computational cost of the 1-step method and is still accurate for moderate sample sizes. Only the 2-step method will be used for the real data application given in the next section.

\begin{table}[H]
\centering
\begin{tabular}{cccccc}
Parameter & True Value & Mean Estimate & MSE & Bias & Variance \\
 \hline
$A(1,1)$ & 0.3 & 0.358 & 0.005 & 0.058 & 0.002 \\
$A(1,2)$ & -0.5 & -0.338 & 0.283 & 0.162 & 0.002 \\
$A(2,1)$ & 0.2 & 0.186 & 0.002 & -0.014 & 0.002 \\
$A(2,2)$ & 0.7 & 0.680 & 0.003 & -0.020 & 0.003 \\
\hline
$B(1)$ & -0.5 & -0.411 & 0.001 & -0.011 & 0.001 \\
$B(2)$ & 0.6 & 0.601 & 0.001 & 0.001 & 0.001 \\
\hline
$C(1)$ & 0.2 & 0.014 & 0.036 & -0.186 & 0.001 \\
$C(2)$ & 0.4 & 0.427 & 0.004 & 0.027 & 0.003 \\
\hline
$R(1,2)$ & -0.2 & -0.199 & 0.001 & -0.001 & 0.001 \\
\end{tabular}
\caption{Estimates Results for the 1-step Method. The mean is taken over the 1000 simulations}\label{table_1step}
\end{table}

\begin{table}[H]
\centering
\begin{tabular}{cccccc}
Parameter & True Value & Mean Estimate & MSE & Bias & Variance \\
 \hline
$A(1,1)$ & 0.3 & 0.298 & 0.002 & -0.002 & 0.002\\
$A(1,2)$ & -0.5 & -0.503 & 0.003 & -0.003 & 0.003 \\
$A(2,1)$ & 0.2 & 0.200 & 0.002 & 0.000 & 0.002 \\
$A(2,2)$ & 0.7 & 0.700 & 0.003 & 0.000 & 0.003 \\
\hline
$B(1)$ & -0.5 & -0.400 & 0.001 & 0.100 & 0.001 \\
$B(2)$ & 0.6 & 0.600 & 0.001 & 0.000 & 0.001 \\
\hline
$C(1)$ & 0.2 & 0.203 & 0.003 & 0.003 & 0.003 \\
$C(2)$ & 0.4 & 0.401 & 0.003 & 0.001 & 0.003 \\
\hline
$R(1,2)$ & -0.2 & -0.199 & 0.001 & 0.001 & 0.001 \\
\end{tabular}
\caption{Estimates Results for the 2-steps Method. The mean is taken over the 1000 simulations}\label{table_2steps}
\end{table}

\subsection{Fishing data}\label{applipli}

The following example relies on a dataset of fish catches provided publicly by the website of the International Council for the Exploration of the Sea (ICES), Copenhagen. More precisely, this dataset contains the Catch Per Unit Effort (CPUE) per length per subarea from the North Sea International Bottom Trawl Survey (NS-IBTS), and was downloaded the $29^{\text{th}}$ may 2024 \citep{ICES2024}.

\medskip

The CPUE corresponds here to the number of catches per hour, and is somehow an indicator of abundance for each fish species.

\medskip

The downloaded NS-IBTS dataset thus contains the CPUE observed for each subarea in the North Sea (see Figure \ref{map_NS_IBTS}) in Area 1, from 1971 to 2024. There are nine fishes species studied in this dataset, scientifically named \emph{Clupea Harengus, Gadus Morhua, Melanogrammus Aeglefinus, Merlangius Merlangus, Pleuronectes Platessa, Pollachius Virens, Scomber Scombrus, Sprattus Sprattus} and \emph{Trisopterus Esmarkii}, which are subdivided in different length classes. 

\medskip
	\begin{center}
	\begin{minipage}{0.6\textwidth}
	\begin{figure}[H]
	\centering
	\includegraphics[scale=0.7]{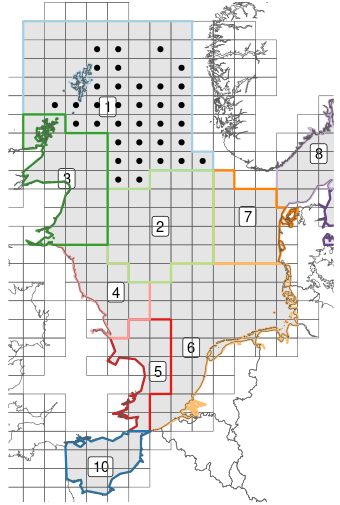}
	\caption[font=scriptsize]{Areas and subareas of the NS-IBTS dataset. The dots correspond to the 36 subareas used in our study.\label{map_NS_IBTS}}
	\end{figure}
	\end{minipage}
	\end{center}
	
\medskip

For practical reasons, our study focuses on the time interval 1975-2024, and we first aggregate the CPUEs allocated by length for each species, then we remove the subareas containing too much missing data, which leaves us with 36 subareas (see Figure \ref{map_NS_IBTS}), corresponding to 36 paths over time for our absence-presence process. Moreover, we compute the first tiercile of the CPUEs for each species over time and space, in order to set the absence-presence index to 0 if the CPUE observed is smaller than this tiercile, and we set the absence-presence index to 1 otherwise. Hence, a value 0 does not necessarily mean that the species is absent, but that the species is somehow less abundant than usual.

Finally, we impute the few remaining missing data at random, according to a Bernoulli distribution $\mathcal{B}(p)$ where $p$ is the proportion of ones obtained in the whole dataset for the species for which we want to impute the absence-presence.

We decide to focus on the following three particular species: \emph{Gadus Morhua, Merlangius Merlangus and Sprattus Sprattus}. The processes of absence-presence for the six other species are used here as exogenous variables. The model we fit to our data is thus
	\begin{equation}
	\forall t \in \Z, \forall j \in \left\lbrace1, \ldots,36\right\rbrace, \forall i \in \left\lbrace 1,2,3 \right\rbrace, \ Y_{i,j,t} = \ind_{(0,+\infty)} \left( \lambda_{i,j,t} + \varepsilon_{i,j,t} \right),
	\end{equation}
where $\lambda_{j,t} = A \cdot Y_{j,t-1} + B \cdot X_{j,t-1} +C$ with $A \in \R^{3\times 3}, \ B \in \R^{3\times 6}, X_{j,t-1} \in \left\lbrace 0,1 \right\rbrace^{6}$ and $C \in \R^3$. $X_{j,t}$ corresponds to the vector of covariates, i.e. the absence-presence process of the six non-studied species, and the $\varepsilon_{j,t}$ are i.i.d with normal distribution $\mathcal{N}(0,R)$, $R \in \R^{3\times 3}$ being a correlation matrix.

By doing so, we actually obtain a sub-model of the complete model with nine species, but with no covariates.

\medskip

Due to to the long time required by the 1-step estimation strategy, we only estimate the values of parameters $A,  B, C$ and $R$ with the 2-steps method. Estimation results are stored in Table \ref{Table_NS_IBTS}.

\medskip

Since we are interested in the interpretation of the model's parameters, we then want to test the significancy of the estimated parameters. Denoting 
	\begin{equation}
	\theta = \left( A(1,1),\ldots, A(3,3),B(1,1),\ldots, B(3,6), C(1),C(2),C(3),R(2,1),R(3,1),R(3,2) \right)',
	\end{equation}
$\theta_0$ the true value of the parameter $\theta$ and $\hat{\theta}$ its estimation, we have by Theorem \ref{theo_cons_panel_global}
	\begin{equation}
	\sqrt{nT} \left( \hat{\theta} - \theta_0 \right) \overset{\mathcal{L}}{\longrightarrow} \mathcal{N}(0,V_{\theta_0})
	\end{equation}
where $V_{\theta_0}$ is an unknown covariance matrix. However, the quantiles of the distribution $\mathcal{N}(0,V_{\theta_0})$ can be approximated by a bootstrap estimation (see \citet{diciccio1996bootstrap} for a review of bootstrap confidence intervals). To do so, we simulate a thousand absence-presence processes with 36 paths and duration 50 with the estimated parameter $\hat{\theta}$. These simulations then give us a thousand estimates of $\hat{\theta}$, denoted $\hat{\theta}^{*}_1,\ldots,\hat{\theta}^{*}_{1000}$.

Computing the empirical quantiles of these latter estimates, we obtain good approximations of the quantiles of the limit distribution $\mathcal{N}(0,V_{\theta_0})$.

We can then deduce easily some confidence intervals for the coefficients of parameter $\theta_0$. These results are given in Table \ref{Table_NS_IBTS}.

\medskip

Finally, we propose some interpretations concerning the results obtained. For example, coefficient $A(3,2) \simeq 0.393$ being significantly positive, it suggests that the presence of species 2 (\emph{Merlangius merlangus}) at time $t-1$ has a positive impact on the presence of species 2 (\emph{Sprattus sprattus}) at time $t$. Indeed, it appears in a survey (see \citet{CEFAS2024}, consulted online the $10^{\mathrm{th}}$ of June, 2024) conducted by the Center for Environment Fishes and Aquaculture Science (CEFAS), monitoring fishes stomachs, that species \emph{Sprattus Sprattus} is a particular predator for the eggs of species like \emph{Merlangius merlangus}.

Coefficient $A(2,3) \simeq 0.448$ is however significantly positive, meaning that the presence of species \emph{
Sprattus sprattus} at time $t-1$ has also a positive impact on the presence of species \emph{Merlangius merlangus} at time $t$, which might seem in contradiction with the previous interpretation. This actually discloses an underlying phenomenon that allows these two species to live in symbiosis. This insight is indeed confirmed by the estimated coefficient $R(2,3) \simeq 0.158$, which is also significantly positive.

The model proposed in this paper also allows to give some insights about fishes interactions that might not be obvious at first sight. This is the case for species \emph{Scomber scombrus} and \emph{Merlangius merlangus}, as coefficient $B(2,5) \simeq -0.184$ is significantly negative. This means that \emph{Scomber scombrus} has a negative impact on the presence of \emph{Merlangius merlangus}, although there is no direct indication that these two predator species are in competition.

\medskip

Many other interpretations could be deduced from Table \ref{Table_NS_IBTS}, but this is actually out of the scope of this paper, where the numerical study is just here given as an illustration.

\begin{table}[H]
\centering
\begin{tabular}{ccc}
\hline
Parameter & Estimate & 95\% Confidence Interval \\
\hline
$A(1,1)$ & 0.841 & [0.715;0.982]\\
$A(1,2)$ & 0.313 & [0.173;0.443]\\
$A(1,3)$ & 0.218 & [-0.017;0.430]\\
$A(2,1)$ & 0.144 & [0.005;0.290]\\
$A(2,2)$ & 1.004 & [0.861;1.135]\\
$A(2,3)$ & 0.448 & [0.197;0.679]\\
$A(3,1)$ & 0.117 & [-0.062;0.297]\\
$A(3,2)$ & 0.393 & [0.195;0.576]\\
$A(3,3)$ & 0.000 & [-0.002;0.002]\\
\hline
$B(1,1)$ & -0.066 & [-0.201;0.079]\\
$B(1,2)$ & 0.122 & [-0.016;0.272]\\
$B(1,3)$ & 0.065 & [-0.074;0.205]\\
$B(1,4)$ & 0.250 & [0.086;0.401]\\
$B(1,5)$ & 0.060 & [-0.094;0.223]\\
$B(1,6)$ & 0.103 & [-0.043;0.266]\\
$B(2,1)$ & 0.513 & [0.343;0.656]\\
$B(2,2)$ & 0.233 & [0.091;0.379]\\
$B(2,3)$ & 0.198 & [0.046;0.334]\\
$B(2,4)$ & -0.198 & [-0.352;-0.041]\\
$B(2,5)$ & -0.184 & [-0.355;-0.003]\\
$B(2,6)$ & 0.094 & [-0.052;0.241]\\
$B(3,1)$ & 0.004 & [-0.025;0.027]\\
$B(3,2)$ & 0.023 & [-0.122;0.152]\\
$B(3,3)$ & -0.022 & [-0.166;0.119]\\
$B(3,4)$ & -0.009 & [-0.053;0.056]\\
$B(3,5)$ & -0.017 & [-0.102;0.091]\\
$B(3,6)$ & 0.093 & [-0.122;0.250]\\
\hline
$C(1)$ & -0.644 & [-0.823;-0.452]\\
$C(2)$ & -0.708 & [-0.873;-0.505]\\
$C(3)$ & -1.697 & [-1.889;-1.442]\\
\hline
$R(1,2)$ & 0.229 & [0.200;0.407]\\
$R(1,3)$ & -0.020 & [-0.134;0.134]\\
$R(2,3)$ & 0.158 & [0.041;0.381]\\
\hline
\end{tabular}
\caption{Estimation results for the NS-IBTS dataset.}\label{Table_NS_IBTS}
\end{table}

\section{Conclusion}

The first aim of this paper was to provide a general framework for modeling multiple binary time series which can be compatible with 
many practical applications. In particular, we developed a general stationary theory as well as several asymptotic results for parametric estimation in multivariate autoregressive probit models. 
Our setup is very flexible since it allows to consider various kind of exogenous covariates, a single path or a longitudinal analysis and some pseudo-likelihood methods which are easy to implement numerically. A solid theoretical background is also provided  and makes a bridge between some important results in probability theory and the asymptotics of pseudo-likelihood estimators for panel data. 

Let us mention that our inference results could also be extended to other parametric models such as the multiple logistic model used for instance in \citet{guanche2014autoregressive} or \citet{sebastian2010testing}. However, the use of a multivariate Gaussian distribution
for the noise component $\varepsilon_t$ is a natural choice and allows to model many possible interactions between the binary coordinates.
A more challenging extension concerns semi-parametric inference of the autoregressive parameters, when the probability distribution of the noise term $\varepsilon$ is not specified. In the case of binary choice models, corresponding to $k=1$ in our context, \citet{manski1975maximum}  proposed the so-called maximum score estimator which was studied for a dynamic model by \citet{de2011dynamic}. A multivariate version of such estimator could be interesting.

Another theoretical limitation of our study is the stationarity assumption imposed to the exogenous variables considered in our model.  
Indeed, many time series found for instance in macroeconomics or climate analysis cannot be considered as stationary.
However, in a non-stationary framework, the asymptotic results could strongly depend on the dynamic of the covariates. See for instance \citet{debaly2023mixing}
for explosive covariates in count time series and also \citet{park2000nonstationary} for unit-root covariates in binary choice model.

Another way of improvement would be to extend the pairwise approach to composite likelihood obtained from sub-vectors of higher dimension such as triplets. We could then get more efficient estimators but a compromise between efficiency and numerical accuracy has to be done.   

\section{Supplementary material}

\subsection{Proof of Theorem \ref{theo_stationarity}}

We first remind Theorem 3 from \citet{truquet2023strong}.

	\begin{theo}\label{theo_truquet_strong}
	Let $(\zeta_t)_{t\in \Z}$ be a stationary and ergodic process valued in $\R^n$, and $E$ be a finite or countable set. For all $z \in \R^n$, we consider an application $F_z : E \longrightarrow E$ such that the mapping
		\begin{equation}
		\begin{array}{l}
		E \times \R^n \longrightarrow E \\
		(x,z) \longmapsto F_z(x)
		\end{array}
		\end{equation}
	is measurable. For all $t,s$ in $\Z$ with $s \leqslant t$, we denote the random mapping from $E$ to $E$
		\begin{equation}
		F_s^t = F_{\zeta_t} \circ \cdots \circ F_{\zeta_s},
		\end{equation}
	and we assume there exists $m \in \N^*$ such that
		\begin{equation}
		\P \left( \mathrm{Card} \left( F_1^m(E) \right) =1 \right) >0.
		\end{equation}
	Then there exists a unique stationary and ergodic process $(Y_t)_{t \in \Z}$ such that
		\begin{equation}
		\forall t \in \Z, \ Y_t = F_{\zeta_t}(Y_{t-1}).
		\end{equation}
	Furthermore, $(Y_t)_{t\in \Z}$ admits a Bernoulli shift representation
		\begin{equation}
		\forall t \in \Z, \ Y_t = G(\zeta_t,\zeta_{t-1},\ldots )
		\end{equation}
	where $G : \left(\R^n \right)^{\N} \longrightarrow E$ is a measurable mapping.
	\end{theo}
	
	\medskip

	\begin{enumerate}[labelindent=0pt,labelwidth=!,wide,label=\arabic*.]
	\item We set $E = \left( \left\lbrace 0,1 \right\rbrace^k \right)^p$ and we define the random mapping
		\begin{equation}
		F_{\zeta_t}(y_1,\ldots y_p) = \left( \ind_{g(y_1,\ldots ,y_{p-1},\zeta_t) >0} , y_1,\ldots , y_{p-1} \right),
		\end{equation}
	where for $y \in \R^k$, $\ind_{y>0}$ denotes the vector of $\R^k$ with binary coordinates $\ind_{y_i >0}, \ 1 \leqslant i \leqslant k$.
	
		\begin{enumerate}[label=$\triangleright$]
		\item We first prove that
			\begin{equation}
			\P \left( \mathrm{Card}\left( F_1^p(E) \right) =1 \right) >0.
			\end{equation}
		To this end, note that
			\begin{equation}\label{eq_single_element}
			\bigcap_{i=1}^k C_i^{j_i} \subset \left\lbrace \mathrm{Card}\left( F_1^p(E) \right) =1 \right\rbrace.
			\end{equation}
		Indeed, on the event $ \cap_{i=1}^k C_i^{j_i} $, we have for any $(y_1,\ldots , y_p) \in E$
			\begin{equation}
			F_{\zeta_1}(y_1,\ldots ,y_p) = (j,y_1,\ldots , y_{p-1}),
			\end{equation}
		where $j=(j_1,\ldots , j_k)$. We then have
			\begin{align}
			F_1^2(y_1,\ldots , y_p)  & = F_{\zeta_2}(j,y_1,\ldots ,y_{p-1}) \nonumber \\
			 & = \left( \ind_{g(j,y_1,\ldots,y_{p-1},\zeta_2) >0} , j,y_1,\ldots ,y_{p-2} \right) \nonumber \\
			 & = (j,j,y_1,\ldots ,y_{p-2}).
			\end{align}
		Iterating this latter operation, we obtain $F_1^p(y_1,\ldots,y_p)=(j,\ldots,j)$, thus on the event $\cap_{i=1}^k C_i^{j_i}$, we have $F_1^p(E)=\left\lbrace j \right\rbrace$ and \eqref{eq_single_element} follows.
		
		By Assumption \textbf{A2}, we then have $ \P \left( \mathrm{Card}\left( F_1^p (E) \right) =1 \right) >0$, and Theorem \ref{theo_truquet_strong} gives the existence of an ergodic process $(Z_t)_{t\in\Z}$ such that
			\begin{equation}
			\forall t \in \Z, \ Z_t = F_{\zeta_t}(Z_{t-1}).
			\end{equation}
		\item We then define for all $t \in \Z$
			\begin{equation}
			Y_t = Z_t(1).
			\end{equation}
		$(Y_t)_{t\in \Z}$ is an ergodic process and for all $t \in \Z$
			\begin{equation}
			Y_t = \ind_{g(Z_{t-1}(1),\ldots,Z_{t-1}(p),\zeta_t)>0} = \ind_{g(Y_{t-1},Z_{t-1}(2),\ldots,Z_{t-1}(p),\zeta_t)>0}.
			\end{equation}
		Yet, one can show by induction that for $t\in\Z$ and $l \in \left\lbrace2,\ldots,p\right\rbrace, \ Z_t(l) = Z_{t-l+1}(1)$, thus
			\begin{equation}
			\forall t \in \Z, \ Y_t= \ind_{g(Y_{t-1},\ldots,Y_{t-p},\zeta_t)>0}.
			\end{equation}
		Furthermore, Theorem \ref{theo_truquet_strong} specifies the existence of a measurable mapping
			\begin{equation}
			G : \left( \R^n \right)^{\N} \longrightarrow E
			\end{equation}
		such that for all $t \in \Z, \ Z_t= G (\zeta_t,\zeta_{t-1},\ldots)$. Taking $H$ the composition $\pi_1 \circ G$, where $\pi_1$ denotes the projection on the first coordinate, we then obtain the Bernoulli shift representation \eqref{eq_Bernoulli_shift}.
		\item We finally prove the uniqueness of process $(Y_t)_{t\in \Z}$. Let $\left( Y_t'\right)_{t\in\Z}$ be a stationary and ergodic process satisfying
			\begin{equation}
			\forall t \in \Z, \ Y_t' = \ind_{g(Y_{t-1}',\ldots,Y_{t-p}',\zeta_t)>0}
			\end{equation}
		and denote $W_t=\left(Y_t',\ldots,Y_{t-p+1}' \right)$. The process $(W_t)_{t\in\Z}$ is ergodic, and for all $t \in \Z$
			\begin{align}
			F_{\zeta_t}(W_{t-1}) & = \left( \ind_{g(Y_{t-1}',\ldots,Y_{t-p}',\zeta_t)>0},Y_{t-1}',\ldots ,Y_{t-p+1}' \right) \nonumber \\
			 & =(Y_t',\ldots,Y_{t-p+1}') \nonumber \\
			 & = W_t.
			\end{align}
		Since $(Z_t)_{t\in\Z}$ is the only process satisfying this latter recursion, we have
			\begin{equation}
			\forall t \in \Z, \ W_t = Z_t,
			\end{equation}
		and in particular
			\begin{equation}
			\forall t \in \Z, \ Y_t=Y_t',
			\end{equation}
		which proves the uniqueness of process $(Y_t)_{t\in \Z}$.
		
		The process $\left( Y_t,\zeta_t\right)_{t\in \Z}$ is then ergodic, since we can write for all $t \in \Z$
			\begin{equation}
			(Y_t,\zeta_t) = \left( H(\zeta_t,\zeta_{t-1},\ldots ),\zeta_t \right).
			\end{equation}
		\end{enumerate}
	\item We already established in the first point that $\cap_{i=1}^k C_i^{j_i} \subset \left\lbrace \mathrm{Card}\left(F_1^p(E)=1\right) \right\rbrace$, thus
			\begin{equation}
			\E \left( \ind_{\bigcap_{i=1}^k C_i^{j_i}} \ \mid \ \mathcal{F}_0 \right) \leqslant \E \left( \ind_{\mathrm{Card}(F_1^p(E))=1} \ \mid \mathcal{F}_0 \right)
			\end{equation}
		i.e.
			\begin{equation}
			\P \left( \cap_{i=1}^k C_i^{j_i} \ \mid \mathcal{F}_0 \right) \leqslant \P \left( \mathrm{Card} \left( F_1^p(E) \right) = 1 \ \mid \mathcal{F}_0 \right),
			\end{equation}
		hence
			\begin{equation}
			\P \left( \mathrm{Card} \left( F_1^p(E) \right) = 1 \ \mid \mathcal{F}_0 \right) >0 \ \mathrm{a.s.}
			\end{equation}
		For simplicity of notations, we set $A_t =  \left\lbrace \mathrm{Card}(F_{t-p+1}^t(E)) =1 \right\rbrace$. We have
			\begin{equation}
			\ind_{A_t}=1 \iff \exists y_0 \in E, \ \forall y \in E, F_{\zeta_t} \circ \cdots \circ F_{\zeta_{t-p+1}}(y) = y_0.
			\end{equation}
		For any $y \in E$, the application
			\begin{equation}
			G_y : \begin{array}[t]{l}
			\left( \R^n \right)^p \longrightarrow E \\
			(e_1,\ldots,e_p) \longmapsto F_{e_p} \circ \cdots \circ F_{e_1} (y)
			\end{array}
			\end{equation}
		is measurable, we thus have
			\begin{align}
			\ind_{A_t} & = \sum_{y_0 \in E} \prod_{y\in E} \ind_{G_y^{-1}\left( \left\lbrace y_0 \right\rbrace \right)} \left( \zeta_t,\ldots,\zeta_{t-p+1} \right) \nonumber \\
			 & = H \left( \zeta_t^{-}\right)
			\end{align}
		for some measurable function $H : \left( \R^n \right)^{\N} \longrightarrow \left\lbrace 0,1 \right\rbrace$, where $\tau$ denotes the shift operator on the sequences.
		
		For $t \in \Z$ fixed, we have by Birkhoff's ergodic theorem
			\begin{equation}\label{eq_Birkhoff_N}
			\dfrac{1}{S} \sum_{s=0}^{S-1} \ind_{A_{t-s}} = \dfrac{1}{S} \sum_{s=0}^{S-1} H \left( \tau^s \zeta_t^{-} \right) \overset{S \to +\infty}{\longrightarrow} \P\left( A_t \ \mid \mathcal{I}_t \right) \quad \mathrm{a.s.}
			\end{equation}
		where $\mathcal{I}_t = \left\lbrace \left( \zeta_t^- \right)^{-1} (I) : I \in \mathcal{B}(\R^s)^{\otimes \N}, \ \tau^{-1}(I)=I \right\rbrace$ and $\zeta_t^-=\left( \zeta_{t-j} \right)_{j \geqslant 0}$.
		
		Yet, by Lemma \ref{lem_aux_birkhoff} below, we have $\mathcal{I}_t \subset \cap_{s \geqslant 0} \mathcal{F}_{t-s}$. In particular, $\mathcal{I}_t \subset \mathcal{F}_{t-p}$ and for all $J \in \mathcal{I}_t$
			\begin{equation}
			\P (A_t \cap J) = \E \left[ \P(A_t \mid \mathcal{I}_t)\ind_{J} \right] = \E \left[ \P (A_t \mid \mathcal{F}_{t-p}) \ind_{J} \right].
			\end{equation}
		The choice $J= \left\lbrace \P (A_t \mid \mathcal{I}_t) =0 \right\rbrace$ then gives
			\begin{equation}
			\E \left[ \P(A_t \mid \mathcal{I}_t)\ind_{J} \right] = 0,
			\end{equation}
		and so $\P(J)=0$ since $\P(A_t \mid \mathcal{F}_{t-p})>0$ almost surely.
		
		Indeed, for all $K \in \mathcal{F}_{t-p}$, $K = \left( \zeta_{t-p}^- \right)^{-1}(L)$, we have
			\begin{align}
			\E \left( \ind_{A_t} \ind_K \right) & = \E \left[ \P (A_t \mid \mathcal{F}_{t-p}) \ind_K \right] \nonumber \\
			 & = \E \left[ H\left( \zeta_t^- \right) \ind_L \left( \zeta_{t-p}^- \right) \right] \nonumber \\
			 & = \E \left[ H\left( \zeta_p^- \right) \ind_L \left( \zeta_0^- \right) \right] \quad \text{by stationarity} \nonumber \\
			 & =\E \left[ \P \left( A_p \mid \mathcal{F}_0 \right)\ind_L \left( \zeta_0^-\right) \right].
			\end{align}
		The choice $K = \left\lbrace \P \left(A_t \mid \mathcal{F}_{t-p} \right) = 0 \right\rbrace$ leads to
			\begin{equation}
			\E \left[ \P \left( A_p \mid \mathcal{F}_0 \right) \ind_L \left( \zeta_0^- \right) \right] = 0,
			\end{equation}
		and since we have already proved that $\P \left( A_p \mid \mathcal{F}_0 \right) > 0$ almost surely, we deduce
			\begin{equation}
			\P \left( \zeta_0^- \in L \right) = \P \left( \zeta_{t-p}^- \in L \right) = \P (K) = 0,
			\end{equation}
		so $\P \left( A_t \mid \mathcal{F}_{t-p}\right) >0$ almost surely.
		
		It follows that $\P \left( A_t \mid \mathcal{I}_t \right) >0$ almost surely, and by equation \eqref{eq_Birkhoff_N}
			\begin{equation}
			\P \left( \underset{s \to +\infty}{\limsup} \ A_{t-s} \right) =1.
			\end{equation}
		Moreover, if $\ind_{A_{t-s}}(\omega)=1$ for a certain value of $\omega$, we have for $T \geqslant p+s$
			\begin{align}
			F_{t-T}^t (y)(\omega) & = F_{t-s+1}^t \left( F_{t-s-p+1}^{t-s} \left( F_{t-T}^{t-s-p}(y) \right) \right)(\omega) \nonumber \\
			 & = F_{t-s-p+1}^t (y_0)(\omega)
			\end{align}
		where $y_0 \in E$ is arbitrary. Then $\lim_{n\rightarrow \infty} F_{t-n}^t(y)$ exists a.s. and does not depend on $y$.
		
		The end of the proof is similar to the one of Theorem 3 in \citet{truquet2023strong}.
	\end{enumerate}
	
	\begin{lem} \label{lem_aux_birkhoff}
	For all $s \geqslant 0$, $\mathcal{I}_t \subset \mathcal{F}_{t-s}.$
	\end{lem}
	
	\begin{proof}
	We actually show by induction on $j$ that
		\begin{equation}
		\forall j \in \N, \ \forall I \in \mathcal{B}(\R^s)^{\otimes \N}, \ \tau^{-1}(I)=I \Longrightarrow \left(  \zeta_t^-\right)^{-1}(I) = \left( \zeta_{t-j}^- \right)^{-1}(I).
		\end{equation}
		
		\begin{enumerate}[label=$\triangleright$,labelindent=0pt,labelwidth=!,wide,itemsep=-3pt]
		\item For $j=0$, there is nothing to prove.
		\item Assume the property true for a fixed $j \in \N$. We thus have for all $I \in \mathcal{B}(\R^s)^{\otimes \N}$ such that $\tau^{-1}(I)=I$
			\begin{align}
			 \left( \zeta_{t-j}^- \right)^{-1}(I)&=
			  \left\lbrace \omega : \left( \zeta_{t-j}(\omega),\zeta_{t-j-1}(\omega),\ldots \right) \in I \right\rbrace \nonumber \\
			 & = \left\lbrace \omega : \left( \zeta_{t-j}(\omega),\zeta_{t-j-1}(\omega),\ldots \right) \in \tau^{-1} I \right\rbrace \nonumber \\
			 & = \left\lbrace \omega : \left( \zeta_{t-j-1}(\omega),\zeta_{t-j-2}(\omega),\ldots \right) \in I \right\rbrace \nonumber \\
			 & = \left(\zeta_{t-j-1}^- \right)^{-1}(I),
			\end{align}
		which concludes the induction.
		\end{enumerate}
	\end{proof}
	
\subsection{Proof of Corollary \ref{cor_stationarity}}

Note first that $\sigma \left( X_s,\varepsilon_s \mid s \leqslant t-1 \right) = \sigma \left( \zeta_{t-1}^-, X_{t-1} \right)$.

The goal is to prove that assumption \textbf{A3} is satisfied. To this end, we set $j_1=\ldots=j_k=1$ and for $t\in \left\lbrace1,\ldots,p\right\rbrace$
	\begin{equation}
	B_t= \bigcap_{i=1}^k \left\lbrace \varepsilon_{i,t} > - \underset{y_1,\ldots,y_p}{\min} \ h_i (y_1,\ldots,y_p,X_{t-1}) \right\rbrace,
	\end{equation}
so that we have $\displaystyle{\cap_{i=1}^k C_i^{j_i} = \cap_{t=1}^p B_t}$.

By Assumption \textbf{A4}, we have 
	\begin{equation}
	\P \left( B_t \mid \sigma \left( \zeta_{t-1}^-,X_{t-1} \right) \right) >0 \ \mathrm{a.s.}
	\end{equation}
Indeed, for any $C \in \sigma \left( \zeta_{t-1}^-,X_{t-1} \right), \ C = \left( \zeta_{t-1}^-,X_{t-1} \right)^{-1}(B)$, we have
	\begin{align}
	\E \left( \ind_{B_t} \ind_C \right) & = \E \left[ \prod_{i=1}^k \ind_{\left\lbrace\varepsilon_{i,t} > -\underset{y_1,\ldots,y_p}{\min} \ h_i(y_1,\ldots,y_p,X_{t-1})\right\rbrace} \ind_C \right] \nonumber \\
	 & =\E \left( \ind_{K}(X_{t-1},\varepsilon_t)\ind_C\right)
	\end{align}
where $K= \left\lbrace (x,e) \in \R^d \times \R^k : \forall i \in \left\lbrace 1,\ldots,k\right\rbrace, \ x_i > - \underset{y_1,\ldots,y_p}{\min} \ h_i(y_1,\ldots,y_p,e) \right\rbrace$. Thus
	\begin{align}
	\E\left( \ind_{B_t} \ind_C \right) & = \int \ind_K \left( X_{t-1}(\omega),\varepsilon_t(\omega) \right) \ind_B \left( \zeta_{t-1}^-(\omega),X_{t-1}(\omega) \right) \P (\d \omega) \nonumber \\
	 & = \int \ind_K (x_{t-1},e) \ind_B (s^-,x_{t-1}) \P_{\left( \varepsilon_t,\zeta_{t-1}^-,X_{t-1} \right)} (\d e, \d s^-, \d x_{t-1}) \nonumber \\
	 & =\iint \ind_K (x_{t-1},e) \P_{\varepsilon_t \mid (\zeta_{t-1}^-,X_{t-1})}(\d e \mid s^-,x_{t-1}) \ind_B (s^-,x_{t-1}) \P_{(\zeta_{t-1}^-,X_{t-1})}(\d s^-,\d x_{t-1}) \nonumber
	 \shortintertext{where $ \P_{\varepsilon_t \mid (\zeta_{t-1}^-,X_{t-1})}$ is a conditional distribution of $\varepsilon_t$ given $(\zeta_{t-1}^-,X_{t-1})$ \nonumber}
	  & = \int \P_{\varepsilon_t \mid(\zeta_{t-1}^-,X_{t-1})}(K_{x_{t-1}} \mid s^-,x_{t-1}) \ind_B(s^-,x_{t-1}) \P_{(\zeta_{t-1}^-,X_{t-1})}(\d s^-,\d x_{t-1}) 
	  \end{align}
with
	\begin{equation}
	K_{x_{t-1}} = \left\lbrace e \in \R^k : \forall i \in \left\lbrace 1,\ldots,k\right\rbrace, \ e_i > - \underset{y_1,\ldots,y_p}{\min} h_i (y_1,\ldots,y_p,x_{t-1}) \right\rbrace = \cap_{i=1}^k \varphi_{i,x_{t-1}}^{-1} \left( (0,+\infty) \right)
	\end{equation}
and $\varphi_{i,x_{t-1}}(e) = e_i + \underset{y_1,\ldots,y_p}{\min} h_i(y_1,\ldots,y_p,x_{t-1})$ is trivially continuous. In particular, $K_{x_{t-1}}$ is a non-empty open set, and by assumption \textbf{A4} we have $\P_{\left( \zeta_{t-1}^-,X_{t-1} \right)}-$almost surely
	\begin{equation}
	\P_{\varepsilon_t \mid \left(\zeta_{t-1}^-,X_{t-1}\right)}(K_{x_{t-1}} \mid s^-,x_{t-1}) >0.
	\end{equation}
If we set $C = \left\lbrace \omega : \P \left( B_t \mid \sigma (\zeta_{t-1}^-,X_{t-1}) \right)(\omega) =0 \right\rbrace$, and suppose by absurd that $\P(C)>0$, we then get the contradiction
	\begin{equation}
	0 = \E (\ind_{B_t} \ind_C) =\int \P_{ \varepsilon_t \mid \left( \zeta_{t-1}^-,X_{t-1} \right)}(K_{x_{t-1}} \mid s^-,x_{t-1}) \ind_B(s^-,x_{t-1}) \P_{\left( \zeta_{t-1}^-,X_{t-1} \right)}(\d s^-,\d x_{t-1}) >0.
	\end{equation}
Thus $\P(C)=0$, i.e.
	\begin{equation}\label{eq_U1}
	\P \left( B_t \mid \sigma (\zeta_{t-1}^-,X_{t-1}) \right) >0 \quad \mathrm{a.s.}
	\end{equation}
We then set for $j \in \left\lbrace 1,\ldots ,p\right\rbrace$
	\begin{equation}
	U_j = \P \left( \cap_{t=1}^j B_t \mid \mathcal{F}_0 \right),
	\end{equation}
and we show by induction on $j$ that
	\begin{equation}
	\forall j \in \left\lbrace 1, \ldots,p\right\rbrace, \ U_j >0 \quad \mathrm{a.s.}
	\end{equation}
	\begin{enumerate}[label=$\triangleright$,labelindent=0pt,labelwidth=!,wide,itemsep=-3pt]
	\item We already have $U_1 >0$ almost surely by equation \eqref{eq_U1}.
	\item Assume that $U_j >0$ almost surely for a fixed $j \in \left\lbrace 1,\ldots,p\right\rbrace$, and take $C \in \mathcal{F}_0$, we have
		\begin{align}
		\E \left( \ind_{\cap_{t=1}^{j+1} B_j} \ind_C \right) & = \E \left( \ind_{B_{j+1}} \ind_{\cap_{j=1}^t B_j} \ind C \right) \nonumber \\
		 & =\E \left( \P(B_{j+1} \mid \mathcal{F}_j) \ind_{\cap_{j=1}^t B_j}\ind_C \right),
		\end{align}
	thus $U_{j+1} = \E \left( \P(B_{j+1} \mid \mathcal{F}_j) \ind_{\cap_{j=1}^t B_j} \mid \mathcal{F}_0  \right)$. However, $\P(B_{j+1} \mid \mathcal{F}_j)>0$ almost surely by equation \eqref{eq_U1}, and $\E \left( \ind_{\cap_{t=1}^j B_t} \mid \mathcal{F}_0 \right) >0$ almost surely by assumption. Thus we have $U_{j+1} >0$ almost surely, which concludes the induction.
	\end{enumerate}
We conclude that
	\begin{equation}
	U_p = \P \left( \cap_{t=1}^p B_t \mid \mathcal{F}_0 \right) = \P \left( \cap_{i=1}^k C_i^{j_i} \mid \mathcal{F}_0 \right) >0 \quad \mathrm{a.s}
	\end{equation}
and assumption \textbf{A3} is verified.
	
\subsection{Proof of Theorem \ref{theo_cons_single_global}}

\subsubsection{Proof for estimator $\hat{\theta}_{\PL}$}

The consistency, the stochastic expansion, and the asymptotic normality of the estimator $\hat{\theta}$ is obtained by using Theorem \ref{theo_pfanzagl} and Theorem \ref{theo_straumann}. We thus have to check the assumptions \emph{\textbf{H1}} to \emph{\textbf{H9}} mentioned in these latter theorems for the process $(Z_t)_{t\in \Z}$ defined by
	\begin{equation}
	Z_t = \left(Y_t,\ldots , Y_{t-p},X_{t-1} \right),
	\end{equation}
and the family of applications $(m_{\theta})_{\theta \in \Theta}$ defined by
	\begin{equation}\label{eq_m_theta_single}
	m_{\theta}(Z_t)= \log \left( \int_{\R^k} \prod_{i=1}^k \ind_{I_Y{i,t}}(\lambda_{i,t}+e_i) \varphi_R(e) \d e \right).
	\end{equation}

	\begin{enumerate}[label=$\triangleright$,labelindent=0pt,labelwidth=!,wide]
	\item \emph{\textbf{H1}} is satisfied since $(Y_t)_{t\in \Z}$ and $(X_t)_{t\in \Z}$are both ergodic by assumption \textbf{B1} and corollary \ref{cor_stationarity}.
	\item \emph{\textbf{H2}} is satisfied by assumption \textbf{B4}.
	\item In order to check assumption \textbf{\emph{H3}}, we first prove the continuity of the application $\theta \longmapsto m_{\theta}(Z_0)$ for any fixed value of $Z_0$. Here we have
		\begin{equation}
		m_{\theta}(Z_0) = \log (f_{\theta}(Z_0))
		\end{equation}
	where
		\begin{equation}
		f_{\theta} (Z_0) = \int_{\R^k} \prod_{i=1}^k \ind_{I_{Y_{i,0}}-\lambda_{i,0}(\gamma)}(x_i) \varphi_R (x) \d x.
		\end{equation}
	Let us fix $\theta_0 = (\gamma_0,R_0) \in \Theta$, for any $i \in \left\lbrace 1, \ldots, k \right\rbrace$, we denote
		\begin{equation}
		H_i = \left\lbrace x=(x_1,\ldots,x_k) \in \R^k : x_i = -\lambda_{i,0}(\gamma_0) \right\rbrace
		\end{equation}
	and $H = \cup_{i=1}^k H_i$. $H$ is negligible for the Lebesgue's measure, since it is a finite union of hyperplanes, and it is easy to verify that 
		\begin{equation}
		\prod_{i=1}^k \ind_{I_{Y_{i,0}}-\lambda_{i,0}(\gamma)}(x_i) \varphi_R (x) \underset{\theta \to \theta_0}{\longrightarrow} \prod_{i=1}^k \ind_{I_{Y_{i,0}}-\lambda_{i,0}(\gamma_0)}(x_i) \varphi_{R_0} (x)
		\end{equation}
	for any $x \notin H$.
	
	Furthermore, according to Lemma \ref{lem_Phi_R_bounds} below, it is possible to find positive constants $m_{\det}$ and $c_2$ such that for any $\theta \in \Theta$ and any $x \in \R^k$
		\begin{equation}
		\prod_{i=1}^k \ind_{I_{Y_{i,0}}-\lambda_{i,0}(\gamma)}(x_i) \varphi_R (x) \leqslant \dfrac{1}{\sqrt{(2\pi)^k m_{\det}}}\exp \left( -\frac{c_2}{2} \|x\|_2^2 \right),
		\end{equation}
	which is integrable on $\R^k$. Thus, by continuity under the integral sign, the application $\theta \longmapsto f_{\theta}(Z_0)$ is continuous, and it follows that $\theta \longmapsto m_{\theta}(Z_0)$ is also continuous.
	
	In order prove that $\E \left( \underset{\theta \in \Theta}{\sup} \ \left| m_{\theta}(Z_0) \right| \right) < +\infty$, we use once again Lemma \ref{lem_Phi_R_bounds} to obtain positive constants $M_{\det}$ an $c_1$ such that for any $\theta \in \Theta$ and $x \in \R^k$
		\begin{align}
		f_{\theta}(Z_0) & \geqslant \dfrac{1}{\sqrt{(2\pi)^k M_{\det}}} \int_{\R^k} \prod_{i=1}^k \ind_{I_{Y_{i,0}}-\lambda_{i,0}(\gamma)}(x_i) \exp \left( - \frac{c_1}{2} \|x\|_2^2 \right) \ d x \nonumber \\
		& \geqslant \dfrac{1}{\sqrt{c_1 M_{\det}}} \prod_{i=1}^k \min \left\lbrace \Phi \left(-\sqrt{c_1} \lambda_{i,0}(\gamma) \right), 1- \Phi \left(-\sqrt{c_1} \lambda_{i,0}(\gamma) \right)\right\rbrace \nonumber \\
		\shortintertext{with the substitution $u=\sqrt{c_1}x$} \nonumber
		& \geqslant \dfrac{1}{\sqrt{c_1 M_{\det}}} \exp \left(-dc_1 \| \lambda_0(\gamma) \|_2^2 \right)
		\end{align}
	where $d$ is the positive constant mentioned in Lemma \ref{lem_exp_ineq}. For any $\theta \in \Theta$, we thus have
		\begin{align}
		\left| m_{\theta}(Z_0) \right| & = - \log \left( f_{\theta}(Z_0) \right) \nonumber \\
		& \leqslant \dfrac{1}{2}\log(c_1M_{\det}) + dc_1 \|\lambda_0(\gamma)\|_2^2 \nonumber \\
		& \leqslant \dfrac{1}{2}\log(c_1M_{\det}) + dc_1 \left( M_{\mathcal{A}} + M_{\mathcal{B}}\|X_{-1}\|_2 \right)^2
		\end{align}
	for some positive constants $M_{\mathcal{A}}$ and $M_{\mathcal{B}}$ obtained by compactness of $\Theta$. Since it is assumed that $X_{-1} \in L^2$, this latter quantity is integrable and we have
		\begin{equation}
		\E \left( \underset{\theta \in \Theta}{\sup} \ \left| m_{\theta}(Z_0) \right| \right) < +\infty.
		\end{equation}
	\begin{lem}\label{lem_Phi_R_bounds}
		There exist constants $M_{\mathrm{det}},m_{\mathrm{det}},c_1,c_2 >0$ such that for all $R \in \mathcal{R}$ and all $x \in \R^k$
			\begin{equation}
			\dfrac{1}{\sqrt{(2\pi)^k M_{\mathrm{det}}}} \times \exp \left( -\dfrac{c_1}{2}\| x\|_2^2 \right) \leqslant \varphi_R (x) \leqslant \dfrac{1}{\sqrt{(2\pi)^k m_{\mathrm{det}}}} \times \exp \left( -\dfrac{c_2}{2}\| x\|_2^2 \right).
			\end{equation}
		\end{lem}
		
		\begin{lemproof}[ of Lemma \ref{lem_Phi_R_bounds}]
		By continuity of the determinant, and by compactness of $\mathcal{R}$, there exist $m_{\mathrm{det}}, M_{\mathrm{det}} >0$ such that
			\begin{equation}
			\forall R \in \mathcal{R}, \ m_{\mathrm{det}} \leqslant \det (R) \leqslant M_{\det}.
			\end{equation}
		Furthermore, take any $R \in \mathcal{R}$, we have for all $x \in \R^k$
			\begin{align}\label{eq_def_Phi_R}
			\varphi_R(x) & = \dfrac{1}{\sqrt{(2\pi)^k \det(R)}}\exp \left(-\frac{1}{2} \langle R^{-1} x , x \rangle \right) \nonumber \\
			& = \dfrac{1}{\sqrt{(2\pi)^k \det(R)}}\exp \left(-\frac{1}{2} \| x\|_{R^{-1}}^2 \right)
			\end{align}
		where $\| \cdot \|_{R^{-1}}$ denotes the norm associated with the matrix $R^{-1}$. By equivalence of the norms, there exist $c(R), C(R) >0$ such that
			\begin{equation}
			\forall x \in \R^k, \ c(R)\times  \|x\|_2 \leqslant \|x \|_{R^{-1}} \leqslant C(R) \times \| x\|_2.
			\end{equation}
		Without loss of generality, we can suppose that
			\begin{equation}
			c(R) = \sup \ \left\lbrace c>0 : \forall x \in \R^k, \ c\|x\|_2 \leqslant \|x \|_{R^{-1}} \right\rbrace
			\end{equation}
		and
			\begin{equation}
			C(R) = \inf \ \left\lbrace C > 0 : \forall x \in \R^k, \ C\|x\|_2 \geqslant \|x \|_{R^{-1}} \right\rbrace.
			\end{equation}
		Let us show that there exists $c>0$ such that
			\begin{equation}
			\forall R \in \mathcal{R}, \ c \leqslant c(R).
			\end{equation}
		By absurd, assume that for all integer $n \geqslant 1$, there exists $R_n \in \mathcal{R}$ such that
			\begin{equation}
			c(R_n) < \dfrac{1}{n}.
			\end{equation}
		Since $\mathcal{R}$ is compact, it is possible to extract a subsequence $\left( R_{\varphi(n)} \right)_{n \geqslant 1}$ that converges toward $\tilde{R} \in \mathcal{R}$. Yet, for all $n \geqslant 1$, we have $c(R_{\varphi(n)}) < \dfrac{1}{\varphi(n)}$, thus there exists $x_{\varphi(n)} \in \mathcal{S}^1$ such that
			\begin{equation}
			\| x_{\varphi(n)} \|_{R^{-1}} < \dfrac{1}{\varphi(n)}.
			\end{equation}
		By compactness once again, we can extract a subsequence $\left( x_{\varphi \circ \psi (n)} \right)_{n \geqslant 1 }$ that converges toward $\tilde{x} \in \mathcal{S}^1$. Thus, for all $ n\geqslant 1$
			\begin{equation}
			\| x_{\varphi \circ \psi (n)} \|_{R_{\varphi \circ \psi(n)}^{-1}} < \dfrac{1}{\varphi \circ \psi (n)} \underset{n \to +\infty}{\longrightarrow} 0.
			\end{equation}
		However, by continuity of the inverse, we also have $\| x_{\varphi \circ \psi (n)} \|_{R_{\varphi \circ \psi(n)}^{-1}} \underset{n \to +\infty}{\longrightarrow} \| \tilde{x} \|_{\tilde{R}^{-1}}$. Hence $\tilde{x} =0$, which is a contradiction.
		
		Therefore, there exists indeed a constant $c >0$ such that
			\begin{equation}
			\forall R \in \mathcal{R}, \ c \leqslant c(R).
			\end{equation}
		One can show the same way that there exists a constant $C >0$ such that
			\begin{equation}
			\forall R \in \mathcal{R}, \ C(R) \leqslant C.
			\end{equation}
		Combining with equality \eqref{eq_def_Phi_R}, we finally have
			\begin{equation}
			\dfrac{1}{\sqrt{(2\pi)^k M_{\det}}} \times \exp \left( -\frac{C^2}{2} \| x\|_2^2 \right) \leqslant \varphi_R (x) \leqslant \dfrac{1}{\sqrt{(2\pi)^k m_{\det}}} \times \exp \left( -\frac{c^2}{2} \| x\|_2^2 \right),
			\end{equation}
		which is the desired result with $c_1 = C^2$ and $c_2=c^2$.
		\end{lemproof}
	\item We now check {\bf H4}. For $\theta \in \Theta$ fixed, setting for any $t \in \Z, \ Y_t^- =(Y_t,Y_{t-1},\ldots )$, we have
		\begin{equation}
		\E \left( m_{\theta} (Z_0) \right) - \E \left( m_{\theta_0} (Z_0) \right) = \E \left( \log \left( \dfrac{\psi_{\theta}(Y_0,Y_{-1}^-,X_{-1})}{\psi_{\theta_0}(Y_0,Y_{-1}^-,X_{-1})} \right) \right),
		\end{equation}
	where by assumption \textbf{B2}, $\displaystyle{\psi_{\theta}(y_0,Y_{-1}^-,X_{-1}) = \int_{\R^k} \prod_{i=1}^k \ind_{I_{y_{i,0}}}(\lambda_{i,0}(\gamma)+x_i) \varphi_R(x) \d x}$ is the density of the conditional distribution $\mathcal{L}_{\theta}\left( Y_0 \mid Y_{-1}^-,X_{-1} \right)$ with respect to the measure $\displaystyle{\mu = \sum_{y\in \left\lbrace 0,1 \right\rbrace^k} \delta_y}$ on $\left\lbrace 0, 1 \right\rbrace^k$.
		
	Using the fact that
		\begin{equation}
		\forall x >0, \ \log (x) \leqslant 2\left( \sqrt{x} -1 \right),
		\end{equation}
	with equality if and only if $x=1$, we thus have
		\begin{small}
		\begin{align}
		\E \left( m_{\theta} (Z_0) \right) - \E \left( m_{\theta_0} (Z_0) \right) & = \iint \log \left( \dfrac{\psi_{\theta}(y_0,y_{-1}^-,x_{-1})}{\psi_{\theta_0}(y_0,y_{-1}^-,x_{-1})} \right) \psi_{\theta_0}(y_0,y_{-1}^-,x_{-1}) \ \mu (\d y_0) \P_{(Y_{-1}^-,X_{-1})}(\d y_{-1}^-, \d x_{-1}) \nonumber \\
		& \leqslant \iint 2 \sqrt{\psi_{\theta}(y_0,y_{-1}^-,x_{-1})}\times \sqrt{\psi_{\theta_0}(y_0,y_{-1}^-,x_{-1})} \ \mu (\d y_0) -2 \P_{(Y_{-1}^-,X_{-1})}(\d y_{-1}^-, \d x_{-1}) \nonumber \\
		& = - \iint \left( \sqrt{\psi_{\theta}(y_0,y_{-1}^-,x_{-1})}- \sqrt{\psi_{\theta_0}(y_0,y_{-1}^-,x_{-1})} \right)^2 \ \mu (\d y_0) \P_{(Y_{-1}^-,X_{-1})}(\d y_{-1}^-, \d x_{-1}) \nonumber \\
		& \leqslant 0,
		\end{align}
		\end{small}
	hence
		\begin{equation}
		\forall \theta \in \Theta, \ \E \left( m_{\theta} (Z_0) \right) \leqslant \E \left( m_{\theta_0} (Z_0) \right).
		\end{equation}
	Furthermore, the equality case implies that $\P_{(Y_{-1}^-,X_{-1})}$-almost surely
		\begin{equation}
		\int \left[ \log \left( \dfrac{\psi_{\theta}(y_0,y_{-1}^-,x_{-1})}{\psi_{\theta_0}(y_0,y_{-1}^-,x_{-1})} \right) -2 \left( \sqrt{\dfrac{\psi_{\theta}(y_0,y_{-1}^-,x_{-1})}{\psi_{\theta_0}(y_0,y_{-1}^-,x_{-1})}} -1 \right) \right] \psi_{\theta_0}(y_0,y_{-1}^-,x_{-1}) \mu (\d y_0) =0,
		\end{equation}
	thus, $\P_{\left(Y_{-1}^-,X_{-1} \right)}$-almost surely, we have for any $y_0 \in \left\lbrace 0, 1 \right\rbrace^k$
		\begin{equation}
		\left[ \log \left( \dfrac{\psi_{\theta}(y_0,y_{-1}^-,x_{-1})}{\psi_{\theta_0}(y_0,y_{-1}^-,x_{-1})} \right) -2 \left( \sqrt{\dfrac{\psi_{\theta}(y_0,y_{-1}^-,x_{-1})}{\psi_{\theta_0}(y_0,y_{-1}^-,x_{-1})}} -1 \right) \right] \psi_{\theta_0}(y_0,y_{-1}^-,x_{-1}) = 0.
		\end{equation}
	Since $\psi_{\theta_0}(y_0,y_{-1}^-,x_{-1})$ is necessarily positive, whatever the values of $y_0,y_{-1}^-$ and $x_{-1}$, we have $\P_{(Y_{-1}^-,X_{-1})}$-almost surely
		\begin{equation}
		\forall y_0 \in \left\lbrace 0, 1\right\rbrace^k, \ \log \left( \dfrac{\psi_{\theta}(y_0,y_{-1}^-,x_{-1})}{\psi_{\theta_0}(y_0,y_{-1}^-,x_{-1})} \right) -2 \left( \sqrt{\dfrac{\psi_{\theta}(y_0,y_{-1}^-,x_{-1})}{\psi_{\theta_0}(y_0,y_{-1}^-,x_{-1})}} -1 \right) = 0,
		\end{equation}
	which leads to
		\begin{equation}
		\forall y_0 \in \left\lbrace 0, 1 \right\rbrace^k , \ \psi_{\theta}(y_0,y_{-1}^-,x_{-1}) = \psi_{\theta_0}(y_0,y_{-1}^-,x_{-1})
		\end{equation}
	$\P_{(Y_{-1}^-,X_{-1})}$-almost surely. Equivalently, this latter equation means that the conditional distributions $\mathcal{L}_{\theta}\left( Y_0,\mid Y_{-1}^-,X_{-1}\right)$ and $\mathcal{L}_{\theta_0}\left( Y_0,\mid Y_{-1}^-,X_{-1}\right)$ are almost surely identical.
		
	In particular, we have almost surely
		\begin{equation}
		\forall i \in \left\lbrace 1, \ldots , k\right\rbrace, \ \forall y \in \left\lbrace 0 ,1 \right\rbrace, \ \P_{\theta} \left( Y_{i,0}=y \mid Y_{-1}^-,X_{-1} \right) = \P_{\theta_0} \left( Y_{i,0}=y \mid Y_{-1}^-,X_{-1} \right),
		\end{equation}
	thus we have almost surely
		\begin{align}
		 & \forall i \in \left\lbrace 1, \ldots k \right\rbrace, \ \P_{\theta} \left( Y_{i,0}=1 \mid Y_{-1}^-,X_{-1} \right) = \P_{\theta_0} \left( Y_{i,0}=1 \mid Y_{-1}^-,X_{-1} \right) \nonumber \\
	\mathrm{i.e.} \quad  & \forall i \in \left\lbrace 1, \ldots k \right\rbrace, \ \Phi (-\lambda_{i,0}(\gamma)) = \Phi(-\lambda_{i,0}(\gamma_0)) \nonumber \\
		\mathrm{i.e.} \quad  & \forall i \in \left\lbrace 1, \ldots k \right\rbrace, \ \lambda_{i,0}(\gamma) = \lambda_{i,0}(\gamma_0).
		\end{align}	
	If we denote $A_1^{(0)},\ldots A_p^{(0)}, B^{(0)}$ and $R^{(0)}$ the values of the matrices composing $\theta_0$, this latter equation leads to
		\begin{equation}\label{eq_identification_reg}
		\sum_{l=1}^p \underbrace{\left(A_l-A_l^{(0)} \right)}_{M_l} Y_{-l} + \underbrace{\left(B-B^{(0)}\right)}_{N} X_{-1}=0 \quad \mathrm{a.s.}
		\end{equation}		
	By assumption \textbf{B5}, we immediately have $N=0$.
		
	Furthermore, since for any $y_{-1},\ldots ,y_{-p}$ in $\left\lbrace 0,1 \right\rbrace^{k}$,
		\begin{equation}
		\P \left( \left\lbrace (Y_{-1},\ldots , Y_{-p}) = (y_{-1},\ldots , y_{-p}) \right\rbrace \right) > 0,
		\end{equation}
	we have for any $y_{-1}, \ldots y_{-p}$ in $ \left\lbrace 0,1 \right\rbrace^{k}$
		\begin{equation}
		M_1 \cdot y_{-1} = -\sum_{l=2}^p M_l \cdot y_{-l}.
		\end{equation}
	This equality is true in the particular case where $y_{-2}=\ldots=y_{-p}=0$, and then we have
		\begin{equation}
		\forall y_{-1} \in \left\lbrace 0,1 \right\rbrace^k, \ M_1 \cdot y_{-1}=0.
		\end{equation}
	This latter equality is thus true for any vector in the canonical base of $\R^k$, and it follows that $M_1=0$.
		
	A simple induction then leads to
		\begin{equation}
		\forall l \in \left\lbrace 1, \ldots , p \right\rbrace, \ M_l = 0.
		\end{equation}
		
		\begin{note}
		In the specific case where the processes $(X_t)_{t\in \Z}$ and $(\varepsilon_t)_{t \in \Z}$ are supposed independent, and assumption \textbf{B5} is replaced by assumption \textbf{B5'}, the previous proof can be modified from equality \eqref{eq_identification_reg}.
		
		In this case, let us suppose by absurd that $M_1 \neq 0$. Then, there exists $v \in \R^k \setminus \left\lbrace 0 \right\rbrace$ such that
			\begin{align}
			1 & = \P \left(v'\cdot Y_{-1} = G(X_{-1}^-,\varepsilon_{-2}^-) \right) \nonumber \\
			\shortintertext{where $G$ is a measurable application} \nonumber
			& = \int \int \ind_{\left\lbrace G(x_{-1}^-, e_{-2}^-) \right\rbrace}(v'\cdot y_{-1}) \P_{Y_{-1} \mid (X_{-1}^-,\varepsilon_{-2}^-)=(x_{-1}^-,e_{-2}^-)}(\d y_{-1}) \P_{(X_{-1}^-,\varepsilon_{-2}^-)}(\d x_{-1}^- , \d e_{-2}^-).
			\end{align}
		Thus, for $\P_{(X_{-1}^-,\varepsilon_{-2}^-)}$-almost all $(x_{-1}^-,e_{-2}^-)$ we have
			\begin{equation}
			\int \ind_{\left\lbrace G(x_{-1}^-, e_{-2}^-) \right\rbrace}(v'\cdot y_{-1}) \P_{Y_{-1} \mid (X_{-1}^-,\varepsilon_{-2}^-)=(x_{-1}^-,e_{-2}^-)}(\d y_{-1}) = 1.
			\end{equation}
		There exists in particular an element $(x_{-1}^-,e_{-2}^-)$ such that for $\P_{Y_{-1} \mid (X_{-1}^-,\varepsilon_{-2}^-)=(x_{-1}^-,e_{-2}^-)}$-almost all $y \in \left\lbrace 0,1 \right\rbrace^k$,
			\begin{equation}\label{eq_identification_H}
			v'\cdot y = G(x_{-1}^-,e_{-2}^-).
			\end{equation}
		However, the conditional distribution $\mathcal{L}_{\theta}(Y_{-1} \mid X_{-1}^-,\varepsilon_{-2}^-)$ is fully supported by $\left\lbrace 0,1 \right\rbrace^k$. Denoting $H(X_{-1}^-,\varepsilon_{-2}^-) = \sum_{l=1}^p A_l Y_{-l-1} + BX_{-2}$, we have indeed for any adequate measurable set $C$ and any $y_{-1} \in \left\lbrace 0,1 \right\rbrace^k$
			\begin{align}
			\P(Y_{-1} = y_{-1} , (X_{-1}^-, \varepsilon_{-2}^-) \in C) & = \P \left( \varepsilon_{-1} \in \prod_{i=1}^k I_{y_{i,-1}-H(X_{-1}^-,\varepsilon_{-2}^-)}, (X_{-1}^-, \varepsilon_{-2}^-) \in C)  \right) \nonumber \\
			& = \int \int_{\R^k} \prod_{i=1}^k \ind_{I_{y_{i,-1}}-H(x_{-1}^-,e_{-2}^-)}(e_{-1}) \P_{\varepsilon_{-1} \mid \left(X_{-1}^-, \varepsilon_{-2}^-\right)=(x_{-1}^-,e_{-2}^-)} (\d e_{i,-1}) \nonumber \\
			& \hspace*{3cm} \times \ind_{C}(x_{-1}^-,e_{-2}^-) \P_{\left( X_{-1}^-,\varepsilon_{-2}^- \right)}(\d x_{-1}^-,\d e_{-2}^-) \nonumber \\
			& = \int \int_{\R^k} \prod_{i=1}^k \ind_{I_{y_{i,-1}}-H(x_{-1}^-,e_{-2}^-)}(e_{i,-1}) \P_{\varepsilon_{-1}}(\d e_{-1}) \nonumber \\
			& \hspace{3cm} \times \ind_{C}(x_{-1}^-,e_{-2}^-) \P_{\left( X_{-1}^-,\varepsilon_{-2}^- \right)}(\d x_{-1}^-,\d e_{-2}^-) \\
			\end{align}
		since $\varepsilon_{-1}$ is independent from $\left( X_{-1}^-,\varepsilon_{-2}^- \right)$ from the assumption of independence between $(X_{t})_{t \in Z}$ and $(\varepsilon_t)_{t\in \Z}$ and assumptions \textbf{B2}. Thus
			\begin{equation}
			\P\left(Y_{-1} = y_{-1} \mid \left(X_{-1}^-,\varepsilon_{-2}^-\right) = \left(x_{-1}^-,e_{-2}^-\right) \right) = \int_{\R^k} \prod_{i=1}^k \ind_{I_{y_{i,-1}}-G(x_{-1}^-,e_{-2}^-)}(e_{i,-1}) \P_{\varepsilon_{-1}}(\d e_{-1}) >0.
			\end{equation}
		Choosing $y=0$ in \eqref{eq_identification_H}, we obtain $G(x_{-1}^-,e_{-2}^-)=0$, and it follows that for any $y \in \left\lbrace 0,1\right\rbrace^k$
			\begin{equation}
			v'\cdot y =0,
			\end{equation}
		which is impossible since $v \neq 0$. Thus, $M_1 =0$, and an immediate induction leads to
			\begin{equation}
			\forall l \in \left\lbrace1, \ldots,p\right\rbrace, \ M_{l}=0.
			\end{equation}
		Equation \eqref{eq_identification_reg} then becomes
			\begin{equation}
			N\cdot X_{-1} = 0,
			\end{equation}
		and assumption \textbf{B5'} gives $N=0$.
		\end{note}
		
	In the end, we have shown that the equality $\E \left( m_{\theta}(Z_0) \right) = \E \left( m_{\theta_0}(Z_0) \right)$ implies $B=B^{(0)}$ and $A_l^{(0)}=A_l$ for any $l \in \left\lbrace 1 , \ldots , p \right\rbrace$.
		
	We still have to prove that this latter equality also implies that $R=R^{(0)}$. We already mentioned that the conditional distributions $\mathcal{L}_{\theta} (Y_0 \mid Y_{-1}^-,X_{-1})$ and $\mathcal{L}_{\theta_0}(Y_0 \mid Y_{-1}^-,X_{-1})$ are almost surely identical. Therefore, it is also the case for the distributions of the couples, i.e. for any $i<j$ in $\left\lbrace 1, \ldots , k\right\rbrace$
		\begin{equation}
		\mathcal{L}_{\theta} \left( (Y_{i,0},Y_{j,0}) \mid Y_{-1}^-, X_{-1} \right) \quad \text{and} \quad \mathcal{L}_{\theta_0} \left( (Y_{i,0},Y_{j,0}) \mid Y_{-1}^-, X_{-1} \right)
		\end{equation}
	are almost surely identical. Note that if $k=1$, there is nothing to prove since the covariance matrix $R$ is actually a correlation matrix.
		
	In particular, we have almost surely
		\begin{equation}
			\P_{\theta} \left( Y_{i,0}=0 , Y_{j,0}=0 \mid Y_{-1}^- , X_{-1}  \right) = \P_{\theta_0} \left( Y_{i,0}=0 , Y_{j,0}=0 \mid Y_{-1}^- , X_{-1}  \right)
		 \end{equation}
	i.e.
		\begin{equation}
		\int_{-\infty}^{-\lambda_{i,0}(\gamma)} \Phi \left( \dfrac{-\lambda_{j,0}(\gamma)-R(i,j)x_i}{\sqrt{1-R(i,j)^2}} \right) \varphi (x_i) \d x_i = \int_{-\infty}^{-\lambda_{i,0}(\gamma_0)} \Phi \left( \dfrac{-\lambda_{j,0}(\gamma_0)-R^{(0)}(i,j)x_i}{\sqrt{1-R^{(0)}(i,j)^2}} \right) \varphi (x_i) \d x_i.
		\end{equation}
	Since we have already proved that $\lambda_{0}(\gamma) = \lambda_{0}(\gamma_0)$, we thus have almost surely
		\begin{equation} \label{eq_ident_r_single}
		\int_{-\infty}^{-\lambda_{i,0}(\gamma_0)} \Phi \left( \dfrac{-\lambda_{j,0}(\gamma_0)-R(i,j)x_i}{\sqrt{1-R(i,j)^2}} \right) \varphi (x_i) \d x_i = \int_{-\infty}^{-\lambda_{i,0}(\gamma_0)} \Phi \left( \dfrac{-\lambda_{j,0}(\gamma_0)-R^{(0)}(i,j)x_i}{\sqrt{1-R^{(0)}(i,j)^2}} \right) \varphi (x_i) \d x_i.
		\end{equation}
	However, for any $\omega_i, \ \omega_j$ in $\R$, the mapping
		\begin{equation}
			f : \begin{array}[t]{l}
			(-1,1) \longrightarrow \R \\
			r \longmapsto \displaystyle{ \int_{-\infty}^{\omega_i} \Phi \left( \dfrac{\omega_j-rx_i}{\sqrt{1-r^2}} \right) \varphi (x_i) \d x_i}
			\end{array}
		\end{equation}
	is an increasing function (see Lemma \ref{lem_increasing_f}). We then deduce that $R(i,j)=R^{(0)}(i,j)$ for any $i<j$ in $\left\lbrace 1, \ldots k \right\rbrace$, thus $R=R^{(0)}$.
		
	As a consequence, we do have
		\begin{equation}
		\forall \theta \in \Theta, \ \E \left( m_{\theta}(Z_0) \right) \leqslant \E \left( m_{\theta_0}(Z_0) \right),
		\end{equation}
	with equality if and only if $\theta=\theta_0$, so \emph{\textbf{H4}} is satisfied.
	
	\begin{lem}\label{lem_increasing_f}
	Take any $\omega_i, \ \omega_j$ in $\R$, the function
		\begin{equation}
		f : \begin{array}[t]{l}
		(-1,1) \longrightarrow \R \\
		r \longmapsto \displaystyle{\int_{-\infty}^{\omega_i}} \Phi \left( \dfrac{\omega_j-r e_i}{\sqrt{1-r^2}} \right) \varphi (e_i) \d e_i
		\end{array}
		\end{equation}
	is strictly increasing.
	\end{lem}
	
	\begin{lemproof}[ of Lemma \ref{lem_increasing_f}]
	We first apply the theorem of differentiation of a function defined by an integral to $f$. To do so, we verify the following assertions.
		\begin{enumerate}[label=$\triangleright$]
		\item For any $r \in (-1,1)$ and any $x_i \in (-\infty , \omega_i]$, we have
			\begin{equation}
			\left| \Phi \left( \dfrac{\omega_j-r x_i}{\sqrt{1-r^2}} \right) \varphi (x_i) \right| \leqslant \varphi(x_i),
			\end{equation}
		thus the integrand function defining $f$ is integrable on $]-\infty , \omega_i]$.
		\item For any $x_i \in (-\infty , \omega_i]$, the application
			\begin{equation}
			r \longmapsto \Phi \left( \dfrac{\omega_j-r x_i}{\sqrt{1-r^2}} \right) \varphi (x_i)
			\end{equation}
		is differentiable, and its derivative is equal to
			\begin{equation}
			\dfrac{r\omega_j - x_i}{\left(1-r^2 \right)^{3/2}} \times \varphi \left( \dfrac{\omega_j-r x_i}{\sqrt{1-r^2}} \right) \varphi(x_i).
			\end{equation}
		Without loss of generality, we can take $r$ in a compact subset of $(-1,1)$, and this latter function can be bounded by an integrable function.
		\end{enumerate}
	The function $f$ is thus differentiable on $(-1,1)$ and we have
		\begin{equation}
		\forall r \in (-1,1), \ \dot{f}(r) = \int_{-\infty}^{\omega_i} \dfrac{r\omega_j - x_i}{\left(1-r^2 \right)^{3/2}} \times \varphi \left( \dfrac{\omega_j-r x_i}{\sqrt{1-r^2}} \right) \varphi(x_i) \d x_i.
		\end{equation}
	Denoting $\alpha = \dfrac{\omega_j}{\sqrt{1-r^2}}, \ \beta = \dfrac{r}{\sqrt{1-r^2}}, \ \gamma = \dfrac{r\alpha}{1-r^2}$ and $\delta = \dfrac{1}{(1-r^2)^{3/2}}$, we get
		\begin{equation}
		\dot{f}(r) = \gamma \underbrace{\int_{-\infty}^{\omega_i} \varphi (\alpha - \beta x_i) \varphi (x_i) \d x_i}_{J} - \delta \underbrace{\int_{-\infty}^{\omega_i} x_i \varphi (\alpha - \beta x_i) \varphi (x_i) \d x_i}_{I}.
		\end{equation}
	Yet
		\begin{align}
		I & = -\int_{-\infty}^{\omega_i} \varphi (\alpha - \beta x_i) \dot{\varphi}(x_i) \d x_i \nonumber \\
		 & = - \left( \left[ \varphi (\alpha - \beta x_i)\varphi(x_i) \right]_{-\infty}^{\omega_i} + \beta \int_{-\infty}^{\omega_i} \dot{\varphi}(\alpha-\beta x_i) \varphi (x_i) \d x_i \right) \nonumber \\ 
		 & = - \varphi (\alpha - \beta \omega_i) \varphi(\omega_i) + \beta \int_{-\infty}^{\omega_i} (\alpha - \beta x_i) \varphi (\alpha-\beta x_i) \varphi (x_i) \d x_i,
		\end{align}
	hence
		\begin{equation}
		I =-\varphi(\omega_i)\varphi(\alpha-\beta \omega_i) - \beta^2I + \alpha \beta J 
		\end{equation}
	i.e.
		\begin{equation}
		I = \dfrac{-\varphi(\omega_i) \varphi(\alpha-\beta \omega_i)}{1+\beta^2} + \dfrac{\alpha \beta}{1+\beta^2} J.
		\end{equation}
	We thus have
		\begin{align}
		\dot{f}(r) & = \gamma J + \dfrac{\delta \varphi(\omega_i)\varphi(\alpha-\beta\omega_i)}{1+\beta^2}- \dfrac{\alpha \beta}{1+\beta^2}J \nonumber \\
		 & = \dfrac{\delta \varphi(\omega_i)\varphi(\alpha-\beta\omega_i)}{1+\beta^2} + \left( \gamma - \dfrac{\alpha \beta \delta}{1+\beta^2} \right) J,
		\end{align}
	but in the other hand $\dfrac{\delta}{1+\beta^2}= \dfrac{1}{\sqrt{1-r^2}}$ and
		\begin{equation}
		\gamma - \dfrac{\alpha \beta \delta}{1+\beta^2}  = \dfrac{r\alpha}{1-r^2}- \dfrac{\alpha r}{1-r^2} = 0.
		\end{equation}
	Thus $\dot{f}(r) = \dfrac{\varphi(\omega_i)\varphi(\alpha-\beta \omega_i)}{\sqrt{1-r^2}}>0$ and $f$ is strictly increasing.
	\end{lemproof}
	
	\item Assumption \emph{\textbf{H5}} is already assumed here.
	\item In order to verify assumptions \textbf{\emph{H6}}, and \textbf{\emph{H7}}, we denote once again for any $\theta \in \Theta$ and any $Z_0$
			\begin{equation}
			m_{\theta}(Z_0) = \log \left( f_{\theta}(Z_0) \right),
			\end{equation}
		where
		\begin{equation}
		f_{\theta}(Z_0) = \int_{\R^k} \prod_{i=1}^k \ind_{I_{Y_{i,0}}-\lambda_{i,0}(\gamma)}(x_i) \varphi_R(x) \d x.
		\end{equation}
	Let $\theta_{\alpha}$ and $\theta_{\beta}$ be two coordinates of the vector $\theta$. If the partial derivatives exist, we actually have
		\begin{equation}
		\dfrac{\partial m_{\theta}}{\partial \theta_{\alpha}}(Z_0) = \dfrac{\partial f_{\theta}}{\partial \theta_{\alpha}}(Z_0) \times \dfrac{1}{f_{\theta}(Z_0)}
		\end{equation}
	and
		\begin{equation}
		\dfrac{\partial^2 m_{\theta}}{\partial \theta_{\beta} \partial_{\alpha}}(Z_0) = \dfrac{\partial^2 f_{\theta}}{\partial \theta_{\beta} \partial_{\alpha}}(Z_0) \times \dfrac{1}{f_{\theta}(Z_0)} - \dfrac{\partial f_{\theta}}{\partial \theta_{\beta}}(Z_0) \times \dfrac{\partial f_{\theta}}{\partial \theta_{\alpha}}(Z_0) \times \dfrac{1}{f_{\theta}(Z_0)^2}.
		\end{equation}
	For future reasons, we will establish the stronger following result, which implies immediately the verification of assumption \textbf{\emph{H6}} and \textbf{\emph{H7}}.
		\begin{lem} \label{lem_strong_bound_deriv}
		Let $q \in ( 1, +\infty)$, the application $\theta \longmapsto m_{\theta}(Z_0)$ is of class $\mathcal{C}^2$ on $\Theta$ and we have 
			\begin{equation}
			\E \left( \underset{\theta \in \Theta}{\sup} \ \left| \dfrac{\partial m_{\theta}}{\partial \theta_{\alpha}}(Z_0) \right|^q \right)  < +\infty \quad \text{and} \quad \E \left( \underset{\theta \in \Theta}{\sup} \ \left| \dfrac{\partial^2 m_{\theta}}{\partial \theta_{\beta} \partial \theta_{\alpha}}(Z_0) \right|^q \right)  < +\infty.
			\end{equation}
		\end{lem}
	The proof of Lemma \ref{lem_strong_bound_deriv} is not exactly straightforward, and requires the following preliminary result.
	
		\begin{lem}\label{lem_bound_deriv_f}
		The application $\theta \longmapsto f_{\theta}(Z_0)$ is of class $\mathcal{C}^2$ on $\Theta$.
		
		Moreover, there exist two polynomials $P,Q$ with positive coefficients, such that for any $\theta \in \Theta$, we have
			\begin{equation}\label{bound_1_deriv_f}
			\left| \dfrac{\partial f_{\theta}}{\partial \theta_{\alpha}}(Z_0) \right| \leqslant \int_{\R^k} \prod_{i=1}^k \ind_{I_{Y_{i,0}}-\lambda_{i,0}(\gamma)}(x_i) \left[ P(\|x\|_2) Q(\|X_{-1}\|_2) \right] \varphi_R(x) \d x
			\end{equation}
		and
			\begin{equation}\label{bound_2_deriv_f}
			\left| \dfrac{\partial^2 f_{\theta}}{\partial \theta_{\beta} \partial \theta_{\alpha}}(Z_0) \right| \leqslant \int_{\R^k} \prod_{i=1}^k \ind_{I_{Y_{i,0}}-\lambda_{i,0}(\gamma)}(x_i) \left[ P(\|x\|_2) Q(\|X_{-1}\|_2) \right] \varphi_R(x) \d x
			\end{equation}
		\end{lem}
		
		\begin{lemproof}[ of Lemma \ref{lem_bound_deriv_f}]
		Due to the number of different cases to deal with, the entire proof of this lemma is quite repetitive and tedious. We only show here one specific case, but the other ones can be treated with the same arguments.
		
		Assume $\theta_{\alpha} = B(i_0,j_0)$, we actually have
			\begin{align}
			f_{\theta}(Z_0) & = \int_{\R^k} \prod_{i=1}^k \ind_{I_{Y_{i,0}}-\lambda_{i,0}(\gamma)}(x_i) \varphi_R(x) \d x \nonumber \\
			& = \int_{\R^{k-1}} \prod_{i\neq i_0} \ind_{I_{Y_{i,0}}-\lambda_{i,0}(\gamma)}(x_i) \dfrac{\exp \left( -\frac{1}{2} \langle \overline{R^{-1}}^{i_0} \overline{x}^{i_0} , \overline{x}^{i_0} \rangle \right)}{\sqrt{(2\pi)^k \det (R)}} \nonumber \\
			& \hspace*{2cm} \times \int_{\R} \ind_{I_{Y_{i_0,0}}-\lambda_{i_0,0}(\gamma)}(x_{i_0}) \exp \left( -\frac{1}{2} R^{-1}(i_0,i_0)x_{i_0}^2 - \sum_{i \neq i_0} R^{-1}(i_0,i) x_{i_0} x_i \right) \d x_{i_0} \d \overline{x}^{i_0}
			\end{align}
		where
			\begin{enumerate}[label=$\bullet$]
			\item $\overline{x}^{i_0}$ denotes the vector $x$ deprived of the $i_0^{\mathrm{th}}$ coordinate,
			\item $\overline{R^{-1}}^{i_0}$ denotes the matrix $R^{-1}$ deprived of the $i_0^{\mathrm{th}}$ line and $i_0^{\mathrm{th}}$ column.
			\end{enumerate}
		Denoting $g_{\theta,Z_0}^{(i_0)}(\overline{x}^{i_0})$ the integrand of this latter integral, we have if $X_{-1}(j_0)=0$
			\begin{equation}
			\dfrac{\partial}{\partial \theta_{\alpha}} \left\lbrace g_{\theta,Z_0}^{(i_0)}(\overline{x}^{i_0}) \right\rbrace = 0
			\end{equation}
		and 
			\begin{equation}
			\dfrac{\partial}{\partial \theta_{\alpha}} \left\lbrace g_{\theta,Z_0}^{(i_0)}(\overline{x}^{i_0}) \right\rbrace = (2Y_{i_0}-1) X_{-1}(j_0) \prod_{i \neq i_0} \ind_{I_{Y_{i,0}}-\lambda_{i,0}(\gamma)}(x_i) \varphi_R \left( x_{\lambda_{i_0,0}(\gamma)} \right)
			\end{equation}
		otherwise, where $ x_{\lambda_{i_0,0}(\gamma)} = \left( x_1,\ldots,x_{i_0-1}, -\lambda_{i_0,0}(\gamma),x_{i_0+1},\ldots , x_k\right)$. Using Lemma \ref{lem_Phi_R_bounds}, there exist positive constants $m_{\det}, c_20$ such that for any $\theta \in \Theta$ and any $\overline{x}^{i_0} \in \R^{k-1}$
			\begin{align}
			\left| \dfrac{\partial}{\partial \theta_{\alpha}} \left\lbrace g_{\theta,Z_0}^{(i_0)}(\overline{x}^{i_0}) \right\rbrace \right| & \leqslant \dfrac{\| X_{-1} \|_2}{\sqrt{(2\pi)^k m_{\det}}} \exp \left(- \frac{c_2}{2} \| x_{\lambda_{i_0,0}(\gamma)} \|_2^2 \right) \nonumber \\
			& \leqslant \dfrac{\| X_{-1} \|_2}{\sqrt{(2\pi)^k m_{\det}}} \exp \left(-\dfrac{c_2}{2} \| \overline{x}^{i_0} \|_2^2 \right)
			\end{align}
		which is integrable on $\R^{k-1}$. Therefore, the partial derivative $\dfrac{\partial f_{\theta}}{\partial \theta_{\alpha}}(Z_0)$ exists, and is either equal to zero or 
			\begin{equation}
			\dfrac{\partial f_{\theta}}{\partial \theta_{\alpha}}(Z_0) = (2 Y_{i_0}-1) X_{-1}(j_0) \int_{\R^{k-1}} \prod_{i \neq i_0} \ind_{I_{Y_{i,0}-\lambda_{i,0}(\gamma)}}(x_i) \varphi_R \left( x_{\lambda_{i_0,0}(\gamma)} \right) \d \overline{x}^{i_0}.
			\end{equation}
		In order to obtain the inequality \eqref{bound_1_deriv_f}, we can use the fundamental theorem of calculus, and we obtain for any $\theta \in \Theta$
			\begin{align}
			\varphi_R \left( x_{\lambda_{i_0,0}(\gamma)} \right) & = \int_{-\infty}^{-\lambda_{i_0,0}(\gamma)} \dfrac{\d }{\d x_{i_0}} \left\lbrace \varphi_R (x) \right\rbrace \d x_{i_0} \nonumber \\
			& = \left| \int_{I_{Y_{i_0,0}}-\lambda_{i_0,0}(\gamma)} \dfrac{\d }{\d x_{i_0}} \left\lbrace \varphi_R (x) \right\rbrace \d x_{i_0}  \right| \nonumber \\
			& \leqslant \int_{I_{Y_{i_0,0}}-\lambda_{i_0,0}(\gamma)} \left| \langle R^{-1}(i_0,\cdot),x \rangle \varphi_R (x) \right| \d x_{i_0} \nonumber \\
			& \leqslant \int_{I_{Y_{i_0,0}}-\lambda_{i_0,0}(\gamma)} M_{\mathcal{R}^{-1}} \|x \|_2 \varphi_R(x) \d x_{i_0}
			\end{align}
		where $M_{\mathcal{R}^{-1}}$ is a positive constant, obtained by compactness of $\mathcal{R}$. Therefore, for any $\theta \in \Theta$
			\begin{equation}
			\left| \dfrac{\partial f_{\theta}}{\partial \theta_{\alpha}}(Z_0) \right| \leqslant M_{\mathcal{R}^{-1}} \int_{\R^{k}} \prod_{i=1}^k \ind_{I_{Y_{i,0}}-\lambda_{i,0}(\gamma)}(x_i) \times \|x \|_2 \times \|X_{-1}\|_2 \times \varphi_R(x) \d x,
			\end{equation}
		which leads to the desired result.
		
		Every other partial derivative, up to order 2, can be bounded the same way, and we obtain Lemma \ref{lem_bound_deriv_f}.
		\end{lemproof}
	
		\begin{lemproof}[ of Lemma \ref{lem_strong_bound_deriv}]
		Take any coordinates $\theta_{\alpha}, \theta_{\beta}$ of the vector $\theta$. By Lemma \ref{lem_bound_deriv_f}, the application $\theta \longmapsto m_{\theta}(Z_0)$ is of class $\mathcal{C}^2$ and we have
			\begin{equation}
		\dfrac{\partial m_{\theta}}{\partial \theta_{\alpha}}(Z_0) = \dfrac{\partial f_{\theta}}{\partial \theta_{\alpha}}(Z_0) \times \dfrac{1}{f_{\theta}(Z_0)}
			\end{equation}
		and
			\begin{equation}
		\dfrac{\partial^2 m_{\theta}}{\partial \theta_{\beta} \partial_{\alpha}}(Z_0) = \dfrac{\partial^2 f_{\theta}}{\partial \theta_{\beta} \partial_{\alpha}}(Z_0) \times \dfrac{1}{f_{\theta}(Z_0)} - \dfrac{\partial f_{\theta}}{\partial \theta_{\beta}}(Z_0) \times \dfrac{\partial f_{\theta}}{\partial \theta_{\alpha}}(Z_0) \times \dfrac{1}{f_{\theta}(Z_0)^2}.
			\end{equation}
		It is however possible to apply Hölder's inequality with respect to the measure of density $\varphi_R$ to inequality \eqref{bound_1_deriv_f}, and we obtain for any $\theta \in \Theta$
			\begin{align} 
			\left| \dfrac{\partial f_{\theta}}{\partial \theta_{\alpha}}(Z_0) \right| & \leqslant C \int_{\R^k} \prod_{i=1}^k \ind_{I_{Y_{i,0}}-\lambda_{i,0}(\gamma)}(x_i) \left[ P(\|x\|_2) Q(\|X_{-1}\|_2) \right] \varphi_R(x) \d x \nonumber  \\
			& \leqslant C \left( \underbrace{\int_{\R^k} \prod_{i=1}^k \ind_{I_{Y_{i,0}}-\lambda_{i,0}(\gamma)}(x_i) \varphi_R (x) \d x}_{f_{\theta}(Z_0)} \right)^{\frac{1}{v}} \times \left( \underbrace{\int_{\R^k} \left[P(\|x\|_2) Q(\|X_{-1}\|_2) \right]^w \varphi_R(x) \d x}_{I_{\theta}(X_{-1})} \right)^{\frac{1}{w}}
			\end{align}
		where $\dfrac{1}{v}+ \dfrac{1}{w}=1$, and $v,w >1$ are to be determined later.
		Using Lemma \ref{lem_Phi_R_bounds}, we have
			\begin{equation}
			I_{\theta}(X_{-1}) \leqslant Q(\|X_{-1}\|_2)^w \times \dfrac{1}{\sqrt{(2\pi)^k m_{\det}}} \int_{\R^k} P(\|x\|_2)^w \exp \left( - \frac{c_2}{2} \|x\|_2^2\right) \d x.
			\end{equation}
		Since the Gaussian distribution admits moments of any order, one can find a positive constant $M_w$ depending only on $w$ such that for any $\theta \in \Theta$
			\begin{equation}
			I_{\theta}(X_{-1})^{\frac{1}{w}} \leqslant M_w Q(\|X_{-1}\|_2).
			\end{equation}
		Thus, we have for any $\theta \in \Theta$
			\begin{equation}
			\left|\dfrac{\partial m_{\theta}}{\partial \theta_{\alpha}}(Z_0) \right| \leqslant C M_w Q(\|X_{-1} \|_2) \times \left(\int_{\R^k} \prod_{i=1}^k \ind_{I_{Y_{i,0}}-\lambda_{i,0}(\gamma)}(x_i) \varphi_R(x) \d x \right)^{\frac{1}{v}-1}.
			\end{equation}
		Yet, by Lemma \ref{lem_Phi_R_bounds}, there exist positive constants $M_{\det}$ and $c_1$ such that for any $\theta \in \Theta$
			\begin{align}
			\int_{\R^k} \prod_{i=1}^k \ind_{I_{Y_{i,0}}-\lambda_{i,0}(\gamma)}(x_i) \varphi_R(x) \d x & \geqslant \dfrac{1}{\sqrt{M_{\det}}} \int_{\R^k} \prod_{i=1}^k \ind_{I_{Y_{i,0}}-\lambda_{i,0}(\gamma)}(x_i) \exp \left(- \frac{c_1}{2} \|x\|_2^2 \right) \d x \nonumber \\ 
			& =  \dfrac{1}{\sqrt{c_1M_{\det}}} \int_{\R^k} \prod_{i=1}^k \ind_{I_{Y_{i,0}}-\sqrt{c_1}\lambda_{i,0}(\gamma)}(x_i) \exp \left(- \frac{\|x\|_2^2}{2} \right) \d x \nonumber \\
			& \geqslant \dfrac{1}{\sqrt{c_1M_{\det}}} \prod_{i=1}^k \min \left\lbrace \Phi\left(-\sqrt{c_1}\lambda_{i,0}(\gamma) \right), 1- \Phi\left(-\sqrt{c_1}\lambda_{i,0}(\gamma) \right)  \right\rbrace \nonumber \\
			& \geqslant \dfrac{d}{\sqrt{c_1M_{\det}}} \exp \left(-c_1 \|\lambda_0(\gamma)\|_2^2 \right) \nonumber
			\shortintertext{where $d$ is the positive constant from Lemma \ref{lem_exp_ineq}} \nonumber 
			& \geqslant \dfrac{d}{\sqrt{c_1M_{\det}}} \exp \left(-c_1 \left( M_{\mathcal{A}} + M_{\mathcal{B}} \|X_{-1}\|_2 \right)^2 \right)
			\end{align}
		for some positive constants $M_{\mathcal{A}}$ and $M_{\mathcal{B}}$, obtained by compactness of $\Theta$. We finally obtain that for any $\theta \in \Theta$
			\begin{equation}
			\left|\dfrac{\partial m_{\theta}}{\partial \theta_{\alpha}}(Z_0) \right| \leqslant C M_w Q(\|X_{-1}\|_2)  \times \left( \dfrac{d}{\sqrt{c_1M_{\det}}} \right)^{\frac{1}{v}-1} \exp \left( c_1 \left(1- \frac{1}{v}\right) \left( M_{\mathcal{A}} + M_{\mathcal{B}} \|X_{-1}\|_2 \right)^2 \right),
			\end{equation}
		and for any $q >1$
			\begin{multline}
			\underset{\theta \in \Theta}{\sup} \ \left|\dfrac{\partial m_{\theta}}{\partial \theta_{\alpha}}(Z_0) \right|^q \leqslant C^q M_w^q Q(\|X_{-1}\|_2)^q \times \left( \dfrac{d}{\sqrt{c_1M_{\det}}} \right)^{q\left(\frac{1}{v}-1\right)} \\ \times \exp \left(qc_1 \left(1- \frac{1}{v}\right) \left( M_{\mathcal{A}} + M_{\mathcal{B}} \|X_{-1}\|_2 \right)^2 \right).
			\end{multline}
		However, $v >1$ is arbitrary, and we can choose it such that $qc_1\left(1- \frac{1}{v}\right)M_{\mathcal{B}}^2 < \kappa$, where $\kappa$ satisfies \eqref{eq_kappa_cons_1} We thus have for any $q >1$
			\begin{equation}
			\E \left(\underset{\theta \in \Theta}{\sup} \ \left|\dfrac{\partial m_{\theta}}{\partial \theta_{\alpha}}(Z_0) \right|^q \right) < +\infty.
			\end{equation}
		Using inequality \eqref{bound_2_deriv_f} from Lemma \ref{lem_bound_deriv_f}, we can show similarly that for any $q >1$
			\begin{equation}
			\E \left(\underset{\theta \in \Theta}{\sup} \ \left| \dfrac{\partial^2 m_{\theta}}{\partial \theta_{\alpha} \partial \theta_{\beta}}(Z_0)\right|^q \right) < +\infty.
			\end{equation}
		\end{lemproof}
	\item Let us now check {\bf H8}. According to Lemma \ref{lem_strong_bound_deriv}, the application
		\begin{equation}
		\theta \longmapsto \E \left( m_{\theta}(Z_0) \right)
		\end{equation}
	is twice differentiable. Thus, we can write the Taylor expansion for all $\theta \in \Theta$
		\begin{equation}
		\E \left( m_{\theta}(Z_0) \right) = \E \left( m_{\theta_0}(Z_0) \right) + \E \left( \nabla m_{\theta_0}(Z_0) \right) \cdot (\theta - \theta_0) + (\theta-\theta_0)' \cdot \E \left( \nabla^2 m_{\theta_0}(Z_0) \right) \cdot (\theta - \theta_0) + o\left( \| \theta-\theta_0 \|_2^2 \right).
		\end{equation}
	Additionally, since $\theta_0$ is the unique maximum of the latter function, we actually have
		\begin{equation}
		\E \left( m_{\theta}(Z_0) \right) = \E \left( m_{\theta_0}(Z_0) \right) + (\theta-\theta_0)' \cdot \E \left( \nabla^2 m_{\theta_0}(Z_0) \right) \cdot (\theta - \theta_0) + o\left( \| \theta-\theta_0 \|_2^2 \right),
		\end{equation}
	and the matrix $\E \left( \nabla^2 m_{\theta_0}(Z_0) \right)$ is thus symmetric semi positive definite. It is then invertible if and only if it is positive definite. As usual for regular conditional likelihood estimators, one can show that non-invertibility of $M:=\E\left(\nabla^2m_{\theta_0}(Z_0)\right)$ is equivalent to the existence of a linear relation between the partial derivatives of $\theta\mapsto f_{\theta}(Z_0)$ at point $\theta_0$. Let us define 
 $$W_t=(\mathds{1}',Y_{t-1}',\ldots,Y_{t-p}',X_{t-1}')',$$
 where $\mathds{1}$ is the vector of $\R^k$ with all the coordinates equal to $1$ and $D$ be the matrix of size $k\times (p+1)k+d$ such that $\lambda_t(\gamma)=D W_t$. Moreover if $D_1,\ldots,D_k$ denote the row vectors of $D$, we have $\lambda_{i,t}(\gamma_0)=D_i'W_t$ for $1\leq i\leq d$. When $M$ is non-invertible, one can find some vectors $\alpha_1,\ldots,\alpha_k$ of $\R^{(p+1)k+d}$ and some real numbers $\kappa_{m,\ell}$ for  $1\leq m<\ell\leq k$ such that a.s.
 \begin{eqnarray*}
0&=&\sum_{i=1}^k (2 Y_{i,1}-1)\alpha_i'W_0\int_{\prod_{j\neq i}\left\{I_{Y_{j,1}}-D_j'W_0\right\}}\phi_R\left(x_1,\ldots,x_{i-1},B_i'W_0,x_{i+1},\ldots,x_k\right)\prod_{j\neq i}dx_j\\
&+& \int_{\prod_{i=1}^k\left\{I_{Y_{i,1}}-Z_i'W_0\right\}}\sum_{1\leq m<\ell\leq k}\kappa_{m,\ell}\frac{\partial \phi_R}{\partial R_{m,\ell}}(x)dx.
\end{eqnarray*}
Since $\P\left(Y_t=y\vert\mathcal{F}_{t-1}\right)$ is positive a.s. for any value of $y\in\{0,1\}^k$, we then get
\begin{eqnarray*}
0&=&\sum_{i=1}^k (2 y_{i,1}-1)\alpha_i'W_0\int_{\prod_{j\neq i}\left\{I_{y_{j,1}}-D_j'W_0\right\}}\phi_R\left(x_1,\ldots,x_{i-1},B_i'W_0,x_{i+1},\ldots,x_k\right)\prod_{j\neq i}dx_j\\\nonumber
&+& \int_{\prod_{i=1}^k\left\{I_{y_{i,1}}-Z_i'W_0\right\}}\sum_{1\leq m<\ell\leq k}\kappa_{m,\ell}\frac{\partial \phi_R}{\partial R_{m,\ell}}(x)dx.
\end{eqnarray*}
almost surely. But summing these relations over $y_{j,1}\in\{0,1\}$ for $j\neq i$, we obtain
$$\alpha_i'W_0=0,\quad 1\leq i\leq k, \mbox{ a.s.}$$
As for the proof for consistency with the identification of $\theta_0$, one can show that $\alpha_i=0$ for $1\leq i\leq k$. Summing now the previous relations over $y_{\ell,1}$ for $\ell\neq i,j$ for some $1\leq i<j\leq k$, we get
$$\kappa_{i,j}\int_{I_{y_{i,1}}-D_i'W_0}\int_{I_{y_{j,1}}-D_j'W_0}\frac{\partial \phi_R}{\partial R_{i,j}}(x_i,x_j)dx_idx_j=0\mbox{ a.s.}$$
However, the latter double integral cannot vanish because it is the derivative of a monotonic function (it is equal to the function $f$ defined in Lemma \ref{lem_increasing_f} if we set $y_{i,1}=y_{j,1}=0$, $r=R_{i,j}$ and $\omega_{\ell}=-D_{\ell}'W_0$ for $\ell=i,j$). Then $\kappa_{i,j}=0$ and Assumption \emph{\textbf{H8}} is satisfied.
	\item Finally, denoting once again
		\begin{equation}
		\psi_{\theta}(y_t , Y_{t-1}^-,X_{t-1}) = \int_{\R^k} \prod_{i=1}^k \ind_{I_{y_{i,t}}}(\lambda_{i,t}(\gamma) + e_i) \varphi_{R}(e) \d e
		\end{equation}
	the density of the conditional distribution $\mathcal{L}_{\theta} \left( Y_t \mid \mathcal{F}_{t-1} \right)$ with respect to the measure $\mu = \sum_{y \in \left\lbrace 0 , 1 \right\rbrace^k} \delta_y$, we have
		\begin{equation}
		\nabla m_{\theta}(Z_t) = \dfrac{\nabla \psi_{\theta}(Y_t,Y_{t-1}^-,X_{t-1})}{\psi_{\theta}(Y_t,Y_{t-1}^-,X_{t-1})}.
		\end{equation}
	Differentiating with respect to $\theta$ both sides of
		\begin{equation}
		\int \psi_{\theta}(y_t,Y_{t-1}^-,X_{t-1}) \mu (\d y_t) = 1,
		\end{equation}
	we obtain for all $\theta \in \Theta$
		\begin{equation}
		\int \nabla \psi_{\theta}(y_t,Y_{t-1},X_{t-1}) \mu (\d y_{t})=0.
		\end{equation}
	Note that this differentiation is licit according to what we established in the previous points, since $\psi_{\theta}(Y_t,Y_{t-1}^-,X_{t-1}) = f_{\theta}(Z_t)$.
		
	Thus, for any positive measurable function $h$ defined on $ \left(\R^d\right)^{\N} \times \left(\R^k\right)^{\N}$, we have
		\begin{align}
		\E \left( h (X_{t-1}^-,\varepsilon_{t-1}^-) \nabla m_{\theta_0}(Z_t) \right) & = \int \int \nabla m_{\theta_0}(y_t,\ldots y_{t-p},x_{t-1}) \psi_{\theta_0}(y_t) \mu (\d y_t) \nonumber \\
		& \hspace*{2cm} h(x_{t-1}^-,e_{t-1}^-) \P_{(X_{t-1}^-,\varepsilon_{t-1}^-)}(\d x_{t-1}^- , \d e_{t-1}^-) \nonumber \\
		& = \int \int \nabla \psi_{\theta_0}(y_t,y_{t-1}^-,x_{t-1}) \mu (\d y_t) \nonumber \\
		& \hspace*{2cm} h(x_{t-1}^-,e_{t-1}^-) \P_{(X_{t-1}^-,\varepsilon_{t-1}^-)}(\d x_{t-1}^- , \d e_{t-1}^-) \nonumber \\
		&  = 0.
		\end{align}
	Thus, $\E \left( \nabla m_{\theta_0}(Z_t) \mid \mathcal{F}_{t-1} \right)=0$, i.e. Assumption \emph{\textbf{H9}} is satisfied.
	\end{enumerate}

\subsubsection{Proof for estimator $\hat{\gamma}$}	

The proof for estimator $\hat{\gamma}$ is quite similar to the one for $\hat{\theta}_{\PL}$. Once again we have to check assumptions \emph{\textbf{H1}} to \emph{\textbf{H9}} from Theorem \ref{theo_pfanzagl} and Theorem \ref{theo_straumann}, for the same process $(Z_t)_{t\in \Z}$ and the family of applications $(m_{\gamma})_{\gamma \in \Gamma}$ defined by
	\begin{equation}
	m_{\gamma}(Z_t) =\sum_{i=1}^k Y_{i,t}\log \left[ \Phi (\lambda_{i,t}(\gamma)) \right] + (1-Y_{i,t}) \log \left[ \Phi (-\lambda_{i,t}(\gamma) \right].
	\end{equation}
	\begin{enumerate}[label=$\triangleright$,labelindent=0pt,labelwidth=!,wide]
	\item \emph{\textbf{H1}} and \emph{\textbf{H2}} are satisfied for the same reasons as before.
	\item The application $\gamma \longmapsto m_{\gamma}(Z_0)$ is continuous almost surely by a simple composition of continuous functions.
		
	The verification of the second point of \emph{\textbf{H3}} relies once again on a compactness argument and the following inequality (see Lemma \ref{lem_exp_ineq})
		\begin{equation} \label{ineq_exp_d}
		\forall x \in \R, \ 1-\Phi(x) \geqslant d \e^{-x^2}
		\end{equation}
	for some constant $d >0$.

	Indeed, there exist positive constants $M_{\mathcal{A}}$ and $M_{\mathcal{B}}$ such that for all $\gamma \in \Gamma$ and all $i \in \left\lbrace 1, \ldots , k\right\rbrace$
		\begin{equation}
		\left| \lambda_{i,0}(\gamma) \right| \leqslant M_{\mathcal{A}} + M_{\mathcal{B}}\| X_{-1}\|_2.
		\end{equation}
	We thus have for all $\gamma \in \Gamma$
		\begin{align}
		\left| m_{\gamma}(Z_0) \right| & = -\left\lbrace\sum_{i=1}^k Y_{i,0} \log \left[ \Phi (\lambda_{i,0}(\gamma)) \right] + (1-Y_{i,0}) \log \left[ \Phi (-\lambda_{i,0}(\gamma)) \right] \right\rbrace \nonumber \\
		 & \leqslant -2k \log \left[ \Phi \left(-M_{\mathcal{A}}-M_{\mathcal{B}}\|X_{-1}\|_2\right) \right] \nonumber \\
		 & =-2k \log \left[ 1- \Phi \left(M_{\mathcal{A}}+M_{\mathcal{B}}\|X_{-1}\|_2\right) \right].
		\end{align}	
	Then, from inequality \eqref{ineq_exp_d}, we obtain
	\begin{align}
		\left| m_{\gamma} (Z_0) \right| &  \leqslant -2k \log \left( d\exp \left( - (M_{\mathcal{A}}+M_{\mathcal{B}} \| X_{-1} \|_2)^2 \right) \right) \nonumber \\
		 & \leqslant -2k \log (d) + 2k \left( M_{\mathcal{A}}+M_{\mathcal{B}} \| X_{-1} \|_2 \right)^2.
		\end{align}
	Hence
		\begin{equation}
		\E \left( \underset{\gamma \in \Gamma}{\sup} \ \left| m_{\gamma}(Z_0)\right| \right) < +\infty
		\end{equation}
	and assumption \emph{\textbf{H3}} is satisfied.
	\item Using the same notations as before, we have
		\begin{equation}
		\E \left( m_{\gamma}(Z_0) \right) - \E \left( m_{\gamma_0}(Z_0) \right) = \sum_{i=1}^k \E \left( \log \left( \dfrac{\psi_{\gamma}(Y_{i,0},Y_{-1}^-,X_{-1})}{\psi_{\gamma_0}(Y_{i,0},Y_{-1}^-,X_{-1})} \right) \right),
		\end{equation}
	where by assumption \textbf{B2}, $\psi_{\gamma} ( y_{i,0},Y_{-1},X_{-1}) = \Phi (\lambda_{i,0}(\gamma))^{y_{i,0}} \Phi (-\lambda_{i,0}(\gamma))^{1-y_{i,0}}$ is the density of the conditional distribution $\mathcal{L}_{\gamma} (Y_{i,0} \mid Y_{-1}^- , X_{-1})$ with respect to the measure $\mu = \delta_0 + \delta_1$ on $\left\lbrace 0, 1 \right\rbrace$.
		
	The same reasoning as before ensures that we have $ \E \left( m_{\gamma}(Z_0) \right) \leqslant \E \left( m_{\gamma_0}(Z_0) \right)$ with equality if and only if $\gamma=\gamma_0$, so assumption \emph{\textbf{H4}} is verified.
	\item Assumption \emph{\textbf{H5}} is already assumed here.
	\item A simple computation gives for all $(i_1,i_2) \in \left\lbrace 1, \ldots , k \right\rbrace^2$ and all $l \in \left\lbrace 1 , \ldots ,p\right\rbrace$
		\begin{equation}
		\dfrac{\partial m_{\gamma}(Z_0)}{\partial A_l(i_1,i_2)} = Y_{-l,i_2}\left( Y_{i_1,0} \dfrac{\varphi (\lambda_{i_1,0}(\gamma))}{\Phi (\lambda_{i_1,0}(\gamma))} + (Y_{i_1,0}-1) \dfrac{\varphi (-\lambda_{i_1,0}(\gamma))}{\Phi (-\lambda_{i_1,0}(\gamma))} \right),
		\end{equation}
	which leads to
		\begin{equation}
		\left| \dfrac{\partial m_{\gamma}(Z_0)}{\partial A_l(i_1,i_2)} \right|^2 =  Y_{-l,i_2} \times \dfrac{\varphi (\lambda_{i_1,0}(\gamma))^2}{\Phi(\lambda_{i_1,0}(\gamma))^2} \ind_{\lambda_{i_1,0}(\gamma) >0} + Y_{-l,i_2} \times \dfrac{\varphi (-\lambda_{i_1,0}(\gamma))^2}{\Phi(-\lambda_{i_1,0}(\gamma))^2} \ind_{\lambda_{i_1,0}(\gamma) \leqslant 0}.
		\end{equation}
	Yet, for any $x \geqslant 0$, we have $\Phi (x) \geqslant \dfrac{1}{2}$, thus $\dfrac{\varphi(x)^2}{\Phi(x)^2} \leqslant 4$, and it follows
		\begin{equation}
		\underset{\gamma \in \Gamma}{\sup} \ \left| \dfrac{\partial m_{\gamma}(Z_0)}{\partial A_l(i_1,i_2)} \right|^2 \leqslant 4.
		\end{equation}
	Similarly, we have for any $i_1 \in \left\lbrace 1, \ldots , k\right\rbrace$ and $i_2 \in \left\lbrace 1, \ldots , d \right\rbrace$
		\begin{equation}\label{eq_H6_theo2}
		\dfrac{\partial m_{\gamma}(Z_0)}{\partial B(i_1,i_2)} = X_{i_2,-1} \left( Y_{i_1,0} \dfrac{\varphi (\lambda_{i_1,0}(\gamma))}{\Phi (\lambda_{i_1,0}(\gamma))} + (Y_{i_1,0}-1) \dfrac{\varphi (-\lambda_{i_1,0}(\gamma))}{\Phi (-\lambda_{i_1,0}(\gamma))} \right),
		\end{equation}
	and we obtain
		\begin{equation}
		\underset{\gamma \in \Gamma}{\sup} \ \left| \dfrac{\partial m_{\gamma}(Z_0)}{\partial B(i_1,i_2)} \right|^2 \leqslant 4 \|X_{-1}\|_2^2.
		\end{equation}
	Since $\E \left( \|X_{-1}\|_2^2 \right) <+\infty$, we have
		\begin{equation}
		\E \left( \underset{\gamma \in \Gamma}{\sup} \ \| \nabla m_{\gamma}(Z_0)\|_2^2 \right) <+\infty,
		\end{equation}
	and \emph{\textbf{H6}} is satisfied.
	\item Setting for any $x \in \R$, $g(x)=\dfrac{\varphi(x)}{\Phi(x)}$, and using equality \eqref{eq_H6_theo2}, we have for any $i_1,i_2,i_3,i_4$ in $\left\lbrace 1, \ldots , k \right\rbrace$ and any $l, l'$ in $\left\lbrace 1, \ldots , p \right\rbrace$
		\begin{align}
		\dfrac{\partial^2 m_{\gamma}(Z_0)}{\partial A_{l'}(i_3,i_4) \partial A_l(i_1,i_2)}  & = Y_{-l,i_2} \left( Y_{i_1,0} \dfrac{\partial g (\lambda_{i_1,0}(\gamma))}{\partial A_{l'}(i_3,i_4)} +(Y_{i_1,0}-1) \dfrac{\partial g(-\lambda_{i_1,0}(\gamma))}{\partial A_{l'}(i_3,i_4)} \right) \nonumber \\
		 & = \left\lbrace 
		 \begin{array}{l}
		 Y_{-l,i_2} Y_{-l',i_4} \left( Y_{i_1,0} \cdot g'(\lambda_{i_1,0}(\gamma)) - (Y_{i_1,0}-1)\cdot g'(\lambda_{i_1,0}(\gamma)) \right) \quad \text{if} \ i_3=i_1 \\
		 0 \quad \text{otherwise}.
		 \end{array}
		 \right.
		\end{align}
	Yet for any $x \geqslant 0$
		\begin{equation}
		 \left| g'(x) \right| = \dfrac{\e^{-\frac{x^2}{2}}}{2\pi} \times \dfrac{x\sqrt{2\pi} \Phi(x) + \e^{-\frac{x^2}{2}}}{\Phi(x)^2} \leqslant \dfrac{2}{\pi} \times \left(1+x\sqrt{2\pi} \right),
		\end{equation}
	thus we obtain
		\begin{equation}
		\left| \dfrac{\partial^2 m_{\gamma}(Z_0)}{\partial A_{l'}(i_3,i_4) \partial A_l(i_1,i_2)} \right| \leqslant \dfrac{2}{\pi}  \left(1+\left| \lambda_{i_1,0}(\gamma) \right|\sqrt{2\pi} \right).
		\end{equation}
	By a compactness argument, there exist two positive constants $M_{\mathcal{A}}$ and $M_{\mathcal{B}}$ such that for any $\gamma \in \Gamma$
		\begin{equation}
		\left| \dfrac{\partial^2 m_{\gamma}(Z_0)}{\partial A_{l'}(i_3,i_4) \partial A_l(i_1,i_2)} \right| \leqslant \dfrac{2}{\pi} \left( 1+ \sqrt{2\pi} \left( M_{\mathcal{A}}+ M_{\mathcal{B}} \|X_{-1} \|_2 \right)\right),
		\end{equation}
	and since $\E \left( \| X_{-1} \|_2 \right) < +\infty$, we obtain
		\begin{equation}
		\E \left(  \underset{\gamma \in \Gamma}{\sup} \ \left| \dfrac{\partial^2 m_{\gamma}(Z_0)}{\partial A_{l'}(i_3,i_4) \partial A_l(i_1,i_2)} \right| \right) < + \infty.
		\end{equation}
	Similarly, we have for any $i_1,i_2,i_3$ in $\left\lbrace 1 , \ldots , k\right\rbrace$, any $i_4 \in \left\lbrace 1 , \ldots ,d \right\rbrace$, and any $l \in \left\lbrace 1, \ldots , p \right\rbrace$
		\begin{equation}
		\dfrac{\partial^2 m_{\gamma}(Z_0)}{\partial B(i_3,i_4)\partial A(i_1,i_2)} = \left\lbrace
		\begin{array}{l}
		Y_{-l,i_2}X_{i_4,-1} \left( Y_{i_1,0} \cdot g'(\lambda_{i_1,0}(\gamma)) - (Y_{i_1,0}-1)\cdot g'(-\lambda_{i_1,0}(\gamma)) \right) \quad \text{if} \ i_1 = i_3 \\
		0 \quad \text{otherwise.}
		\end{array}
		\right.
		\end{equation}
	With the same compactness argument as before we obtain that for any $\gamma \in \Gamma$
		\begin{align}
		\left| \dfrac{\partial^2 m_{\gamma}(Z_0)}{\partial B(i_3,i_4)\partial A(i_1,i_2)} \right| & \leqslant \|X_{-1}\|_2 \cdot \left| g'(\lambda_{i_1,0}(\gamma))\right| \nonumber \\
		 & \leqslant\dfrac{2\|X_{-1} \|_2}{\pi} \left( 1+ \left| \lambda_{i_1,0}(\gamma) \right|\cdot \sqrt{2\pi} \right) \nonumber \\
		 & \leqslant \dfrac{2\|X_{-1}\|_2}{\pi} \left( 1+\sqrt{2\pi}(M_{\mathcal{A}} + M_{\mathcal{B}}\|X_{-1}\|_2 \right)
		\end{align}
	for some positive constants $M_{\mathcal{A}}$ and $M_{\mathcal{B}}$.
		
	Since $\E \left( \| X_{-1}\|_2^2 \right) < + \infty$, we then obtain
		\begin{equation}
		\E \left( \underset{\gamma \in \Gamma}{\sup} \ \left| \dfrac{\partial^2 m_{\gamma}(Z_0)}{\partial B(i_3,i_4)\partial A(i_1,i_2)} \right| \right) < + \infty.
		\end{equation}
	Finally, we have for any $i_1,i_3$ in $\left\lbrace 1, \ldots , k\right\rbrace$ and any $i_2,i_4$ in $\left\lbrace 1 , \ldots , d \right\rbrace$
		\begin{equation}
		\dfrac{\partial^2 m_{\gamma}(Z_0)}{\partial B(i_3,i_4)\partial B(i_1,i_2)} =\left\lbrace
		\begin{array}{l}
		X_{i_4,-1}X_{i_2,-1} \left( Y_{i_1,0}\cdot g'(\lambda_{i_1,0}(\gamma)) - (Y_{i_1,0}-1) \cdot g'(-\lambda_{i_1,0}(\gamma)) \right) \quad \text{if} \ i_3=i_1 \\
		0 \quad \text{otherwise}.
		\end{array}
		\right.
		\end{equation}
	Yet, for any $x \geqslant 0$, we have $\Phi(x) \geqslant \dfrac{\e^{-\frac{x^2}{2}}}{2}$, thus
		\begin{equation}
		\forall x \geqslant 0, \left| g'(x) \right| = \left| -\dfrac{\e^{-\frac{x^2}{2}}}{2\pi} \times \dfrac{x\sqrt{2\pi}\Phi(x)+\e^{-\frac{x^2}{2}}}{\Phi(x)^2} \right| \leqslant \dfrac{2}{\pi} + \sqrt{\dfrac{2}{\pi}}x\e^{-\frac{x^2}{2}}.
		\end{equation}
	Furthermore, since $x\e^{-\frac{x^2}{2}} \underset{x \to +\infty}{\longrightarrow} 0$, there exists $x_0 >0$ such that for all $x \geqslant x_0$, $x\e^{-\frac{x^2}{2}} \leqslant 1$. Hence, for all $\gamma \in \Gamma$ we have
		\begin{equation}
		\left| \dfrac{\partial^2 m_{\gamma}(Z_0)}{\partial B(i_3,i_4)\partial B(i_1,i_2)} \right| \leqslant \| X_{-1}\|_2^2 \cdot \left(\dfrac{2}{\pi} + \sqrt{\dfrac{2}{\pi}}(x_0+1) \right).
		\end{equation}
	Since $\E \left( \| X_{-1} \|_2^2 \right) <+ \infty$, we then have 
		\begin{equation}
		\E \left( \underset{\gamma \in \Gamma}{\sup} \ \left| \dfrac{\partial^2 m_{\gamma}(Z_0)}{\partial B(i_3,i_4)\partial B(i_1,i_2)} \right| \right) < + \infty,
		\end{equation}
	and it yields
		\begin{equation}
		\E \left( \underset{\gamma \in \Gamma}{\sup} \ \| \nabla^2 m_{\gamma}(Z_0) \|_2 \right) < +\infty,
		\end{equation}
	so assumption \emph{\textbf{H7}} is satisfied.
	\item Assumption \emph{\textbf{H8}} can be verified the exact same way as in the proof for estimator $\hat{\theta}_{\PL}$.
	\item The verification of \emph{\textbf{H9}} follows the same principle as in subsection 6.3.1, yet we have to verify that the application $\int \psi_{\gamma}(y_{i,t},Y_{t-1},X_{t-1}) \mu (\d y_{i,t})$ is differentiable with respect to $\gamma$, where
		\begin{equation}
		\psi_{\gamma}(y_{i,t},Y_{t-1}^-,X_{t-1}) = \Phi(\lambda_{i,t}(\gamma))^{y_{i,t}}\Phi(-\lambda_{i,t}(\gamma))^{1-y_{i,t}}
		\end{equation}
	is the density of the conditional distribution $\mathcal{L}_{\gamma}\left( Y_{i,t} \mid Y_{t-1}^-,X_{t-1} \right)$ with respect to the measure $\mu = \delta_0 + \delta_1$ on $\left\lbrace 0, 1 \right\rbrace$.
		
	To do so, note that for any $l \in \left\lbrace 1, \ldots , p\right\rbrace$, and any $i_1,i_2$ in $\left\lbrace 1,\ldots, k\right\rbrace$, we have
		\begin{equation}
		\dfrac{\partial \psi_{\gamma}(y_{i,t},Y_{t-1}^-,X_{t-1})}{\partial A_l(i_1,i_2)} = \left( Y_{i_2,t-l}\varphi (\lambda_{i,t}(\gamma)) \right) \ind_{\lambda_{i,t}(\gamma)>0} - \left( Y_{i_2,t-l}\varphi (-\lambda_{i,t}(\gamma)) \right) \ind_{\lambda_{i,t}(\gamma)\leqslant 0}
		\end{equation}
	thus for any $\gamma \in \Gamma$
		\begin{equation}
		\left| \dfrac{\partial \psi_{\gamma}(y_{i,t},Y_{t-1}^-,X_{t-1})}{\partial A_l(i_1,i_2)} \right| \leqslant 2.
		\end{equation}
	Similarly, we have for any $i_1 \in \left\lbrace 1, \ldots , k \right\rbrace$ and $i_2 \in \left\lbrace 1, \ldots, d \right\rbrace$
		\begin{equation}
		\dfrac{\partial \psi_{\gamma}(y_{i,t},Y_{t-1}^-,X_{t-1})}{\partial B(i_1,i_2)} =  \left( X_{i_2,t-1}\varphi (\lambda_{i,t}(\gamma)) \right) \ind_{\lambda_{i,t}(\gamma)>0} - \left( X_{i_2,t-1}\varphi (-\lambda_{i,t}(\gamma)) \right) \ind_{\lambda_{i,t}(\gamma)\leqslant 0},
		\end{equation}
	thus for any $\gamma \in \Gamma$
		\begin{equation}
		\left| \dfrac{\partial \psi_{\gamma}(y_{i,t},Y_{t-1}^-,X_{t-1})}{\partial B(i_1,i_2)} \right| \leqslant 2 \|X_{t-1}\|_2.
		\end{equation}
	The differentiation of $\gamma \longmapsto \int \psi_{\gamma}(y_{i,t},Y_{t-1}^-,X_{t-1}) \mu (\d y_{i,t})$ is thus licit.
	\end{enumerate}

\subsubsection{Proof for estimator $\hat{r}$}

For any couple $(i,j) \in \left\lbrace 1, \ldots ,k\right\rbrace$ with $i<j$, $r_{i,j}$ belongs to a compact set $S_{i,j}$ of $[-1,1]$, and setting for $s \in S_{i,j}$
	\begin{equation}
	\ell_{i,j}(\gamma,s,Z_t) = \log \left( \int_{I_{Y_{i,t}}-\lambda_{i,t}(\gamma)} \Phi \left[ (2 Y_{i,t}-1) \dfrac{\lambda_{j,t}(\gamma)+sx_i}{\sqrt{1-s^2}} \right] \varphi (x_i) \d x_i \right),
	\end{equation}
we actually have
	\begin{equation}
	\hat{r}_{i,j} = \underset{s \in S_{i,j}}{\argmax} \ \sum_{t=p+1}^T \ell_{i,j}(\hat{\gamma},s,Z_t).
	\end{equation}
We first deal with the consistency of the estimator $\hat{r}_{i,j}, \ i<j$. To this end, we check that assumptions \emph{\textbf{I1}} to \emph{\textbf{I3}} of Proposition \ref{prop8_debaly}
	\begin{enumerate}[label=$\triangleright$,labelindent=0pt,labelwidth=!,wide]
	\item Assumption \emph{\textbf{I1}} is satisfied since $(Z_t)_{t\in \Z}$ is ergodic.
	\item Using the continuity of $\gamma\mapsto\lambda_{i,t}(\gamma)$, $\Phi$, the boundedness of $\Phi$ and Lebesgue's dominated convergence theorem, one can deduce that the mapping
		\begin{equation}
		(\gamma,s)\mapsto \int_{I_{Y_{i,t}}-\lambda_{i,t}(\gamma)}\Phi\left[(2Y_{i,t}-1)\dfrac{\lambda_{j,t}(\gamma)-s x_i}{\sqrt{1-s^2}}\right]\varphi(x_i) \d x_i
		\end{equation}
	is almost surely continuous.
	\item We then show that 
		\begin{equation} \label{eq_cons_r_single}
		\E \left( \underset{\gamma \in \Gamma, s \in S_{i,j}}{\sup} \ \left| \ell_{i,j}(\gamma,s,Z_0) \right| \right) < +\infty.
		\end{equation}
	By Lemma \ref{lem_exp_ineq}, there exists a positive number $d$ such that
		\begin{equation}
		\forall x \in \R, \ \Phi(x)\geqslant d\exp\left(-x^2\right).
		\end{equation}
	Thus, for any $\gamma \in \Gamma$ and any $s \in S_{i,j}$
		\begin{equation}
		\Phi\left[(2Y_{i,0}-1)\dfrac{\lambda_{j,0}(\gamma)+s x_i}{\sqrt{1-s^2}}\right]\geqslant d \exp\left(-\dfrac{(\lambda_{j,0}(\gamma)+sx_i)^2}{1-s^2}\right)\geqslant d\exp\left(-\frac{2\lambda_{j,0}(\gamma)^2+2s^2x_i^2}{1-s^2}\right).
		\end{equation}
	Setting $h(s)=\dfrac{2 s^2}{1-s^2}+\dfrac{1}{2}$, we then get
		\begin{eqnarray}
		\ell_{i,j}(\gamma,s,Z_0) & \geqslant & \log(d)-\dfrac{2\lambda_{j,0}(\gamma)^2}{1-s^2}+\log \left( \int_{I_{Y_{i,0}}-\lambda_{i,0}(\gamma)}\exp\left(-\frac{2s^2 x_i^2}{1-s^2}\right)\varphi(x_i)\d x_i \right) \nonumber \\
		& = & \log (d) - \dfrac{2\lambda_{j,0}(\gamma)^2}{1-s^2} - \dfrac{1}{2}\log (2h(s)) + \log \left( \int_{I_{Y_{i,0}}-\sqrt{2h(s)}\lambda_{i,0}(\gamma)} \varphi (x_i) \d x_i \right) \nonumber \\
		& \geqslant & \log (d) - \dfrac{2\lambda_{j,0}(\gamma)^2}{1-s^2} - \dfrac{1}{2}\log (2h(s)) \nonumber \\
		& & \hspace*{2cm} + \log \left( \min \left\lbrace \Phi (-\sqrt{2h(s)}\lambda_{i,0}(\gamma)),1-\Phi(-\sqrt{2h(s)}\lambda_{i,0}(\gamma)) \right\rbrace \right) \nonumber \\
               & \geqslant & 2\log(d)-\dfrac{2\lambda_{j,0}(\gamma)^2}{1-s^2}-\dfrac{1}{2}\log(2h(s))-2h(s)\lambda_{i,0}(\gamma)^2.
		\end{eqnarray}
	By a compactness argument, there exist $M_1,M_2>0$ such that
		\begin{eqnarray}
			\underset{\gamma \in \Gamma, s \in S_{i,j}}{\sup} \ \left| \ell_{i,j}(\gamma,s,Z_0) \right| & = & \underset{\gamma \in \Gamma, s \in S_{i,j}}{\sup} \left\lbrace - \ell_{i,j}(\gamma,s,Z_0) \right\rbrace \nonumber \\
			 & \leqslant & M_1 + M_2 \| X_{-1}\|_2^2
		\end{eqnarray}
	Since $\E \left( \| X_{-1} \|_2^2 \right) < +\infty$ using our assumptions, we get \eqref{eq_cons_r_single}, and assumption \emph{\textbf{I2}} is satisfied.
	\item One can finally show the exact same way as it was done in subsections 6.3.1 and 6.3.2, that for any $ r_{i,j} \in S_{i,j}$
		\begin{equation}
		\E \left(\ell_{i,j}(\gamma_0,r_{i,j},Z_0) \right) \leqslant \E \left( \ell_{i,j}(\gamma,r_{0,i,j},Z_0 \right).
		\end{equation}
	Furthermore, the equality case leads to the equality 
		\begin{equation}
		\int_{-\infty}^{-\lambda_{i,0}(\gamma_0)} \Phi \left( \dfrac{-\lambda_{j,0}(\gamma_0)-r_{i,j}x}{\sqrt{1-r_{i,j}^2}} \right) \varphi (x) \d x = \int_{-\infty}^{-\lambda_{i,0}(\gamma_0)} \Phi \left( \dfrac{-\lambda_{j,0}(\gamma_0)-r_{0,i,j}x}{\sqrt{1-r_{0,i,j}^2}} \right) \varphi (x) \d x 
		\end{equation}
	which is the same as \eqref{eq_ident_r_single}. We thus have $r_{i,j}=r_{0,i,j}$ in this case, and \emph{\textbf{I3}} is satisfied.
	\end{enumerate}
	From Proposition \ref{prop8_debaly}, we then deduce that
		\begin{equation}
		\hat{r} \underset{T \to + \infty}{\longrightarrow} r_0 \quad \mathrm{a.s.}
		\end{equation}
The stochastic expansion of $\sqrt{T-p} \left( \hat{r}-r_0\right)$ can be established using Proposition \ref{prop_sto_exp_refine}. Note that
	\begin{equation}
	\hat{r} = \underset{r \in \mathcal{R}}{\argmax} \ \sum_{p+1 \leqslant t \leqslant T} \ell_{t}(\hat{\gamma},r).
	\end{equation}
where $\ell_t$ is defined in \eqref{eq_def_ell}.

In what follows, we check that conditions \emph{\textbf{I4}} to \emph{\textbf{I8}} of this latter proposition are satisfied.
	\begin{enumerate}[label=$\triangleright$,labelindent=0pt,labelwidth=!,wide]
	\item Condition \emph{\textbf{I4}} is satisfied by assumption.
	\item We then check assumption \textbf{\emph{I6}} and prove that
		\begin{equation} \label{eq_two_moments}
		\E \left( \underset{\theta \in \Theta}{\sup} \ \| \nabla_{1,2} \ell_{0}(\theta) \| \right) < +\infty \quad \text{and} \quad \E \left( \underset{\theta \in \Theta}{\sup} \ \| \nabla_2^2 \ell_{0}(\theta) \| \right) < +\infty.
		\end{equation}
	Using the fact that $\ell_{0}(\theta)$ is a sum of functions involving different parameters, we fix $i$ and $j$ with $1\leqslant i<j\leqslant k$ and set for any $\gamma \in \Gamma$ and $s \in S_{i,j}$ 
		\begin{equation}
		\widetilde{\ell}_{i,j}(\gamma,s,Z_0)=\int_{I_{Y_{i,0}}-\lambda_{i,0}(\gamma)}\Phi\left(g_{j,x}(\gamma,s,Z_0)\right)\varphi(x) \d x,
		\end{equation}
	where $g_{j,x}(\gamma,s,Z_0) = \left( 2 Y_{j,0}-1 \right) \dfrac{\lambda_{j,0}(\gamma)+sx}{\sqrt{1-s^2}}$. This way, we actually have
		\begin{equation}
		\ell_{0}(\theta) = \sum_{1 \leqslant i < j \leqslant k} \log \left( \widetilde{\ell}_{i,j}(\gamma, r_{i,j},Z_0) \right).
		\end{equation}
	Next we denote by $\gamma_{\alpha}$ one of the coordinates of the vector $\gamma$. Note that $\gamma_{\alpha}$ appears in the expression of either $\lambda_{i,0}(\gamma)$ or $\lambda_{j,0}(\gamma)$ but not
in both. First, computing the different partial derivatives, one can find a positive constant $C$ such that
		\begin{equation}
		\left| \dfrac{\partial g_{j,x}}{\partial s}(\gamma,s,Z_0)\right|\leqslant C\left(1 + |x| + \| X_{-1}\|_2\right),\quad \left| \dfrac{\partial g_{j,x}}{\partial \gamma_{\alpha}}(\gamma,s,Z_0) \right| \leqslant C \left( 1 + \| X_{-1} \|_2 \right),
		\end{equation}
	and
		\begin{equation}
		\left| \dfrac{\partial^{2}g_{j,x}}{\partial s \partial \gamma_{\alpha}}(\gamma,s,Z_0)\right| \leqslant C\left(1+\| X_{-1}\|_2 \right), \quad \left| \dfrac{\partial^{2}g_{j,x}}{\partial s^2}(\gamma,s,Z_0)\right| \leqslant C\left( 1 + |x| + \|X_{-1}\|_2 \right).
		\end{equation}
	We next compute the first order partial derivatives of $\widetilde{\ell}_{i,j}(\gamma,s,Z_0)$. More precisely, one can apply the theorem for the differentiability of a function defined by an integral to get 
		\begin{equation}
		\dfrac{\partial \widetilde{\ell}_{i,j}}{\partial s}(\gamma,s,Z_0)=\int_{I_{Y_{i,0}}-\lambda_{i,0}(\gamma)}\varphi\left(g_{j,x}(\gamma,s,Z_0)\right)\varphi(x)\dfrac{\partial g_{j,x}}{\partial s}(\gamma,s,Z_0) \d x.
		\end{equation}
	Moreover, if $\gamma_{\alpha}$ is a coordinate that appears in the expression of $\lambda_{i,0}(\gamma)$
		\begin{equation}
		\dfrac{\partial \tilde{\ell}_{i,j}}{\partial \gamma_{\alpha}}(\gamma,s,Z_0) = (2 Y_{i,0}-1) \dfrac{\partial \lambda_{i,0}}{\partial \gamma_{\alpha}}(\gamma) \times \Phi (g_{j,x}(\gamma,s,Z_0))\varphi (-\lambda_{i,0}(\gamma)).
		\end{equation}
	If $\gamma_{\alpha}$ is a coordinate that appears in the expression of $\lambda_{j,0}(\gamma)$, we have
		\begin{equation}
		\dfrac{\partial \tilde{\ell}_{i,j}}{\partial \gamma_{\alpha}}(\gamma,s,Z_0) = \int_{I_{Y_{i,0}}-\lambda_{i,0}(\gamma)} \varphi (g_{j,x}(\gamma,s,Z_0)) \varphi (x) \times \dfrac{\partial g_{j,x}}{\partial \gamma_{\alpha}}(\gamma,s,Z_0) \ \d x.
		\end{equation}
	For the second order partial derivatives of $\widetilde{\ell}_{i,j}(\gamma,s,Z_0)$, we have
		\begin{eqnarray}\label{eq_2nd_partial_derivative}
		\dfrac{\partial^{2} \widetilde{\ell}_{i,j}}{\partial s^2}(\gamma,s, Z_0) & = & -\int_{I_{Y_{i,0}}-\lambda_{i,0}(\gamma)}g_{j,x}(\gamma,s,Z_0)\varphi\left(g_{j,x}(\gamma,s,Z_0)\right)\varphi(x)\left| \dfrac{\partial g_{j,x}}{\partial s}(\gamma,s,Z_0)\right|^2 \d x \nonumber \\
         & + & \int_{I_{Y_{i,0}}-\lambda_{i,0}(\gamma)}\varphi\left(g_{j,x}(\gamma,s,Z_0)\right)\varphi(x)\dfrac{\partial^{2} g_{j,x}}{\partial s^2}(\gamma,s,Z_0) \d x.
		\end{eqnarray}
	Moreover, if $\gamma_{\alpha}$ is a coordinate of $\gamma$ that appears in the expression of $\lambda_{i,0}(\gamma)$, 
		\begin{equation}
		\dfrac{\partial^{2} \widetilde{\ell}_{i,j}}{\partial s \partial \gamma_{\alpha}}(\gamma,s,Z_0)= \left( 2 Y_{i,0}-1 \right) \varphi\left(g_{j,-\lambda_{i,0}(\gamma)}(\gamma,s,Z_0)\right)\varphi\left(-\lambda_{i,0}(\gamma)\right) \dfrac{\partial g_{j,-\lambda_{i,0}(\gamma)}}{\partial s}(\gamma,s,Z_0) \dfrac{\partial \lambda_{i,0}}{\partial \gamma_{\alpha}}(\gamma).
		\end{equation}
	Note that we can always modify this derivative using the following integral expression that can be obtained from the fundamental theorem of calculus
		\begin{align}
		\left( 2Y_{i,0}-1 \right) \varphi \left( g_{j,-\lambda_{i,0}(\gamma)}(\gamma,s,Z_0) \right) \varphi (-\lambda_{i,0}(\gamma)) & = \int_{I_{Y_{i,0}}-\lambda_{i,0}(\gamma)} \varphi \left( g_{j,x}(\gamma,s,Z_0) \right) \varphi (x) \nonumber \\
		 & \hspace*{2cm} \times \left( \dfrac{(2Y_{j,0}-1)s}{\sqrt{1-s^2}} g_{j,x}(\gamma,s,Z_0) +x \right) \d x.
		\end{align}
	Thus, if $\gamma_{\alpha}$ is a coordinate of $\gamma$ that appears in the expression of $\lambda_{i,0}(\gamma)$
		\begin{align}
		\dfrac{\partial^{2} \widetilde{\ell}_{i,j}}{\partial s \partial \gamma_{\alpha}}(\gamma,s,Z_0) & = \int_{I_{Y_{i,0}}-\lambda_{i,0}(\gamma)} \varphi \left( g_{j,x}(\gamma,s,Z_0) \right) \varphi (x)\left( \dfrac{(2Y_{j,0}-1)s}{\sqrt{1-s^2}} g_{j,x}(\gamma,s,Z_0) +x \right) \nonumber \\ 
		& \hspace*{4cm}\times \dfrac{\partial g_{j,-\lambda_{i,0}(\gamma)}}{\partial s}(\gamma,s,Z_0) \dfrac{\partial \lambda_{i,0}}{\partial \gamma_{\alpha}}(\gamma) \d x
		\end{align}
	On the other hand, when $\gamma_{\alpha}$ is a coordinate of $\gamma$ that appears in the expression of $\lambda_{j,0}(\gamma)$, we have the expression
		\begin{eqnarray}
		\dfrac{\partial^{2} \widetilde{\ell}_{i,j}}{\partial s \partial \gamma_{\alpha}}(\gamma,s,Z_0) & = & \int_{I_{Y_{i,0}}-\lambda_{i,0}(\gamma)} \varphi\left(g_{j,x}(\gamma,s,Z_0)\right)\varphi(x) \dfrac{\partial^2g_{j,x}}{\partial s \partial \gamma_{\alpha}}(\gamma,s,Z_0) \d x \nonumber \\
         & - & \int_{I_{Y_{i,0}}-\lambda_{i,0}(\gamma)}\varphi\left(g_{j,x}(\gamma,s,Z_0)\right)\varphi(x)g_{j,x}(\gamma,s,Z_0) \dfrac{\partial g_{j,x}}{\partial s}(\gamma,s,Z_0) \nonumber \\
         & & \hspace*{5cm} \times \dfrac{\partial g_{j,x}}{\partial\gamma_{\alpha}}(\gamma,s,Z_0) \d x.
		\end{eqnarray}
	Consequently, we can write the following partial derivatives
		\begin{equation}
		\dfrac{\partial^2}{\partial s^2} \log \left( \tilde{\ell}_{i,j}(\gamma,s,Z_0) \right) = \underbrace{\dfrac{\partial^2 \tilde{\ell}_{i,j}}{\partial s^2}(\gamma,s,Z_0) \times \dfrac{1}{\tilde{\ell}_{i,j}(\gamma,s,Z_0)}}_{I_{i,j}^{(1)}(\gamma,s,Z_0)} - \underbrace{\left(\dfrac{\partial \tilde{\ell}_{i,j}}{\partial s}(\gamma,s,Z_0) \times \dfrac{1}{\tilde{\ell}_{i,j}(\gamma,s,Z_0)}\right)^2}_{I_{i,j}^{(2)}(\gamma,s,Z_0)}
		\end{equation}				
	and
		\begin{align}
		\dfrac{\partial^2}{\partial s \partial \gamma_{\alpha}}\log \left( \tilde{\ell}_{i,j}(\gamma,s,Z_0) \right) & = \underbrace{\dfrac{\partial^2 \tilde{\ell}_{i,j}}{\partial s \partial \gamma_{\alpha}}(\gamma,s,Z_0) \times \dfrac{1}{\tilde{\ell}_{i,j}(\gamma,s,Z_0)}}_{I_{i,j}^{(3)}(\gamma,s,Z_0)} \nonumber \\ 
		& \hspace*{1cm}- \underbrace{\dfrac{\partial \tilde{\ell}_{i,j}}{\partial \gamma_{\alpha}}(\gamma,s,Z_0)\times \dfrac{\partial \tilde{\ell}_{i,j}}{\partial s}(\gamma,s,Z_0) \times \dfrac{1}{\tilde{\ell}_{i,j}(\gamma,s,Z_0)^2}}_{I_{i,j}^{(4)}(\gamma,s,Z_0)}.
		\end{align}		
	Our aim is now to bound uniformly on $\Theta$ each quantity $I_{i,j}^{(k)}(\gamma,s,Z_0), \ 1 \leqslant k \leqslant 4$, by a random variable that is integrable.
	
	Using Lemma \ref{lem_exp_ineq} and the upper-bounds for $g_{j,x}$ and its partial derivatives, we can show that the partial derivative \eqref{eq_2nd_partial_derivative} can be bounded, up to a constant, by an integral $I$ defined by 
		\begin{equation}
		I=\int_{I_{Y_{i,0}}-\lambda_{i,0}(\gamma)}\Phi\left(g_{j,x}(\gamma,s,Z_0)\right)\varphi(x)\left(1+\left|x\right|^q+ \|X_{-1}\|_2^q\right) \d x
		\end{equation}
	for some positive integer $q$. Now we apply the H\"older inequality, with $v^{-1}+w^{-1}=1$ and $v>1$ being determined latter. We have 
		\begin{equation}
		I\leqslant \left(\widetilde{\ell}_{i,j}(\gamma,s,Z_0)\right)^{1/v}\left(\int_{\R} \varphi(x)\left(1+\vert x\vert^q\right)^w \d x\right)^{1/w} + \| X_{-1}\|_2^q \widetilde{\ell}_{i,j}(\gamma,s,Z_0).
		\end{equation}
	Using Lemma \ref{lem_low_bound_exp}, there also exist $C_1,C_2>0$ such that
		\begin{equation}
		\underset{\theta\in\Theta}{\inf} \ \widetilde{\ell}_{i,j}(\gamma,s,Z_0)\geqslant C_1\exp\left(-C_2 \|X_{-1}\|_2^2\right).
		\end{equation} 
	By choosing $v$ sufficiently close to one and using our exponential moment assumption, we deduce that $I_{i,j}^{(1)}(\gamma,s,Z_0)$ is bounded uniformly on $\Theta$ by an integrable random variable. 
	
	One can show the same way that each quantity $I_{i,j}^{(k)}(\gamma,s,Z_0), \ 2 \leqslant k \leqslant 3$, is also bounded uniformly on $\Theta$ by an integrable random variable.
	
	We then easily deduce \eqref{eq_two_moments}, and assumption \emph{\textbf{I6}} is satisfied.
	\item The condition
		\begin{equation}
		\E \left( \| \nabla_2 \ell_{0}(\theta_0)\|^2 \right) < +\infty
		\end{equation}
	can be obtained using similar arguments as for getting \eqref{eq_two_moments}, so \emph{\textbf{I5}} is satisfied.
	\item We then prove that $ L_{2,2}=\E \left( \dfrac{\partial^2 \ell_{0}(\theta_0)}{\partial r^2} \right)$ is invertible.
				
	Since the different terms in the sum $\ell_{0} (\theta_0)$ depend on different parameters, $L_{2,2}$ contains only its diagonal terms
		\begin{equation}
		\E \left( \dfrac{\partial^2 \ell_{0}(\theta_0)}{\partial r_{i,j}^2} \right) = \dfrac{\partial^2}{\partial s^2} \left\lbrace \log \left( \tilde{\ell}_{i,j}(\gamma_0,s,Z_0) \right) \right\rbrace_{\vert s=r_{0,i,j}}, \quad 1\leqslant i < j \leqslant k.
		\end{equation}
	It is thus only necessary to check that for any $1 \leqslant i < j \leqslant k$,
		\begin{align}
		s(i,j)=\E \left[ \dfrac{1}{\tilde{\ell}_{i,j}(\gamma_0,r_{0,i,j},Z_0)^2} \left( \dfrac{\partial^2 \tilde{\ell}_{i,j}(\gamma_0,r_{0,i,j},Z_0)}{\partial s^2} \times \tilde{\ell}_{i,j}(\gamma_0,r_{0,i,j},Z_0) - \left(\dfrac{\partial \tilde{\ell}_{i,j}(\gamma_0,r_{0,i,j},Z_0)}{\partial s} \right)^2 \right) \right] < 0.
		\end{align}
	Since $\tilde{\ell}_{i,j}(\gamma_0,r_{0,i,j},Z_0)$ is the density of the pair $(Y_{i,0},Y_{j,0})$ conditionally on $\mathcal{F}_{-1}$, evaluated at point $(Y_{i,0},Y_{j,0})$, we have
		\begin{equation}
		\E \left[ \dfrac{\partial^2 \tilde{\ell}_{i,j}(\gamma_0,r_{0,i,j},Z_0)}{\partial s^2} \right] =0,
		\end{equation}
	hence
		\begin{equation}
		s(i,j) = -\E \left[ \dfrac{1}{\tilde{\ell}_{i,j}(\gamma_0,r_{0,i,j},Z_0)^2} \times \left(\dfrac{\partial \tilde{\ell}_{i,j}(\gamma_0,r_{0,i,j},Z_0)}{\partial s} \right)^2 \right] \leqslant 0.
		\end{equation}
	However, if $s(i,j)=0$, we have $\dfrac{\partial \tilde{\ell}_{i,j}(\gamma_0,r_{0,i,j},Z_0)}{\partial s} = 0$ a.s. and when $Y_{i,0}=Y_{j,0}=0$ we obtain
		\begin{equation}
		\dfrac{\partial}{\partial s} \left( s \longmapsto \int_{-\infty}^{-\lambda_{i,0}} \Phi \left( -\dfrac{\lambda_{j,0}+s x}{\sqrt{1-s^2}} \right) \varphi (x) \d x \right)_{\vert s=r_{0,i,j}} = 0 \quad \mathrm{a.s.}
		\end{equation}
	But we already mention in Lemma \ref{lem_increasing_f} that such a derivative is always positive. Thus $s(i,j) <0$ and $L_{2,2}$ is invertible. Thus, assumption \emph{\textbf{I7}} is satisfied.
		\item We finally prove the following martingale difference property \textbf{\emph{I8}}
		\begin{equation}\label{eq_mart_ell_1}
		\forall t \in \Z, \ \E \left( \nabla_2 \ell_{t}(\theta_0) \ \mid \mathcal{F}_{t-1} \right)=0.
		\end{equation}
	Using the representation of $\ell_{t}(\theta_0)$ as a sum of terms depending on different parameters, we only have to prove that for $1 \leqslant i < j \leqslant k$
		\begin{equation}\label{eq_mart_ell_2}
		\E \left( \dfrac{\dfrac{\partial \tilde{\ell}_{i,j}}{\partial s}(\gamma_0,r_{0,i,j},Z_t)}{\tilde{\ell}_{i,j}(\gamma_0,r_{0,i,j},Z_t)} \bigm\mid \mathcal{F}_{t-1} \right) = 0.
		\end{equation}
	Note that $\tilde{\ell}_{i,j}(\gamma_0,r_0,Z_t)$ is simply the conditional density of the pair $(Y_{i,t},Y_{j,t})$ conditionally to $\mathcal{F}_{t-1}$. Equality \eqref{eq_mart_ell_2} thus holds in a similar way as in the previous points, and assumption \emph{\textbf{I8}} is satisfied.
	\end{enumerate}
	
One can finally state that $\sqrt{T} \left( \hat{\gamma}-\gamma_0 , \hat{r}-r_0 \right)$ has a Gaussian limiting distribution by plugging the result for estimator $\hat{\gamma}$ in equation \eqref{eq_stoc_exp_3}. It yields, using the notation $m_{\gamma}(Z_t)$ defined previously
	\begin{equation}
	\sqrt{T}\left( \hat{r} - r_0 \right) = -\dfrac{L_{2,2}^{-1}}{\sqrt{T}} \sum_{t=p+1}^T \left\lbrace \nabla_2 \ell_{t}(\theta_0) - L_{i,2}J^{-1} \underbrace{\sum_{i=1}^k \dot{h}_{Y_{i,t}}(\lambda_{i,t}(\gamma_0))\nabla_1\lambda_{i,t}(\gamma_0)}_{\nabla_1m_{\gamma_0}(Z_t)} \right\rbrace + o_{\P}(1).
	\end{equation}
Hence
\begin{small}
	\begin{align}
	\sqrt{T}\left( \hat{\gamma}-\gamma_0,\hat{r}-r_0 \right) & = \dfrac{1}{\sqrt{T}} \sum_{t=p+1}^T \left( J^{-1}\nabla_1 m_{\gamma_0}(Z_t) , -L_{2,2}^{-1} \left(\nabla_2 \ell_{t}(\theta_0)-L_{1,2} J^{-1} \nabla_1 m_{\gamma_0}(Z_t)\right) \right) + o_{\P}(1) \nonumber \\
	& = \left(
			\begin{array}{c;{2pt/1pt}c}
			J^{-1} & 0 \\ \hdashline[2pt/2pt]
			0 & -L_{2,2}^{-1}
			\end{array}\right) \cdot \dfrac{1}{\sqrt{T}} \underbrace{\left(\begin{array}{c}
																	\nabla_1 m_{\gamma_0}(Z_t) \\ \hdashline[2pt/2pt]
																	\nabla_2 \ell_{t}(\theta_0) - L_{1,2}J^{-1}\nabla_1 m_{\gamma_0}(Z_t) 																\end{array} \right)}_{W_t}+o_{\P}(1).
	\end{align}
\end{small}
Yet, it has already been established that $\left( \nabla_1 m_{\gamma_0}(Z_t) \right)_{t\in \Z}$ is a martingale difference series relatively to $\left( \mathcal{F}_t \right)_{t\in \Z}$. Then, by equality \eqref{eq_mart_ell_1}, $(W_t)_{t\in \Z}$ is also a martingale difference series relatively to $\left( \mathcal{F}_{t}\right)_{t \in \Z}$. Since we have already showed that $\E \left( \| W_0 \|^2 \right) < +\infty$, Theorem \ref{theo_ibragimov_billingsley} then ensures the Gaussian limiting distribution of $\dfrac{W_t}{\sqrt{T}}$, and by Slutsky's lemma, the same conclusion holds for $\sqrt{T} \left( \hat{\gamma}-\gamma_0,\hat{r}-r_0 \right)$.

	\begin{lem}\label{lem_low_bound_exp}
	Take $a$ and $b$ two real numbers. Setting $\sigma = (1+4a^2)^{-1/2}$, we have the lower bound
		\begin{equation}
		\int_{I_{Y_{i,0}}-\lambda_{i,0}(\gamma)} \Phi (ax+b) \varphi(x) \d x \geqslant d^2 \sigma \exp \left( -2b^2 - \dfrac{\lambda_{i,0}(\gamma)^2}{\sigma^2} \right)
		\end{equation}
	where $d$ is the positive constant mentioned in Lemma \ref{lem_exp_ineq}.
	\end{lem}
	
	\begin{lemproof}[ of Lemma \ref{lem_low_bound_exp}]
	By Lemma \ref{lem_exp_ineq}, there exists a positive constant $d>0$ such that for all $x \in \R$
		\begin{equation}
		\Phi(ax+b) \geqslant d \exp (-(ax+b)^2) \geqslant \exp (-2a^2x^2-2b^2).
		\end{equation}
	We then get
		\begin{eqnarray}
		\int_{I_{Y_{i,0}}-\lambda_{i,0}(\gamma)} \Phi (ax+b) \varphi (x) \d x & \geqslant & d \exp(-2b^2) \int_{I_{Y_{i,0}}-\lambda_{i,0}(\gamma)} \dfrac{1}{\sqrt{2\pi}} \exp \left( - \dfrac{x^2}{2\sigma^2} \right) \d x \nonumber \\
		& \geqslant & d \exp(-2b^2) \sigma \min \left\lbrace \Phi \left( - \dfrac{\lambda_{i,0}(\gamma)}{\sigma} \right), 1-  \Phi \left( - \dfrac{\lambda_{i,0}(\gamma)}{\sigma} \right) \right\rbrace \nonumber \\
		& \geqslant & d^2 \exp (-2b^2) \sigma \exp \left(-\dfrac{\lambda_{i,0}(\gamma)^2}{\sigma^2} \right).
		\end{eqnarray}
	\end{lemproof}

\subsection{Proof of Proposition \ref{prop_law_large_numbers_panel}}

Our aim here is to apply Theorem \ref{theo_giap} from \citet{giap2016multidimensional}. 
For convenience, we also define $(Z_{j,t})_{t\in\Z}$ as independent copies of $(Z_{1,t})_{t\in \Z}$ for any negative $j$. We also set $Z:=(Z_{j,t})_{(j,t)\in \Z^2}$.

The probability space $\Omega$ we work with here is $\R^{\Z \times \Z}$, embedded with the sigma field denoted $\R^{\Z \otimes \Z}$ and generated by the cylinders of the form
	\begin{equation}
	C = \prod_{j,t \in \Z} C_{j,t},
	\end{equation}
where $C_{j,t} \in \mathcal{B}(R)$, such that there exist $\mathcal{J}, \mathcal{T} \subset \Z$ finite and
	\begin{equation}
	\forall (j,t) \notin \mathcal{J}\times \mathcal{T}, \ C_{j,t} = \R.
	\end{equation}
The transformations of $\R^{\Z \times \Z}$ we consider here are
	\begin{equation}
	\tau_1 : \begin{array}[t]{l}
	\R^{\Z \times \Z} \longrightarrow \R^{\Z \times \Z} \\
	\left(x_{j,t}\right)_{j,t \in \Z} \longmapsto \left( x_{j+1,t} \right)_{j,t\in \Z}
	\end{array}
	\quad \text{and} \quad
	\tau_2 : \begin{array}[t]{l}
	\R^{\Z \times \Z} \longrightarrow \R^{\Z \times \Z} \\
	\left(x_{j,t}\right)_{j,t \in \Z} \longmapsto \left( x_{j,t+1} \right)_{j,t\in \Z}.
	\end{array}
	\end{equation}
Finally, the probability measure we consider is $\P_{Z}$
	\begin{equation}
	\forall A = \left(A_{j,t} \right)_{j,t\in \Z} \in \R^{\Z \otimes \Z}, \ \P_{Z} (A) = \P \left( \cap_{j,t \in \Z} Z_{j,t}^{-1}(A_{j,t}) \right).
	\end{equation}

	\begin{enumerate}[label=$\triangleright$,labelindent=0pt,labelwidth=!,wide]
	\item We first check that $\tau_1$ and $\tau_2$ preserve $\P_Z$. To this end, take any $A \in \R^{\Z \otimes \Z}$, we have
		\begin{align}
		\P_Z \left(\tau_1^{-1} (A) \right) & = \P \left( \cap_{j \in \Z} \underbrace{\cap_{t\in\Z} \left\lbrace Z_{j,t} \in A_{j+1,t}\right\rbrace}_{\in \mathcal{G}_j} \right) \nonumber \\
		& = \prod_{j \in \Z} \P \left( \cap_{t\in \Z} \left\lbrace Z_{j,t} \in A_{j+1,t} \right\rbrace \right) \nonumber
		\shortintertext{by independence of the $\left(\mathcal{G}_j \right)_{j \in \Z}$}
		& = \prod_{j\in \Z} \P \left( \cap_{t \in\Z} \left\lbrace Z_{j+1,t} \in A_{j+1,t} \right\rbrace \right) \nonumber
		\shortintertext{since the distribution of $\left(Z_{j,t}\right)_{t\in \Z}$ does not depend on $j$}
		& = \P_Z(A).
		\end{align}
	Hence $\tau_1$ preserves $\P_Z$. On the other hand
		\begin{align}
		\P_Z\left(\tau_2^{-1} \right) & = \P \left( \cap_{j \in \Z} \cap_{t\in \Z} \left\lbrace Z_{j,t} \in A_{j,t+1} \right\rbrace \right) \nonumber \\
		& = \prod_{j \in \Z} \P \left( \cap_{t \in\Z} \left\lbrace Z_{j,t} \in A_{j,t+1} \right\rbrace \right) \nonumber \\
		& = \prod_{j \in \Z} \P \left( \cap_{t\in\Z} \left\lbrace Z_{j,t+1} \in A_{j,t+1} \right\rbrace \right)
		\nonumber
		\shortintertext{by stationarity of $(Z_{j,t})_{t\in \Z}$}
		& = \P_Z(A),
		\end{align}
	so $\tau_2$ also preserves $\P_Z$.
	\item Setting $h:(\R^{d'})^{\Z\times \Z}\rightarrow \R^{d'}$ the mapping defined by $h(z)=z_{0,0}$, one can apply Theorem \ref{theo_giap} to the mapping $f\circ h$ to get
		\begin{equation}
		\lim_{N,T \to +\infty} \dfrac{1}{NT} \sum_{j=1}^N \sum_{t=1}^T f\circ h\circ\tau_1^j\tau_2^t = \E_Z\left(f\circ h \mid \mathcal{I}_{\tau} \right) \quad \P_Z-\mathrm{a.s.}
		\end{equation}
	where $\mathcal{I}_{\tau} = \mathcal{I}_{\tau_1} \cap \mathcal{I}_{\tau_2}$ and for $i \in \left\lbrace 1,2 \right\rbrace, \ \mathcal{I}_{\tau_i} = \left\lbrace A \in \R^{\Z \otimes \Z} \mid \tau_i^{-1}(A) = A \right\rbrace$.
	\item In order to obtain the desired result, it is then sufficient to prove that $\tau_1$ is ergodic, i.e.
		\begin{equation}
		\forall A \in \mathcal{I}_{\tau_1}, \ \P_Z(A) \in \left\lbrace 0,1\right\rbrace.
		\end{equation}
	We actually show that for any $A \in \mathcal{I}_{\tau_1}$
		\begin{equation}\label{eq_I_tau_1}
		\lim_{n \to +\infty} \P_Z (A \cap \tau_1^{-n}(A)) = \P_Z(A)^2.
		\end{equation}
	First, this latter equality is satisfied for any cylinder of the form
		\begin{equation}
		B = \prod_{j,t\in \Z} B_{j,t},
		\end{equation}
	with for all $|j| > J$
		\begin{equation}
		B_{j,t} = \R, \ \forall t \in \Z,
		\end{equation}
	for some fixed $J \in \N$. Indeed, for all $n \geqslant 2J+1$, we have
		\begin{align}
		\P_Z \left( B \cap \tau_1^{-n}(B) \right) & = \P \left( \cap_{j \in \Z} \cap_{t \in \Z} \left\lbrace Z_{j,t} \in B_{j,t} \right\rbrace \bigcap \cap_{j \in \Z} \cap_{t \in \Z} \left\lbrace Z_{j+n,t} \in B_{j,t} \right\rbrace \right) \nonumber \\
		& = \P \left( \cap_{j=-J}^J \underbrace{\cap_{t\in\Z} \left\lbrace Z_{j,t} \in B_{j,t} \right\rbrace}_{\in \mathcal{G}_j} \bigcap \cap_{j'=n-J}^{n+J} \underbrace{\cap_{t\in \Z} \left\lbrace Z_{j',t} \in B_{j'-n,t} \right\rbrace}_{\in \mathcal{G}_{j'}} \right).
		\end{align}
	Yet, the indices $j$ and $j'$ are necessarily different since $ n \geqslant 2J+1$, so by independence of the $\left( \mathcal{G}_j \right)_{j \in \Z}$
		\begin{align}\label{eq_invariance_tau_1}
		\P_Z \left( B \cap \tau_1^{-n}(B) \right) & = \P_Z(B) \times \prod_{j'=n-J}^{n+J} \P \left( \cap_{t\in \Z} \left\lbrace Z_{j',t} \in B_{j'-n,t} \right\rbrace \right) \nonumber \\
		& = \P_Z (B) \times \prod_{j=-J}^J \P \left( \cap_{t \in \Z} \left\lbrace Z_{j+n,t} \in B_{j,t} \right\rbrace \right) \nonumber \\
		& = \P_Z (B) \times \prod_{j=-J}^J \P \left( \cap_{t \in \Z} \left\lbrace Z_{j,t} \in B_{j,t} \right\rbrace \right) \nonumber
		\shortintertext{since the distribution of $(Z_{j,t})_{t\in\Z}$ does not depend on $j$}
		& = \P_Z(B)^2.
		\end{align}
	We now prove equality \eqref{eq_I_tau_1}, and take any $A \in \mathcal{I}_{\tau_1}$ and $\varepsilon >0$. There exists a cylinder $B$ such that $\P_Z(A \Delta B) <\varepsilon$ and we have
		\begin{align}
		\left| \P_Z \left( A \cap \tau_1^{-n}(A) \right) - \P_Z(A)^2 \right| & \leqslant \left| \P_Z \left( A \cap \tau_1^{-n}(A) \right) - \P_Z\left(B \cap \tau_1^{-n}(A) \right) \right| \nonumber \\
		& \hspace*{2cm} + \left|\P_Z\left(B \cap \tau_1^{-n}(A) \right) - \P_Z \left( B \cap \tau_1^{-n}(B) \right) \right| \nonumber \\
		& \hspace*{2cm} + \left|\P_Z \left( B \cap \tau_1^{-n}(B) \right) - \P_Z(B)^2 \right| \nonumber \\
		& \hspace*{2cm} + \left| \P_Z(B)^2 - \P_Z(A)^2 \right|.
		\end{align}
	However
		\begin{equation}
		\left|\P_Z(B)^2-\P_Z(A)^2 \right| \leqslant 2 \times \left|\P_Z(B)-\P_Z(A) \right| \leqslant 2 \varepsilon
		\end{equation}
	and for $n$ large enough, we have according to equality \eqref{eq_invariance_tau_1}
		\begin{equation}
		\left|\P_Z \left( B \cap \tau_1^{-n}(B) \right) - \P_Z(B)^2 \right| =0.
		\end{equation}
	For such a $n$, we	also have since $A \in \mathcal{I}_{\tau_1}$
		\begin{align}
		\left| \P_Z \left( A \cap \tau_1^{-n}(A) \right) - \P_Z\left(B \cap \tau_1^{-n}(A) \right) \right| & = \left| \P_Z(A) - \P_Z(B \cap A)\right| \nonumber \\
		& = \P_Z(A \setminus B) \nonumber \\
		& < \varepsilon.
		\end{align}
	Finally, we have
		\begin{align}
		\left|\P_Z\left(B \cap \tau_1^{-n}(A) \right) - \P_Z \left( B \cap \tau_1^{-n}(B) \right) \right| & \leqslant \P_Z \left( \tau_1^{-n}(A) \Delta \tau_1^{-n}(B) \right) \nonumber \\
		& = \P_Z \left( \tau_1^{-n}\left( A \Delta B \right) \right) \nonumber \\
		& = \P_Z(A \Delta B) \quad \text{since $\tau_1$ preserves $\P_Z$} \nonumber \\
		& < \varepsilon.
		\end{align}
	Thus, for $n \in \N$ large enough
		\begin{equation}
		\left| \P_Z \left( A \cap \tau_1^{-n}(A) \right) - \P_Z(A)^2 \right| < 4 \varepsilon,
		\end{equation}
	hence
		\begin{equation}
		\lim_{n \to +\infty} \P_Z \left( A \cap \tau_1^{-n}(A) \right) = \P_Z(A) = \P_Z(A)^2.
		\end{equation}
	We then easily deduce that $\P_Z(A) \in \left\lbrace 0,1 \right\rbrace$, and we have the desired result.
	\end{enumerate}

\subsection{Proof of Proposition \ref{completion}}
As in the proof of Proposition \ref{prop_law_large_numbers_panel}, we may assume that $Z_{j,t}$ is also defined for a non-positive integer $j$. To this end, it is simply necessary to define $(V_j,\kappa_{j,t},\varepsilon_{j,t})_{t\in\Z}$ for non-positive indices $j$, independently across the $j$ and with the same distribution as $(V_1,\kappa_{1,t},\varepsilon_{1,t})_{t\in\Z}$. 
Set $\widetilde{\kappa}_{j,t}=\left(\eta_t,\kappa_{j,t},\varepsilon_{j,t}\right)$ for $(j,t)\in\N\times \Z$. It is quite clear that there exists a measurable mapping $\overline{G}:\R^{d_2}\times \left(\R^{d_4+d_6+k}\right)^{\N}\rightarrow \R$ such that
$$Z_{j,t}=\overline{G}\left(V_j,\left(\widetilde{\kappa}_{j,t-s}\right)_{s\geq 0}\right).$$
From our independence conditions, it is straightforward to show that the probability distribution of $Z:=\left(Z_{j,t}\right)_{(j,t)\in\Z^2}$ is invariant under the two shifts $\tau_1$ and $\tau_2$ defined in the proof of Proposition \ref{prop_law_large_numbers_panel}. We are going to show that for all Borel subset $A$ of $\R^{\Z\times \Z}$, then
\begin{equation}\label{ob}
\lim\sup_{m\rightarrow \infty}\lim\sup_{n\rightarrow \infty}\P_Z\left(A\cap \tau_1^{-m}\tau_2^{-n}A\right)=\P_Z(A)^2.
\end{equation}
Applying \eqref{ob} to $A\in \mathcal{I}_{\tau_1}\cap \mathcal{I}_{\tau_2}$ an invariant set, we will get $\P_Z(A)\in\{0,1\}$.
We will prove \eqref{ob} when $A$ is cylinder set. The rest of the proof is similar to that of Proposition \ref{prop_law_large_numbers_panel}. For $(j,t)\in\Z^2$ and $v_j\in \R^{d_2}$, set
$$Z_{j,t}(v_j)=\overline{G}\left(v_j,\left(\widetilde{\kappa}_{j,t-s}\right)_{s\geq 0}\right).$$
If $A=\prod_{(j,t)\in\Z^2}A_{j,t}$ with $A_{j,t}=\R$ except if $(j,t)\in J\times K$, for some $J$ and $K$ some finite subsets of $\Z$, we have
$$\P_Z\left(A\cap \tau_1^{-m}\tau_2^{-n} A\right)=\P\left(\cap_{(j,t)\in J\times K}\left\{Z_{j,t}\in A_{j,t},Z_{j+m,t+n}\in A_{j,t}\right\}\right).$$
Conditioning with respect to $V_j=v_j$ for $j\in\Z$, we get for a fixed value of $m$,
\begin{eqnarray*}
&&\lim_{n\rightarrow \infty}\P\left(\cap_{(j,t)\in J\times K}\left\{Z_{j,t}(v_j)\in A_{j,t},Z_{j+m,t+n}(v_{j+m})\in A_{j,t}\right\}\right)\\
&=&\P\left(\cap_{(j,t)\in J\times K}\left\{Z_{j,t}(v_j)\in A_{j,t}\right\}\right)\P\left(\cap_{(j,t)\in J\times K}\left\{Z_{j,t}(v_{j+m})\in A_{j,t}\right\}\right).
\end{eqnarray*}
This property is due to the mixing property of the Bernoulli shifts. See for instance \citet{samorodnitsky2016stochastic}, Corollary $2.2.5.$
Indeed, setting $J_m=J\cup\{m+j: j\in J\}$, $x_{J_m}=\left(x_j\right)_{j\in J_m}$ for $x\in \R^{\Z}$, and for $t\in \Z$, $Z_{J_m,t}=\left(Z_{j,t}(v_j)\right)_{j\in J_m}$, $B_{j,t}=A_{j-m,t}$ if $j-m\in J$ and $B_{j,t}=\R$ if $j-m\notin J$, we have the equality
$$\cap_{(j,t)\in J\times K}\left\{Z_{j,t}(v_j)\in A_{j,t},Z_{j+m,t+n}(v_{j+m})\in A_{j,t}\right\}=\cap_{t\in K}\left\{Z_{J_m,t}\in \prod_{j\in J_m}A_{j,t},Z_{J_m,t+n}\in \prod_{j\in J_m}B_{j,t}\right\}.$$
But since $Z_{J_m,t}\in \sigma\left(\widetilde{\kappa}_{J_m,t-s}: s\geq 0\right)$, one can use the mixing property for this multivariate Bernoulli shift. 

Now if $m$ is large enough, we have $j\neq j'+m$ for any pair $(j,j')\in J^2$. One can integrate the previous equality with respect to the distribution of the $V_j's$ and then use the dominated convergence theorem as well as the independence properties to get
\begin{eqnarray*}
&&\lim_{n\rightarrow \infty}\P\left(\cap_{(j,t)\in J\times K}\left\{Z_{j,t}\in A_{j,t},Z_{j+m,t+n}\in A_{j,t}\right\}\right)\\
&=&\P\left(\cap_{(j,t)\in J\times K}\left\{Z_{j,t}\in A_{j,t}\right\}\right)^2.
\end{eqnarray*}
This yields to \eqref{ob} and completes the proof.$\square$

\subsection{Proof of Theorem $4.3$}

First we note that under Assumption {\bf C1} or {\bf C'1}
and Assumption {\bf C2}, Corollary $3.1$ guarantees the existence of a unique stationary solution $(Y_{j,t})_{t\in\Z}$ to the recursions.

All the asymptotic results follow by checking the assumptions of Theorem \ref{theo_cons_panel} and Theorem
\ref{theo_nor_panel} and their corollaries. Similar assumptions have been already checked for proving Theorem 
$3.2$ and no new arguments are necessary. Details are then omitted.$\square$

\subsection{Additional asymptotic results for general M-estimators}

This section consists in a collection of various results that are used in section 6. Proofs are omitted here.

We begin with a useful inequality about the CDF $\Phi$ of a Gaussian distribution $\mathcal{N}(0,1)$

	\begin{lem}\label{lem_exp_ineq}
	There exist two positive real numbers $c$ and $d$ such that for all $x \in \R$
		\begin{equation}
		\min \left\lbrace \Phi (x), 1- \Phi (x) \right\rbrace \geqslant \dfrac{c}{1+ |x|}\e^{-x^2/2} \geqslant d \e^{-x^2}.
		\end{equation}
	\end{lem}
	
%
%
%
	
\medskip

The following result is derived from a more general result in \citet{pfanzagl1969measurability}, and is about consistency. We propose an alternative proof for our setup.

	\begin{theo}\label{theo_pfanzagl}
	Let $(Z_t)_{t\in\Z}$ be a stochastic process valued in a measurable space $\mathcal{Z}$, whose distribution depends on an unknown parameter $\theta_0 \in \Theta$, where $\Theta$ is a set of $\R^n, \ n \in \N^*$.
		
		For all $\theta \in \Theta$, we consider a measurable application $m_{\theta} : \mathcal{Z} \longrightarrow \R$.
		
		Let us define the assumptions:
		
			\begin{enumerate}[label=\textbf{H\arabic*}]
			\item $(Z_t)_{t\in \Z}$ is ergodic.
			\item $\Theta$ is a compact set.
			\item For all $z \in \mathcal{Z}$, the application $\theta \longmapsto m_{\theta}(z)$ is continuous and $ \E \left(  \underset{\theta \in \Theta}{\sup} \ |m_{\theta}(Z_0)|\right) < +\infty$.
			\item For all $\theta \in \Theta$, $\E \left( m_{\theta}(Z_0) \right) \leqslant \E \left( m_{\theta_0} (Z_0)\right)$ with equality if and only if $\theta=\theta_0$.
			\end{enumerate}
		If assumptions \textbf{H1} to \textbf{H4} are satisfied, then
			\begin{equation}\label{eq_def_theta_hat}
			\hat{\theta}_T := \underset{\theta \in \Theta}{\argmax} \ \sum_{t=1}^T m_{\theta}(Z_t) \underset{T \to +\infty}{\longrightarrow} \theta_0 \quad \mathrm{a.s.}
			\end{equation}
	\end{theo}

The next result deals with asymptotic normality, and is very similar to Theorem 5.6.1 in \citet{straumann2006estimation}. Here, $\| . \|$ denotes any norm on the considered space, for example the infinite norm. Furthermore, for a generic function $f : \R^n \longrightarrow \R$, we denote $\nabla f$ the gradient vector of $f$ and $\nabla^2 f$ its hessian matrix.

	\begin{theo}\label{theo_straumann}
	
	Let us consider a process $(Z_t)_{t \in \Z}$ satisfying conditions \textbf{H1} to \textbf{H4}	 mentioned in theorem \ref{theo_pfanzagl}, and consider $(\mathcal{F}_t)_{t\in \Z}$ a filtration. We assume that for all $z \in \mathcal{Z}$, the application $\theta \longmapsto m_{\theta}(z)$ is of class $\mathcal{C}^2$ on $\Theta$.
	
	Assume also that
	
		\begin{enumerate}[label=\textbf{H\arabic*}]
		\setcounter{enumi}{4}
		\item $\theta_0$ belongs to the interior of $\Theta$.
		\item $ \E \left( \| \nabla m_{\theta_0}(Z_0) \|^2 \right) < +\infty$.
		\item $\E \left( \underset{\theta \in \Theta}{\sup} \ \| \nabla^2 m_{\theta}(Z_0)\| \right) < +\infty$.
		\item The matrix $\E \left( \nabla^2 m_{\theta_0} (Z_0) \right)$ is invertible.
		\item For all $t \in \Z$, $\E \left( \nabla m_{\theta_0}(Z_t) \mid \mathcal{F}_{t-1} \right)=0$, i.e.$\left( \nabla m_{\theta_0}(Z_t) \right)_{t\in \Z}$ is a martingale difference series relatively to $\left( \mathcal{F}_{t-1} \right)_{t\in \Z}$.
		\end{enumerate}
	Then, $\hat{\theta}_{T}:= \underset{\theta \in \Theta}{\argmax} \ \sum_{t=1}^T m_{\theta}(Z_t) \underset{n \to +\infty}{\longrightarrow} \theta_0$ almost surely, and we have the stochastic expansion
		\begin{equation} \label{eq_stoc_exp_gen}
		\E \left( \nabla^2 m_{\theta_0}(Z_0) \right) \cdot \sqrt{T}\left( \hat{\theta}_T - \theta_0 \right) = -\dfrac{1}{\sqrt{T}} \sum_{t=1}^T \nabla m_{\theta_0} (Z_t) +o_{\P}(1).
		\end{equation}
	This leads to the following asymptotic normality
		\begin{equation}
		\sqrt{n}\left( \hat{\theta}_n - \theta_0 \right) \overset{\mathcal{L}}{\longrightarrow} \mathcal{N}\left(0, \E \left( \nabla^2 m_{\theta_0}(Z_0) \right)^{-1} \cdot \mathrm{Var}\left( \nabla m_{\theta_0}(Z_0) \right) \cdot \E \left( \nabla^2 m_{\theta_0}(Z_0) \right)^{-1} \right).
		\end{equation}
	\end{theo}

Proposition 8 of \citet{debaly2023multivariate} states the following refinement of Theorem \ref{theo_pfanzagl}

	\begin{prop}\label{prop8_debaly}
	
	Take $\Gamma$ and $\mathcal{R}$ two compact sets of $\R^p$ and $\R^q$ respectively, and denote $\Theta = \Gamma \times \mathcal{R}$.
	
	Let $(Z_t)_{t\in \Z}$ be a stochastic process valued in a measurable space $\mathcal{Z}$, whose distribution depends on an unknown parameter $\theta_0 = (\gamma_0,r_0) \in \Theta$.
	
	For all $\theta = (\gamma,r) \in \Theta$, we consider a measurable application $m_{\theta} : \mathcal{Z} \longrightarrow \R$.
	
	Let us define the assumptions:
		\begin{enumerate}[label=\textbf{I\arabic*}]
		\item $(Z_t)_{t\in \Z}$ is ergodic.
		\item For all $z \in \mathcal{Z}$, the application $\theta \longmapsto m_{\theta}(z)$ is continuous and $ \E \left(  \underset{\theta \in \Theta}{\sup} \ |m_{\theta}(Z_0)|\right) < +\infty$.
		\item For all $r \in \mathcal{R}$, $\E \left( m_{(\gamma_0,r)}(Z_0) \right) \leqslant \E \left( m_{(\gamma_0,r_0)}(Z_0) \right)$ with equality if and only if $r=r_0$.
		\end{enumerate}
	Finally, we suppose that for any $T \in \N^*$, we have a strongly consistent estimator $\hat{\gamma}_T$ of $\gamma_0$
		\begin{equation}
		\hat{\gamma}_T \underset{T \to +\infty}{\longrightarrow} \gamma_0 \quad \mathrm{a.s.}
		\end{equation}
	Then, if assumptions \textbf{I1} to \textbf{I3} are satisfied, we have
		\begin{equation}
		\hat{r}_T := \underset{r \in \mathcal{R}}{\argmax} \ \sum_{t=1}^T m_{(\hat{\gamma}_T,r)}(Z_t) \underset{T \to +\infty}{\longrightarrow} r_0 \quad \mathrm{a.s.}
		\end{equation}
	\end{prop}

It is also possible to obtain a refinement of Theorem \ref{theo_straumann}, that we detail in the following result.

Here, for a generic function
	\begin{equation}
	f : \begin{array}[t]{l}
	\R^p \times \R^q \longrightarrow \R \\
	(x,y) \longmapsto f(x,y)
	\end{array}
	\end{equation}		
we denote by $\nabla_1 f$ and $\nabla_2 f$ the gradients of function $f$ obtained when we differentiate respectively with respect to $x$ and $y$. Similarly, we denote $\nabla_1^2 f$ and $\nabla_2^2$ the corresponding Hessian matrices, and we also denote 
	\begin{equation}
	\nabla_{1,2} f = \left( \dfrac{\partial^2 f}{\partial x_j \partial y_i} \right)_{ 1 \leqslant i \leqslant q, \ 1 \leqslant j \leqslant p}.
	\end{equation}
	\begin{prop} \label{prop_sto_exp_refine}
	Take $\Gamma$ and $\mathcal{R}$ two compact sets of $\R^p$ and $\R^q$ respectively, and denote $\Theta=\Gamma \times \mathcal{R}$.
	
	Let $(Z_t)_{t\in \Z}$ be a stochastic process satisfying conditions \textbf{I1} to \textbf{I3} mentioned in Proposition \ref{prop8_debaly}, and consider $\left( \mathcal{F}_t \right)_{t \in \Z}$ a filtration. Assume that for all $z \in \mathcal{Z}$, the application $\theta \longmapsto m_{\theta}(z)$ is of class $\mathcal{C}^2$ on $\Theta$. Assume also that
	
	\begin{enumerate}[label=\textbf{I\arabic*}]
		\setcounter{enumi}{3}
		\item $\theta_0$ belongs to the interior of $\Theta$.
		\item $\E \left( \| \nabla_2 m_{\theta_0}(Z_0)\|^2 \right)<+\infty$.
		\item $\E \left( \underset{\theta \in \Theta}{\sup} \ \|\nabla_{1,2} m_{\theta}(Z_0)\| \right) < +\infty$ and $\E \left( \underset{\theta \in \Theta}{\sup} \ \|\nabla_2^2 m_{\theta}(Z_0)\| \right) < + \infty$.
		\item The matrix $\E \left( \nabla_2^2 m_{\theta_0} (Z_0) \right)$ is invertible.
		\item For all $t \in \Z$, $\E \left( \nabla_2 m_{\theta_0}(Z_t) \mid \mathcal{F}_{t-1} \right)=0$, i.e.$\left( \nabla_2 m_{\theta_0}(Z_t) \right)_{t\in \Z}$ is a martingale difference series relatively to $\left( \mathcal{F}_{t-1} \right)_{t\in \Z}$.
		\end{enumerate}
	
	Finally, we suppose that for any $T \in \N^*$, we have a strongly consistent estimator $\hat{\gamma}_T$ of $\gamma_0$ such that
		\begin{equation}
		\sqrt{T} \left( \hat{\gamma}_T - \gamma_0\right) = O_{\P}(1).
		\end{equation}			
	Then $\hat{r}_T := \underset{r \in \mathcal{R}}{\argmax} \ \sum_{t=1}^T m_{(\hat{\gamma}_T,r)}(Z_t) \underset{T \to +\infty}{\longrightarrow} r_0$ almost surely, and we have the stochastic expansion
		\begin{equation}
		\E \left( \nabla_2^2 m_{\theta_0}(Z_0) \right) \cdot \sqrt{T}\left( \hat{r}_T - r_0 \right) = -\dfrac{1}{\sqrt{T}}\sum_{t=1}^T \nabla_2 m_{\theta_0}(Z_t) - \E \left( \nabla_{1,2} m_{\theta_0} \right) \sqrt{T} \left( \hat{\gamma}_T - \gamma_0 \right) + o_{\P}(1).
		\end{equation}
	\end{prop}

	The previous theorems dealing about asymptotic normality are largely based on the well-known Ibragimov-Billingsley's theorem for martingale differencies below. See \citet{billingsley2013convergence}, Theorem 35.12, for a proof. This theorem extends naturally to multivariate random variables.
	
	\begin{theo} \label{theo_ibragimov_billingsley}
	Consider an array $(Y_{k,t})_{k,t\in \N^*}$ an array of real random variables such that for each $k \in \N^*$, $(Y_{k,t})_{t \in \N^*}$ is a difference martingale series with respect to a filtration $\left( \mathcal{G}_{k,t} \right)_{t \in \N}$ (with $\mathcal{G}_{k,0} = \left\lbrace \emptyset , \Omega \right\rbrace$).
	
	Suppose the $Y_{k,t}$'s have second moments, and set
		\begin{equation}
		\sigma_{k,t}^2 = \E \left( Y_{k,t}^2 \mid \mathcal{G}_{k,t-1} \right).
		\end{equation}
	Assume that
		\begin{equation}
		\sum_{t=1}^{+\infty} \sigma_{k,t}^2 \overset{\P}{\longrightarrow} \sigma^2,
		\end{equation}
	where $\sigma$ is a positive constant, and that
		\begin{equation}
		\sum_{t=1}^{+\infty} \E \left( Y_{k,t}^2 \ind_{|Y_{k,t}| \geqslant \varepsilon} \right) \underset{k \to +\infty}{\longrightarrow} 0
		\end{equation}
	for each $\varepsilon>0$. Then
		\begin{equation}
		\sum_{t=1}^{+\infty} Y_{k,t} \overset{\mathcal{L}}{\longrightarrow} \mathcal{N}(0,\sigma^2)
		\end{equation}
	when $k \to +\infty$.
	\end{theo}
	
Let us finally recall Lemma 3.1 from \citet{giap2016multidimensional}, which is a fundamental result in our paper. 

	\begin{theo}\label{theo_giap}
	Consider $(\X, \| \cdot \|)$ a Banach space and let $k$ be a positive integer, and $\theta_1,\ldots,\theta_k$ be commutative measure preserving transformations of a probability space $ (\Omega, \mathcal{A}, \P)$.
	
	Then for any measurable application $f : \Omega \longrightarrow \X$ such that $\E \left( \| f \|  \left( \log^{+} \| f\| \right)^{k-1} \right) < + \infty$, the multiple sequence
		\begin{equation}
		A_{m_1,\ldots,m_k}f = \dfrac{1}{m_1\dots m_k} \sum_{i_1=0}^{m_1-1} \dots \sum_{i_k=0}^{m_k-1} f \left( \theta_1^{i_1}\dots \theta_k^{i_k} \right)
		\end{equation}
	converges almost surely to $\E \left( f \mid \mathcal{I}_{\theta} \right)$ as $\min (m_1,\ldots m_k) \to +\infty$, where $\mathcal{I}_{\theta} = \cap_{i=1}^k \mathcal{I}_{\theta_i}$ and \\ $\mathcal{I}_{\theta_i} = \left\lbrace A \in \mathcal{A} : \theta_i^{-1} (A) = A \right\rbrace$.
	
	Moreover, if $\theta_s$ is ergodic for some $s$ belonging to $\left\lbrace 1, \ldots,k\right\rbrace$, $\E \left( f \mid \mathcal{I}_{\theta} \right) = \E (f)$ almost surely.
	\end{theo}

\subsection{Asymptotic results for $M-$estimators with panel data}
We now state several general results relative to $M$-estimators for panel data. The two first ones deal with the consistency, while the others allow us to obtain central limit theorems.

	\begin{theo}\label{theo_cons_panel} 
	Let $(Z_{j,t})_{j\in \N^*,t \in \Z}$ be a family of random variables valued in $\mathcal{Z}=\R^{d'}$ for some positive integer $d'$, satisfying \eqref{bernshift} and either {\bf C1} or {\bf C'1}.   
	
	Additionally, for a parameter $\theta$ belonging to a compact set $\Theta \subset \R^{d_{\Theta}}$, let us consider a measurable function $m_{\theta} : \mathcal{Z} \longrightarrow \R$, and assume that:
		\begin{enumerate}[label=\textbf{P\arabic*}]
		\item For all $z \in \mathcal{Z}$, the application $ \theta \longmapsto m_{\theta} (z)$ is continuous.
		\item There exists $\delta >0$ such that $ \E \left( \underset{\theta \in \Theta}{ \sup } \ \left| m_{\theta} (Z_{1,1}) \right|^{1+\delta} \right) < +\infty$.                          
		\item There exists $\theta_0 \in \Theta$ such that for all $\theta \in \Theta$, $\E (m_{\theta} (Z_{1,1} )) \leqslant \E ( m_{\theta_0} (Z_{1,1}))$, with equality if and only if $\theta=\theta_0$.
		\end{enumerate}
	Then
		\begin{equation}
		\hat{\theta}_{n,T} := \underset{\theta \in \Theta}{\argmax} \ \sum_{t=1}^T \sum_{j=1}^n m_{\theta}(Z_{j,t}) \underset{n,T \to +\infty}{\longrightarrow} \theta_0 \quad \mathrm{a.s.}
		\end{equation}
	\end{theo}
	It is possible to refine Theorem \ref{theo_cons_panel} in the following case.

	\begin{cor}\label{cor_cons_panel}
	Let $(Z_{j,t})_{j\in \N^*,t \in \Z}$ be a family of random variables valued in $\mathcal{Z}=\R^{d'}$ for some positive integer $d'$, satisfying \eqref{bernshift} and either {\bf C1} or {\bf C'1}. We consider two compact sets $\Gamma$ and $\mathcal{R}$ of $\R^{d_{\Gamma}}$ and $\R^{d_{\mathcal{R}}}$ respectively, and we denote $\Theta = \Gamma \times \mathcal{R}$.
	
	Additionally, for a parameter $\theta = (\gamma,r) \in \Theta$, we consider a measurable function $m_{\theta} : \mathcal{Z} \longrightarrow \R$, and assume that:
		\begin{enumerate}[label=\textbf{Q\arabic*}]
			\item For all $z \in \mathcal{Z}$, the application $\theta \longmapsto m_{\theta}(z)$ is continuous.
		\item $\E \left( \underset{\theta \in \Theta}{\sup} \ \left| m_{\theta}(Z_{1,1}) \right|^{1+ \delta} \right) < +\infty$ for some $\delta >0$.
		\item There exists $\theta_0 = (\gamma_0,r_0) \in \Theta$ such that for all $r \in \mathcal{R}, \E \left( m_{\gamma_{0},r}(Z_{1,1}) \right) \leqslant  \E \left( m_{\gamma_{0},r_0}(Z_{1,1}) \right)$, with equality if and only if $r=r_0$.
		\end{enumerate}
	If $\hat{\gamma}_{n,T}$ is a statistic such that $\hat{\gamma}_{n,T} \underset{n,T \to +\infty}{\longrightarrow} \gamma_0$ almost surely, then
		\begin{equation}
		\hat{r}_{n,T} :=\underset{r \in \mathcal{R}}{\argmax} \ \sum_{t=1}^T\sum_{j=1}^n m_{\hat{\gamma}_{n,T}}(Z_{j,t}) \underset{n,T \to +\infty}{\longrightarrow} r_0 \quad \mathrm{a.s.}
		\end{equation}
	\end{cor}
	
In what follows, for any function $f_{\theta}$ valued in $\R$, depending on a parameter $\theta \in \R^{d_{\Theta}}$, we denote respectively $\nabla f_{\theta}$ and $\nabla^2 f_{\theta}$ the gradient and hessian matrix of $f_{\theta}$ when we differentiate with respect to $\theta$. We also denote by $\mathcal{A}\vee\mathcal{B}$ the sigma-field generated by two sigma-fields $\mathcal{A}$ and $\mathcal{B}$.

Furthermore, $\| \cdot \|$ denotes any norm on the considered vector space.

	\begin{theo}\label{theo_nor_panel}
	Let $(Z_{j,t})_{j\in \N^*,t \in \Z}$ be a family of random variables valued in $\mathcal{Z}=\R^{d'}$ for some positive integer $d'$, satisfying \eqref{bernshift} and either {\bf C1} or {\bf C'1}.
	
	Additionally, for a parameter $\theta$ belonging to a compact set $\Theta \subset \R^{d_{\Theta}}$, let us consider a measurable function $m_{\theta} : \mathcal{Z} \longrightarrow \R$ such that assumptions \textbf{P1} to \textbf{P3} are satisfied.
	
	Assume furthermore that
		\begin{enumerate}[label=\textbf{P\arabic*}]
		\setcounter{enumi}{3}
		\item $\theta_0$ belongs to the interior of $\Theta$.
		\item For all $z \in \mathcal{Z}, \ \theta \longmapsto m_{\theta}(z)$ is of class $\mathcal{C}^2$.
		
		\item $\E \left( \| \nabla m_{\theta_0}(Z_{1,1})\|^{2+\delta} \right) < +\infty$ for some $\delta >0$.
		\item $\E \left( \underset{\theta \in \Theta}{\sup} \ \| \nabla^2 m_{\theta}(Z_{1,1})\|^{1+ \delta} \right) < +\infty$ for some $\delta >0$.
		\item $\E \left(\nabla^2 m_{\theta_0}(Z_{1,1}) \right)$ is invertible.
		\item Denoting for $t\in\Z$, $\mathcal{H}_t=\vee_{\ell\geq 1}\mathcal{F}_{\ell,t}$ where the  $\mathcal{F}_{\ell,t}'$s are defined in \eqref{field}, then $\E\left(\nabla m_{\theta_0}(Z_{j,t})\vert \mathcal{H}_{t-1}\right)=0$ for $t\in \Z$ and $j\geq 1$.
		\end{enumerate}
	Then $\hat{\theta}_{n,T} := \underset{\theta \in \Theta}{\argmax} \ \sum_{t=1}^T \sum_{j=1}^n m_{\theta}(Z_{j,t}) \underset{n,T \to +\infty}{\longrightarrow} \theta_0$ almost surely and we have the stochastic expansion
		\begin{equation}
		\sqrt{nT} \ \E \left( \nabla^2 m_{\theta_0}(Z_{1,1})\right) \cdot \left( \hat{\theta}_{n,t}-\theta_0 \right) = -\dfrac{1}{\sqrt{nT}} \sum_{t=1}^T \sum_{j=1}^n \nabla m_{\theta_0}(Z_{j,t}) + o_{\P}(1).
		\end{equation}
	when $n,T \to \infty$. Therefore, we have the convergence in distribution
		\begin{equation}
		\sqrt{nT} \left(\hat{\theta}_{n,T} - \theta_0 \right) \overset{\mathcal{L}}{\longrightarrow} \mathcal{N}\left( 0 , \E \left( \nabla^2 m_{\theta_0}(Z_{1,1}) \right)^{-1} \cdot \mathrm{Var} \left( \nabla m_{\theta_0}(Z_{1,1}) \right) \cdot \E \left( \nabla^2 m_{\theta_0}(Z_{1,1}) \right)^{-1} \right)
		\end{equation}
	when $n,T \to \infty$.
	\end{theo}

Similarly to Corollary \ref{cor_cons_panel}, it is possible to refine Theorem \ref{theo_nor_panel} in the following case.

Just like in section \ref{sec_param_inf_single}, for a function $f_{(\gamma,r)}$ valued in $\R$, depending on a parameter $(\gamma,r) \in \Gamma \times \mathcal{R} \subset \R^{d_{\Gamma}} \times \R^{d_{\mathcal{R}}}$, we denote $\nabla_1 f$ and $\nabla_2 f$ the gradients obtained by differentiating $f$ with respect to $\gamma$ or $r$ respectively. We denote in a similar way $\nabla_1^2$ and $\nabla_2^2$ the corresponding hessian matrices and
	\begin{equation}
	\nabla_{1,2} f = \left( \dfrac{\partial^2 f}{\partial \gamma \partial r_i} \right)_{1 \leqslant i \leqslant d_{\mathcal{R}} , 1 \leqslant j \leqslant d_{\Gamma}} \in \R^{d_{\mathcal{R}} \times d_{\Gamma}}.
	\end{equation}

	\begin{cor}\label{cor_nor_panel}
 Let $(Z_{j,t})_{j\in \N^*,t \in \Z}$ be a family of random variables valued in $\mathcal{Z}=\R^{d'}$ for some positive integer $d'$, satisfying \eqref{bernshift} and either {\bf C1} or {\bf C'1}. We consider two compact sets $\Gamma$ and $\mathcal{R}$ of $\R^{d_{\Gamma}}$ and $\R^{d_{\mathcal{R}}}$ respectively, and we denote $\Theta = \Gamma \times \mathcal{R}$.
	
	Additionally, for a parameter $\theta =(\gamma,r) \in \Theta$, we consider a measurable function $m_{\theta} : \mathcal{Z} \longrightarrow \R$ such that assumptions \textbf{Q1} to \textbf{Q3} are satisfied.
	
	Assume furthermore that
		\begin{enumerate}[label=\textbf{Q\arabic*}]
		\setcounter{enumi}{3}
		\item $\theta_0$ belongs to the interior of $\Theta$.
		\item For all $z \in \mathcal{Z}, \theta \longmapsto m_{\theta}(Z_{1,1})$ is of class $\mathcal{C}^2$.
		\item $\E \left( \| \nabla_2 m_{\theta_0}(Z_{1,1}) \|^{2+\delta} \right) < +\infty$ for some $\delta >0$.
		\item $\E \left( \underset{\theta \in \Theta}{\sup} \ \| \nabla_{1,2}m_{\theta}(Z_{1,1}) \|^{1+\delta} \right) < +\infty$ for some $\delta >0$.
		\item $\E \left( \underset{\theta \in \Theta}{\sup} \ \| \nabla_{2}^2 m_{\theta}(Z_{1,1}) \|^{1+\delta} \right) < +\infty$ for some $\delta >0$.
		\item $\E \left( \nabla_2^2 m_{\theta_0}(Z_{1,1}) \right)$ is invertible.
		\item  Denoting for $t\in\Z$, $\mathcal{H}_t=\vee_{\ell\geq 1}\mathcal{F}_{\ell,t}$ where the  $\mathcal{F}_{\ell,t}'$s are defined in \eqref{field}, then $\E\left(\nabla_2 m_{\theta_0}(Z_{j,t})\vert \mathcal{H}_{t-1}\right)=0$ for $t\in \Z$ and $j\geq 1$.
		\end{enumerate}
	Finally, we suppose that for any $n,T$ in $\N^*$, we have a strongly consistent estimator $\hat{\gamma}_{n,T}$ of $\gamma_0$ such that
		\begin{equation}
		\sqrt{nT}\left( \hat{\gamma}_{n,T} - \gamma_0 \right) = O_{\P} (1).
		\end{equation}
	when $n,T \to \infty$.
	
	Then $\hat{r}_{n,T} := \underset{r \in \mathcal{R}}{\argmax} \ \sum_{t=1}^T \sum_{j=1}^n m_{\left(\hat{\gamma}_{n,T},r \right)}(Z_{j,t}) \underset{n,T \to +\infty}{\longrightarrow} r_0$ almost surely and we have the stochastic expansion
		\begin{multline}
		\sqrt{nT} \ \E \left( \nabla_2^2 m_{\theta_0}(Z_{1,1}) \right) \left( \hat{r}_{n,T}-r_0 \right) = -\dfrac{1}{\sqrt{nT}} \sum_{t=1}^T \sum_{j=1}^n \nabla_2 m_{\theta_0}(Z_{j,t}) \\ - \E \left( \nabla_{1,2}m_{\theta_0}(Z_{1,1}) \right) \sqrt{nT} \left( \hat{\gamma}_{n,T} - \gamma_0 \right) + o_{\P}(1).
		\end{multline}
	\end{cor}
	
	\begin{note}
	Once again, assumption \emph{\textbf{Q8}} implies assumption \emph{\textbf{Q1}}.
	\end{note}

\subsection{Proof of Theorem \ref{theo_cons_panel}}

The proof of Theorem \ref{theo_cons_panel} is actually an adaptation of the proof of Theorem \ref{theo_pfanzagl} and follows directly from the two following lemmas.

	\begin{lem}\label{lem_panel_unif_law_large_numbers}
	Let $(Z_{j,t})_{j \in \N^*,t \in \Z}$ be a family of $\mathcal{F}_{j,t}$-measurable random variables valued in $\R^{d'}$ for some positive integer $d'$ and satisfying \eqref{bernshift}. Suppose that either Assumption {\bf C1} or Assumption {\bf C'1} holds true.
	
	Additionally, for a parameter $\theta$ belonging to a compact set $\Theta \subset \R^{d_{\Theta}}, \ n \geqslant 1$, let us consider a measurable function $m_{\theta} : Z \longmapsto \R$, and assume that
		\begin{enumerate}[label=\alph*)]
		\item For all $z \in Z$, the application $\theta \longmapsto m_{\theta}(z)$ is continuous.
		\item For some $\delta >0, \ \E \left( \underset{\theta \in \Theta}{\sup} \ |m_{\theta}(Z_{1,1})|^{1+\delta} \right) < +\infty$.
	\end{enumerate}
	Then
		\begin{equation}
		\underset{\theta \in \Theta}{\sup} \ \left| \dfrac{1}{nT} \sum_{t=1}^T \sum_{j=1}^n m_{\theta}(Z_{j,t}) - \E \left( m_{\theta}(Z_{1,1}) \right)\right| \underset{n,T \to +\infty}{\longrightarrow} 0 \quad \mathrm{a.s.}
		\end{equation}
	\end{lem}	

	\begin{lemproof}[ of Lemma \ref{lem_panel_unif_law_large_numbers}]
	We set, for any $z \in Z$ and any $\Delta >0$
		\begin{equation}
		W_{\Delta}(z) := \sup \ \left\lbrace \left| m_{\theta + h}(z) - m_{\theta}(z) \right| : \theta \in \Theta, \ \theta+h \in \Theta, \ \|h\|_2 <\Delta \right\rbrace.
		\end{equation}	
	One can easily show that
		\begin{equation}
		\lim_{\Delta \to 0} \E \left( W_{\Delta}(Z_{1,1}) \right) = 0.
		\end{equation}
	Take any $\varepsilon >0$, according to this latter equality, there exists $\delta >0$ such that
		\begin{equation}
		\Delta \leqslant \delta \Longrightarrow \E \left( W_{\Delta}(Z_{1,1}) \right) < \varepsilon.
		\end{equation}
	In addition, $\Theta \subset \cup_{\theta \in \Theta} \ B(\theta,\delta)$, where $B(\theta,\delta)$ denotes the open Euclidean ball of center $\theta$ and radius $\delta$. So, by compactness of $\Theta$, there exists a finite set $I$ such that $\theta \subset \cup_{i \in I} B(\theta_i, \delta)$.
	
	For all $\theta \in \Theta$ and $n,T \in \N^*$, we denote 
		\begin{equation}
		M_{n,T}(\theta) = \dfrac{1}{nT} \sum_{j=1}^n \sum_{t=1}^T m_{\theta} (Z_{j,t}) \quad \text{and} \quad M(\theta) = \E \left( m_{\theta}(Z_{1,1})\right).
		\end{equation}
	For $\theta \in \Theta$, there exists $i\in I$ such that $\theta \in B(\theta_i,\delta)$ and for all $n,T \in \N^*$
		\begin{align}\label{ineq_MNT_triang}
		\left| M_{n,T}(\theta) - M(\theta) \right| & \leqslant \left| M_{n,T}(\theta) - M_{n,T}(\theta_i) \right| +  \left| M_{n,T}(\theta_i) - M(\theta_i) \right| +  \left| M(\theta_i) - M(\theta) \right| \nonumber \\
		& \leqslant \dfrac{1}{nT} \sum_{j=1}^n \sum_{t=1}^T \left| m_{\theta}(Z_{j,t})-m_{\theta_i}(Z_{j,t}) \right| + \left| M_{n,T}(\theta_i) - M(\theta_i) \right| + \E \left( W_{\delta}(Z_{1,1}) \right) \nonumber \\
		& \leqslant \dfrac{1}{nT} \sum_{j=1}^n \sum_{t=1}^T W_{\delta}(Z_{j,t}) + \underset{i \in I}{\max} \ \left| M_{n,T}(\theta_i)-M(\theta_i) \right| +  \E \left( W_{\delta}(Z_{1,1}) \right).
		\end{align}
		\begin{enumerate}[label=$\triangleright$]
		\item Using Assumption $b)$ and either Proposition \ref{prop_law_large_numbers_panel} or Proposition \ref{completion}, we have for all $i \in I$
			\begin{equation}\label{eq_MNT-M}
			\left| M_{n,T}(\theta_i)-M(\theta_i) \right| \underset{n,T \to +\infty}{\longrightarrow} 0 \quad \mathrm{a.s.}
			\end{equation}

		We thus deduce that
			\begin{equation}
			\underset{i \in I}{\max} \ \left| M_{n,T}(\theta_i)-M(\theta_i) \right| \underset{n,T \to +\infty}{\longrightarrow} 0 \quad \mathrm{a.s.}
			\end{equation}
		\item We then show that
			\begin{equation}\label{eq_conv_MNT}
			\dfrac{1}{nT} \sum_{j=1}^n \sum_{t=1}^T W_{\delta}(Z_{j,t}) \underset{n,T \to +\infty}{\longrightarrow} \E \left( W_{\delta}(Z_{1,1}) \right) \quad \mathrm{a.s.}
			\end{equation}
		We set here, for any $j \in \N^*$ and $t \in \Z$, $H_{j,t} = W_{\delta}(Z_{j,t})$. According to point \emph{b)}, for some $\kappa >0$, we have
			\begin{align}
			\E \left( \left| H_{1,1} \right|^{1+\kappa} \right) & = \E \left( W_{\delta}(Z_{1,1})^{1+\kappa} \right) \nonumber \\
			& \leqslant 2^{1+\kappa} \E \left( \underset{\theta \in \Theta}{\sup} \ \left| m_{\theta}(Z_{1,1})\right|^{1+ \kappa} \right) \nonumber \\
			& < +\infty.
			\end{align} 
		Thus, Proposition \ref{prop_law_large_numbers_panel} or Proposition \ref{completion} gives the convergence \eqref{eq_conv_MNT}. 
		\end{enumerate}
	Inequality \eqref{ineq_MNT_triang} then leads to
		\begin{equation}
		\underset{\theta \in \Theta}{\sup} \ \left| M_{n,T}(\theta)-M(\theta) \right| \leqslant \dfrac{1}{nT}\sum_{j=1}^n \sum_{t=1}^T W_{\delta}(Z_{j,t}) + \underset{i \in I}{\max} \ \left| M_{n,T}(\theta_i)-M(\theta_i) \right| + \E \left( W_{\delta}(Z_{1,1}) \right)
		\end{equation}
	and by \eqref{eq_MNT-M} and \eqref{eq_conv_MNT}, we have almost surely
		\begin{equation}
		\underset{\theta \in \Theta}{\sup} \ \left| M_{n,T}(\theta)-M(\theta) \right| \leqslant 2\E \left( W_{\delta}(Z_{1,1}) \right) < 2\varepsilon
		\end{equation}
	 for $n$ and $T$ large enough, hence
	 	\begin{equation}
	 	\underset{\theta \in \Theta}{\sup} \ \left| \dfrac{1}{nT} \sum_{j=1}^n \sum_{t=1}^T m_{\theta}(Z_{j,t}) - \E \left( m_{\theta}(Z_{1,1}) \right) \right| \underset{n,T \to +\infty}{\longrightarrow} 0 \quad \mathrm{a.s.}
	 	\end{equation}
	\end{lemproof}

	\begin{lem}\label{lem_trad_conv_compact}
	For some $i\in\{1,2\}$, let $(f_m)_{m \in \N^i}$ be a sequence of real-valued functions defined on a compact space $\Theta\subset\R^{d_{\Theta}}$, such that $(f_m)_{m \in \N}$ converges uniformly toward a function $f$ on $\Theta$
		\begin{equation}
		\underset{\theta \in \Theta}{\sup} \ \left| f_m(\theta) - f(\theta) \right| \underset{m \to +\infty}{\longrightarrow} 0.
		\end{equation}
	Assume in addition that functions $f_m, \ m \in \N^i$ are continuous, and that $f$ admits a unique maximum in $\theta_0 \in \Theta$. Then
		\begin{equation}
		\underset{\theta \in \Theta}{\argmax} \ f_m(\theta) \underset{m \to \infty}{\longrightarrow} \theta_0.
		\end{equation}
	\end{lem}
	
		\begin{lemproof}[ of Lemma \ref{lem_trad_conv_compact}]
  Let $\Vert\cdot\Vert$ be an arbitrary norm on $\R^{d_{\Theta}}$.
	For all $m \in \N^i$, let us denote $\hat{\theta}_m := \underset{\theta \in \Theta}{\argmax} \ f_m(\theta)$. Suppose that $\hat{\theta}_m$ does not converge to $\theta_0$.
 There then exists $\varepsilon>0$ and a sequence $(m_k)_{k\in\N}$ valued in $\N^i$, with each coordinate going to $\infty$ and such that $\Vert\hat{\theta}_{m_k}-\theta_0\Vert>\varepsilon$.
 Since $\Theta$ is compact, up to an extraction, one can assume that for some $\theta^*\neq \theta_0 \in \Theta$,
		\begin{equation}
		\hat{\theta}_{m_k} \underset{k \to \infty}{\longrightarrow} \theta^*.
		\end{equation}
	We then have for all $k \in \N$
		\begin{equation}
		\left| f_{m_k}(\hat{\theta}_{m_k}) - f(\theta^*) \right| \leqslant \left| f_{m_k}(\hat{\theta}_{m_k}) -f(\hat{\theta}_{m_k} ) \right| + \left| f(\hat{\theta}_{m_k} ) - f(\theta^*) \right|,
		\end{equation}
	and the right side of this latter inequality converges toward $0$ when $k \to \infty$ by uniform convergence and continuity. Thus
		\begin{equation}
		f_{m_k} \left( \hat{\theta}_{m_k} \right) \underset{k \to \infty}{\longrightarrow} f(\theta^*).
		\end{equation}
	Furthermore, for any $k \in \N$ we have
		\begin{equation}
		f_{m_k}\left( \hat{\theta}_{m_k} \right) \geqslant f_{m_k}(\theta_0),
		\end{equation}
	and considering the limits, it yields $f(\theta^*) \geqslant f(\theta_0)$. Since $\theta_0$ is the unique maximum of $f$ on $\Theta$, we deduce $\theta^*=\theta_0$. From this contradiction, we get the result.
	\end{lemproof}

\subsection{Proof of Theorem \ref{theo_nor_panel}}

Once again, the proof of Theorem \ref{theo_nor_panel} is an adaptation of the proof of Theorem \ref{theo_straumann}. We begin with a preliminary lemma.

	\begin{lem}\label{lem_proof_theo_nor_panel}
	We have the following convergence in distribution
		\begin{equation}
		\dfrac{1}{\sqrt{nT}}\sum_{t=1}^T \sum_{j=1}^n \nabla m_{\theta_0}(Z_{j,t}) \overset{\mathcal{L}}{\longrightarrow} \mathcal{N}\left( 0, \mathrm{Var}(\nabla m_{\theta_0}(Z_{1,1}) \right)
		\end{equation}
	when $n,T \to +\infty$.
	\end{lem}
	
	\begin{lemproof}[ of Lemma \ref{lem_proof_theo_nor_panel}]
	Let us denote for all $n,T$ in $\N^*$
		\begin{equation}
		S_{n,T} = \dfrac{1}{\sqrt{nT}} \sum_{t=1}^T \sum_{j=1}^n \nabla m_{\theta_0}(Z_{j,t}) \quad \text{and} \quad V = \mathrm{Var}(\nabla m_{\theta_0}(Z_{1,1})).
		\end{equation}
	Note that $V$ is well defined according to assumption \textbf{\emph{P6}}.
	
	It is sufficient to show that for any increasing subsequences of integers $\left( n_k \right)_{k\geqslant 1}$ and $\left( T_k \right)_{k\geqslant 1}$, we have 
		\begin{equation}\label{eq_twosub}
		S_{n_k,T_k} \overset{\mathcal{L}}{\longrightarrow} \mathcal{N}(0,V)
		\end{equation}
	when $k\to +\infty$.
	
	Take indeed a function $g$ that is continuous and bounded, we want to show that
		\begin{equation}
		\E (S_{n,T}) \underset{n,T \to +\infty}{\longrightarrow} \E (g(N))
		\end{equation}
	where $N \sim \mathcal{N}(0,V)$. If it was not the case, in the same spirit as in the proof of Theorem \ref{theo_cons_panel}, we could build for some fixed $\varepsilon >0$ two increasing subsequences of integers $(n_k)_{k \geqslant 1}$ and $(T_k)_{k \geqslant 1}$ such that
		\begin{equation}
		\forall k \geqslant 1, \ \left| \E \left( g(S_{n_k,T_k}) \right) - \E \left( g(N) \right) \right| > \varepsilon,
		\end{equation}
	which is in contradiction with \eqref{eq_twosub}.
	
	We thus consider two increasing subsequences of integers $(n_k)_{k \geqslant 1}$ and $(T_k)_{k \geqslant 1}$ and denote for some $x\in \R^{d_{\Theta}}$
		\begin{equation}
		S_{n_k,T_k} = \sum_{t=1}^{T_k} \underbrace{\sum_{j=1}^{n_k} \dfrac{x'\nabla m_{\theta_0}(Z_{j,t})}{\sqrt{n_k T_k}}}_{V_{k,t}},
		\end{equation}
	so, using the Cram\'er-Wold device, we have to prove that
		\begin{equation}
		\sum_{t=1}^{T_k} V_{k,t} \overset{\mathcal{L}}{\longrightarrow} \mathcal{N}(0,x'Vx)
		\end{equation}
	when $k \to +\infty$. To this end, we will apply Theorem \ref{theo_ibragimov_billingsley} to the process $(W_{k,t})_{k,t \in \N^*}$ given for any $k, t \in \N^*$ by
		\begin{equation}
		W_{k,t} = \left\lbrace \begin{array}{l}
		V_{k,t} \quad \text{if} \ t \leqslant T_k \\
		0 \quad \text{if} \ t > T_k
		\end{array}\right. .
		\end{equation}
	The sum $\sum_{t=1}^{+\infty} W_{k,t}$ is necessarily finite for any $k \in \N^*$, and for each couple $(k,t), \ W_{k,t}$ is in $L^2$, since it is the case for $V_{k,t}$ according to assumption \textbf{\emph{P10}}. Thus we can define for any $k,t \in \N^*$
		\begin{equation}
		\sigma_{k,t}^2  = \E \left(W_{k,t}^2\ \mid \ \mathcal{H}_{t-1} \right)
		\end{equation}
	where $\mathcal{H}_{t-1} = \vee_{j=1}^{+\infty} \mathcal{F}_{j,t-1}$. The family of $\sigma$-fields $\left(  \mathcal{H}_{t}\right)_{t \in \N}$ is a filtration, and the sum $\sum_{t=1}^{+ \infty} \sigma_{k,t}^2$ is also finite.

		\begin{enumerate}[labelindent=0pt,labelwidth=!,wide,label=$\triangleright$]
		\item From Assumption \textbf{\emph{P9}}, the random sequence $\left( W_{k,t} \right)_{t \in \N^*}$ is a martingale difference series adapted to the filtration $\left( \mathcal{H}_t \right)_{t \in \N^*}$. 
		\item Let us then show that
			\begin{equation}
			\sum_{t=1}^{+\infty} \sigma_{k,t}^2 \overset{\P}{\longrightarrow} V
			\end{equation}
		when $k \to \infty$. We actually have
			\begin{eqnarray}
			\sum_{t=1}^{\infty} \sigma_{k,t}^2 & = & \sum_{t=1}^{T_k} \E \left( W_{k,t}^2 \ \mid \mathcal{H}_{t-1} \right) \nonumber \\
			& = & \sum_{t=1}^{T_k} \sum_{j=1}^{n_k} \dfrac{1}{n_k T_k} \E \left( \left\vert x'\nabla m_{\theta_0}(Z_{j,t})\right\vert^2 \ \mid \mathcal{H}_{t-1} \right),
			\end{eqnarray}
since for $j\neq \ell$, $Z_{j,t}$ and $Z_{\ell,t}$ are conditionally independent given $\mathcal{H}_{t-1}$ (they depend on $\varepsilon_{j,t}$ and $\varepsilon_{\ell,t}$ and of some variables $\mathcal{H}_{t-1}-$measurable). Moreover setting $\kappa=\left\vert x'\nabla m_{\theta_0}\right\vert^2$, we have
$$\widetilde{Z}_{j,t}:=\E\left(\kappa(Z_{j,t})\vert \mathcal{H}_{t-1}\right)=\int \kappa\circ g\left(X_{j,t-1},u,(\zeta_{t-s})_{s\geq 1}\right)\P_{\varepsilon_{1,1}}(du).$$
One can then observe that $\left(\widetilde{Z}_{j,t}\right)_{j\in \N^{*},t\in \Z}$
is still of the form \eqref{bernshift} with either {\bf C1}
or {\bf C'1} being satisfied. Using the integrability of order $2+\delta$ of the gradient, one can then use a law of large numbers and we get
$$\lim_{k\rightarrow \infty}\frac{1}{n_k T_k}\sum_{j=1}^{n_k}\sum_{t=1}^{T_k}\widetilde{Z}_{j,t}=\E\left(\widetilde{Z}_{1,1}\right)=x'Vx\mbox{ a.s.}$$

		\item We then prove that if $\varepsilon>0$,
			\begin{equation}
			\sum_{t=1}^{+\infty} \E \left(W_{k,t}^2 \ind_{\vert W_{k,t} \vert \geqslant \varepsilon} \right) \underset{k \to +\infty}{\longrightarrow} 0
			\end{equation}
		We have for any $k \in \N$
			\begin{align}
			\sum_{t=1}^{+\infty} \E \left( \vert W_{k,t}\vert^2 \ind_{\vert W_{k,t} \vert \geqslant \varepsilon} \right) & = \sum_{t=1}^{T_k} \E \left( \vert V_{k,t} \vert^2 \ind_{\vert V_{k,t} \vert \geqslant \varepsilon} \right) \nonumber \\
			& =\sum_{t=1}^{T_k} \dfrac{1}{T_k} \E \left( \left\vert\dfrac{1}{\sqrt{n_k}} \sum_{j=1}^{n_k} x'\nabla m_{\theta_0}(Z_{j,t})\right\vert^2 \ind_{ \left\lbrace \left\vert\dfrac{1}{\sqrt{n_k}} \sum_{j=1}^{n_k} x'\nabla m_{\theta_0}(Z_{j,t})\right\vert \geqslant \sqrt{T_k} \varepsilon \right\rbrace} \right) \nonumber \\
			& = \E \left( \left\vert U_k \right\vert^2 \ind_{\left\vert U_k \right\vert^2 \geqslant T_k \varepsilon^2} \right)
			\end{align}
		denoting $U_k = \dfrac{1}{\sqrt{n_k}} \sum_{j=1}^{n_k} \nabla m_{\theta_0}(Z_{j,1})$. This last equality is indeed licit since the distribution of $(Z_{1,t},\ldots,Z_{N_k,t})$ is identical to to the distribution of $(Z_{1,1},\ldots,Z_{N_k,1})$ by stationarity.
It is then enough the show that the sequence $(U_k^2)_{k\geq 1}$ is uniformly integrable. Since $\E U_k^2=x' Vx$, the result will follow from Theorem $3.6$  \citet{billingsley2013convergence} if we show that $U_k$ convergences weakly to a random variable $U$ for which $\E(U^2)=x'Vx$. If {\bf C1} is verified, one can use the central limit theorem for i.i.d. sequences to conclude.
If now {\bf C'1} is satisfied, the result is less obvious.
To get the required convergence, we note that there exists a measurable mapping 
$$\widetilde{g}:\R^k\times \R^{d_2}\times (\R^{d_6}\times\R^k)^{\N}\times (\R^{d_4})^{\N}\rightarrow \R^{d'}$$
such that, setting $\overline{\kappa}_{j,t}=(\varepsilon_{j,t},\kappa_{j,t})$, 
$$Z_{j,1}=\widetilde{g}\left(\varepsilon_{j,1},V_j,(\overline{\kappa}_{j,1-s})_{s\geq 1},(\eta_{1-s})_{s\geq 1}\right).$$
and, using the martingale difference property,
$$\int x'\nabla m_{\theta_0}\circ \widetilde{g}\left(u,V_j,(\overline{\kappa}_{j,1-s})_{s\geq 1},(\eta_{1-s})_{s\geq 1}\right)\P_{\varepsilon_{1,1}}(du)=0$$
almost surely. From the independence assumptions, we deduce that $\E\left(x'\nabla m_{\theta_0}(Z_{j,1})\vert (\eta_s)_{s\geq 0}\right)=0$ a.s.
and that, conditionally on $(\eta_{1-s})_{s\geq 1}$, the variables $x'\nabla m_{\theta_0}(Z_{j,t})$, $j\geq 1$, are i.i.d., centered and the sequence $(U_k)_{k\geq 0}$ converges weakly to a zero-mean Gaussian random variables 
and variance
$$V\left((\eta_{1-s})_{s\geq 1}\right):=\E\left(\left\vert x'\nabla m_{\theta_0}(Z_{1,1})\right\vert^2\vert (\eta_{1-s})_{s\geq 1}\right).$$
We then deduce that the sequence $(U_k)_{k\geq 0}$ converges weakly to a mixture of Gaussian distributions, that is of the form $\sqrt{V\left((\eta_{1-s})_{s\geq 1}\right)}W$ where $W$ follows a standard Gaussian distribution and is independent of the sequence $(\eta_{1-s})_{s\geq 1}$. We then deduce the uniform integrability of the sequence $(U_k)_{k\geq 0}$.
\end{enumerate}
	\end{lemproof}

Let us denote $d_{\Theta}= \dim (\Theta)$. For any couple $(j,t) \in \N^* \times \Z$, and each index $k \in \left\lbrace 1, \ldots ,d_{\Theta} \right\rbrace$, we apply the Taylor-Lagrange formula to $\dfrac{\partial m_{\theta}}{\partial \theta_k}(Z_{j,t})$ between $\hat{\theta}_{n,T}$ and $\theta_0$. There exists $\tilde{\theta}_k$ between $\hat{\theta}_{n,T}$ and $\theta_0$ such that
	\begin{equation}
	\sum_{t=1}^T \sum_{j=1}^n \dfrac{\partial m_{\hat{\theta}_{n,T}}}{\partial \theta_k}(Z_{j,t}) = \sum_{t=1}^T \sum_{j=1}^n \dfrac{\partial m_{\theta_0}}{\partial \theta_k}(Z_{j,t}) + \sum_{t=1}^T \sum_{j=1}^n \sum_{l=1}^{d_{\Theta}} \dfrac{\partial^2 m_{\tilde{\theta}_k}}{\partial \theta_l \partial \theta_k}(Z_{j,t}) \cdot \left( \hat{\theta}_{n,T} - \theta_0 \right)_l.
	\end{equation}
Therefore, if we denote $\nabla^2 m_{\tilde{\theta}}(Z_{j,t})$ the matrix
	\begin{equation}
	\nabla^2 m_{\tilde{\theta}}(Z_{j,t}) \left( \dfrac{\partial^2 m_{\tilde{\theta}_k}}{\partial \theta_l \partial \theta_k}(Z_{j,t}) \right)_{1\leqslant l,k \leqslant p},
	\end{equation}
we get the formula
	\begin{equation}
	\dfrac{1}{nT} \sum_{t=1}^T \sum_{j=1}^n \nabla m_{\hat{\theta}_{n,T}}(Z_{j,t}) = \dfrac{1}{nT}\sum_{t=1}^T \sum_{j=1}^n \nabla m_{\theta_0}(Z_{j,t}) + \dfrac{1}{nT} \sum_{t=1}^T \sum_{j=1}^n \nabla^2 m_{\tilde{\theta}}(Z_{j,t}) \cdot \left( \hat{\theta}_{n,T}-\theta_0 \right).
	\end{equation}
This leads to the equation
	\begin{equation}\label{eq_def_Gnt}
	\sqrt{nT} \left( \E \left( \nabla^2 m_{\theta_0}(Z_{1,1}) \right) + d_{n,T} \right) \cdot \left( \hat{\theta}_{n,T}-\theta_0 \right) = \underbrace{\dfrac{1}{\sqrt{nT}} \sum_{t=1}^T \sum_{j=1}^n \left\lbrace \nabla m_{\hat{\theta}_{n,T}}(Z_{j,t})-\nabla m_{\theta_0}(Z_{j,t}) \right\rbrace}_{G_{n,T}},
	\end{equation}
where
	\begin{equation}
	d_{n,T} = \dfrac{1}{nT} \sum_{t=1}^T \sum_{j=1}^n \nabla^2 m_{\tilde{\theta}}(Z_{j,t}) - \E \left( \nabla^2 m_{\theta_0}(Z_{1,1}) \right).
	\end{equation}
	
	\begin{enumerate}[labelindent=0pt,labelwidth=!,wide,label=$\triangleright$]
	\item We first prove that $d_{n,T} \overset{\P}{\longrightarrow} 0$ when $n,T \to \infty$. We have indeed
		\begin{equation}
		d_{n,T} = \dfrac{1}{nT} \sum_{t=1}^T \sum_{j=1}^n \left\lbrace \nabla^2 m_{\tilde{\theta}}(Z_{j,t}) - \E \left( \nabla^2 m_{\tilde{\theta}}(Z_{1,1}) \right) \right\rbrace + \E \left( \nabla^2 m_{\tilde{\theta}}(Z_{1,1}) - \nabla^2 m_{\theta_0}(Z_{1,1}) \right).
		\end{equation}
	Using Assumptions \textbf{\emph{P5}} and \textbf{\emph{P7}}, it is possible to apply Lemma \ref{lem_panel_unif_law_large_numbers}, and we obtain
		\begin{equation}
		\underset{\theta \in \Theta}{\sup} \left| \dfrac{1}{nT}\sum_{t=1}^T \sum_{j=1}^n \nabla^2 m_{\theta}(Z_{j,t}) - \E \left( \nabla^2 m_{\theta}(Z_{1,1}) \right) \right| \underset{n,T \to \infty}{\longrightarrow} 0 \quad \mathrm{a.s.}
		\end{equation}
	Thus
		\begin{equation}
		\dfrac{1}{nT} \sum_{t=1}^T \sum_{j=1}^n \left\lbrace \nabla^2 m_{\tilde{\theta}}(Z_{j,t}) - \E \left( \nabla^2 m_{\tilde{\theta}}(Z_{1,1}) \right) \right\rbrace \underset{n,T \to +\infty}{\longrightarrow} 0 \quad \mathrm{a.s.}
		\end{equation}
	Furthermore, since $\hat{\theta}_{n,T} \underset{n,T \to \infty}{\longrightarrow} \theta_0$ almost surely, we also have $\tilde{\theta} \underset{n,T \to \infty}{\longrightarrow} \theta_0$ almost surely, and Assumption \textbf{\emph{P7}} ensure that
		\begin{equation}
		\E \left( \nabla^2 m_{\tilde{\theta}}(Z_{1,1}) - \nabla^2 m_{\theta_0}(Z_{1,1}) \right) \underset{n,T \to \infty}{\longrightarrow} 0
		\end{equation}
	by Lebesgue's dominated convergence theorem. Here, we actually prove this latter convergence for any increasing sequences of integers $(n_k)_{k \in \N}$ and $(T_k)_{k \in \N}$. Hence, $d_{n,T} \underset{n,T \to \infty}{\longrightarrow} 0$ almost surely, and in particular $d_{n,T} = o_{\P}(1)$ when $n,T \to +\infty$.
	\item We then show that $G_{n,T} \overset{\mathcal{L}}{\longrightarrow} \mathcal{N} \left(0,\mathrm{Var}\left( \nabla m_{\theta_0}(Z_{1,1}) \right) \right)$ when $n,T \to \infty$.
	
	Since $\hat{\theta}_{n,T}$ is a minimizer of $\theta \longmapsto \sum_{t=1}^T \sum_{j=1}^n m_{\theta}(Z_{j,t})$, and since $\hat{\theta}_{n,T}$ converges almost surely to $\theta_0$ when $n,T \to \infty$, we almost surely have $\hat{\theta}_{n,T} \in \mathring{\Theta}$ for $n$ and $T$ large enough, and thus
		\begin{equation}
		\sum_{t=1}^T \sum_{j=1}^n \nabla m_{\hat{\theta}}(Z_{j,t}) = 0,
		\end{equation}
	hence
		\begin{equation}\label{eq_conv_nabla_m}
		\dfrac{1}{\sqrt{nT}} \sum_{t=1}^T \sum_{j=1}^n \nabla m_{\hat{\theta}_{N,T}}(Z_{j,t}) \underset{n,T \to \infty}{\longrightarrow} 0 \quad \mathrm{a.s.} 
		\end{equation}
	Furthermore, we have by Lemma \ref{lem_proof_theo_nor_panel}
		\begin{equation}
		\dfrac{1}{\sqrt{nT}} \sum_{t=1}^T \sum_{j=1}^n \nabla m_{\theta_0}(Z_{j,t}) \overset{\mathcal{L}}{\longrightarrow} \mathcal{N}\left(0, \mathrm{Var}\left( \nabla m_{\theta_0}(Z_{1,1}) \right) \right)
		\end{equation}
	when $n,T \to \infty$. We thus obtain by Slutsky's lemma
		\begin{equation}
		G_{n,T} \overset{\mathcal{L}}{\longrightarrow}  \mathcal{N}\left(0, \mathrm{Var}\left( \nabla m_{\theta_0}(Z_{1,1}) \right) \right)
		\end{equation}
	when $n,T \to \infty$. Note that once again, this latter convergence is obtained for all increasing subsequences of integers $(N_k)_{k \in \N}$ and $(T_k)_{k \in \N}$.
	\item Let us finally show that $\sqrt{nT} \left( \hat{\theta}-\theta_0 \right) = O_{\P}(1)$ when $n,T \to \infty$ and take any $\varepsilon >0$.
	
	By continuity of the determinant, there exists $\eta >0$ such that for any matrix $x$ satisfying $\| x \| \leqslant \eta$, $\E \left( \nabla^2 m_{\theta_0}(Z_{1,1}) \right) + x$ is invertible. Let us denote
		\begin{equation}
		S_{\eta} : = \underset{\|x\| \leqslant \eta}{\sup} \ \| \left( \E \left( \nabla^2 m_{\theta_0}(Z_{1,1}) \right) +x \right)^{-1} \|,
		\end{equation}
	which is finite by continuity of the inverse. We have for any $M >0$ and $n,T \in \N^*$
		\begin{multline}
		\P \left( \sqrt{nT} \| \hat{\theta}_{n,T} - \theta_0 \| > M \right) = \underbrace{\P \left( \sqrt{nT} \| \hat{\theta}_{n,T} - \theta_0 \| > M , \| d_{n,T} \| \leqslant \eta \right)}_{A_{n,T}(M)} \\
		 + \underbrace{\P \left( \sqrt{nT} \| \hat{\theta}_{n,T} - \theta_0 \| > M, \| d_{n,T} \| > \eta \right)}_{B_{n,T}(M)}.
		\end{multline}
	On the one hand
		\begin{equation}
		A_{n,T}(M) \leqslant \P \left( S_{\eta} \times \| G_{n,T}\| > M , \|d_{n,T}\| \leqslant \eta \right) \leqslant \P \left( \| G_{n,T} \| >  \dfrac{M}{S_{\eta}} \right).
		\end{equation}
	Indeed, if $d_{n,T} \leqslant \eta$, we have
		\begin{align}
		\sqrt{nT} \| \hat{\theta}_{n,T} - \theta_0 \| > M & \iff \| \left( \E \left( \nabla^2 m_{\theta_0}(Z_{1,1}) \right) +d_{n,T} \right)^{-1} \cdot G_{n,T} \| > M \nonumber \\
		& \Longrightarrow \| \left( \E \left( \nabla^2 m_{\theta_0}(Z_{1,1}) \right) +d_{n,T} \right)^{-1} \| \times \| G_{n,T} \| > M \nonumber \\
		& \Longrightarrow S_{\eta} \times \| G_{n,T} \| > M.
		\end{align}
	Since $G_{n,T}$ converges in distribution, $G_{n,T} = O_{\P}(1)$ and there exists $M_0 >0$ and $n_A,T_A \in \N^*$ such that for all $n \geqslant n_A$ and all $T \geqslant T_A$
		\begin{equation}
		A_{n,T}(M_0) \leqslant \P \left( G_{n,T} \geqslant \dfrac{M_0}{S_{\eta}} \right) < \dfrac{\varepsilon}{2}.
		\end{equation}
	On the other hand, $d_{n,T} = o_{\P}(1)$, thus there exists $n_B, T_B \in \N^*$ such that for all $n \geqslant n_B$ and $T \geqslant T_B$
		\begin{equation}
		B_{n,T}(M_0) \leqslant \P ( \|d_{n,T} \| > \eta ) \leqslant \dfrac{\varepsilon}{2}.
		\end{equation}
	Hence, for $n \geqslant \max (n_A,n_B)$ and $T \geqslant \max (T_A,T_B)$ we have
		\begin{equation}
		\P \left( \sqrt{nT} \| \hat{\theta}_{n,T} - \theta_0 \| > M \right) < \varepsilon,
		\end{equation}
	i.e. $\sqrt{nT} \left( \hat{\theta}_{n,T} - \theta_0 \right) = O_{\P}(1)$ when $n,T \to \infty$.
	\end{enumerate}
We deduce from this last point that $d_{n,T} \cdot \sqrt{nT} \left( \hat{\theta}_{n,T} -\theta_0 \right) = o_{\P}(1)$ when $n,T \to + \infty$. Thus, from equations \eqref{eq_def_Gnt} and \eqref{eq_conv_nabla_m}, we have
	\begin{equation}
	\E \left( \nabla^2 m_{\theta_0}(Z_{1,1}) \right) \cdot \sqrt{nT} \left( \hat{\theta}_{n,T} -\theta_0 \right) = -\dfrac{1}{\sqrt{nT}} \sum_{t=1}^T \sum_{j=1}^n \nabla m_{\theta_0}(Z_{j,t}) +o_{\P}(1)
	\end{equation}
when $n,T \to \infty$.

The asymptotic normality then follows from Lemma \ref{lem_proof_theo_nor_panel} and Slutsky's lemma.

\newpage

\printbibliography

\end{document}